\documentclass[11pt,letterpaper]{article}
\usepackage[T1]{fontenc}
\usepackage[margin=1in]{geometry}
\usepackage{amsmath,amssymb,amsthm}
\usepackage{newtxtext,newtxmath}
\usepackage{microtype}
\usepackage{graphicx}
\usepackage{booktabs}
\usepackage{enumitem}
\usepackage{float}
\usepackage{placeins}
\usepackage[font=small,labelfont=bf]{caption}
\usepackage{tikz}
\usetikzlibrary{arrows.meta,positioning}
\usepackage[numbers,sort&compress]{natbib}
\usepackage{xurl}
\usepackage{hyperref}

\theoremstyle{plain}
\newtheorem{theorem}{Theorem}
\newtheorem{corollary}{Corollary}
\newtheorem{proposition}{Proposition}
\newtheorem{lemma}{Lemma}
\theoremstyle{definition}
\newtheorem{definition}{Definition}
\theoremstyle{remark}
\newtheorem{remark}{Remark}

\newcommand{\ds}{\cdots}
\newcommand{\cF}{\mathcal{F}}
\newcommand{\cG}{\mathcal{G}}
\newcommand{\set}[1]{\left\{{#1}\right\}}
\newcommand{\crbk}[1]{\left( #1 \right)}

\title{Asymptotic Optimality of $c\mu$-Rule for LLM Inference Scheduling under Memory Constraint}
\author{Zhenghao Li\\{\small HKUST}
  \and Shengbo Wang\\{\small University of Southern California}
  \and Nian Si\\{\small HKUST}}
\date{}
\hypersetup{
  pdftitle={Asymptotic Optimality of c-mu-Rule for LLM Inference Scheduling under Memory Constraint},
  pdfauthor={Zhenghao Li; Shengbo Wang; Nian Si}
}

\begin{document}
\maketitle
\begin{abstract}
Large language model (LLM) inference, the process by which a pre-trained model receives an input prompt and generates output tokens sequentially, has become increasingly important with the widespread deployment of LLM-based applications. As inference workloads continue to scale, GPU memory has emerged as a critical bottleneck: the long input prompts and generation sequences common in agentic and reasoning workloads lead to substantial key–value (KV) cache usage. In this paper, we study the resulting online scheduling problem from a queueing perspective. Requests arrive stochastically with heterogeneous prompt and output lengths, undergo prefill and autoregressive decode, and compete for a finite KV-cache capacity. Unlike classical queueing models, the duration of a processing step depends on both prefill computation and KV-cache I/O, so service capacity is jointly determined by computation, memory access, and memory occupancy. We develop a fluid model that reduces the system to a multiclass single-server queue with effective service rates capturing computation, memory access, and memory capacity. This fluid analysis yields a $c\mu$ priority rule, where the effective service rate $\mu$ is jointly determined by computational speed, memory bandwidth, and memory capacity. We prove that this rule minimizes cumulative holding cost in the fluid model and, for a finite workload, also minimizes completion time and maximizes full-run output-token throughput. We then construct a feasible scheduling policy for the original finite system and prove that its scaled dynamics converge to the optimal fluid policy. Finally, experiments on a real LLM serving engine using an NVIDIA H800 show that practical implementations of the proposed priority rule substantially reduce average latency while maintaining high throughput relative to existing scheduling policies. The code is available at \url{https://github.com/lzhqw/llm-cmu-scheduling}.

\end{abstract}

\section{Introduction}
The public debut of ChatGPT in late 2022 marked a turning point in the swift adoption of generative artificial intelligence, driven largely by advancements in large language models (LLMs) \citep{brown2020language,openai2023gpt,singh2025openai}. These systems have rapidly integrated into routine workflows across research, education, industry, and entertainment \citep{chatgpt2023,githubcopilot2023,huang2025orlm}, facilitating billions of user interactions to date \citep{chatgpt_users_2023,gordon_2024}. However, alongside significant gains in human productivity, this explosive growth has introduced critical bottlenecks in serving large-scale LLM inference reliably and efficiently.

LLM inference refers to the process by which a pre-trained LLM receives an input prompt and generates output tokens sequentially. Compared with traditional machine learning models, LLM inference is memory and computation intensive, often involves long generation sequences, and can incur substantial latency and resource costs. The importance of efficient inference is further amplified by the growing adoption of reasoning models and test-time compute, which increase the amount of computation required for each request. Deloitte\footnote{\url{https://www.deloitte.com/us/en/insights/industry/technology/technology-media-and-telecom-predictions/2026/compute-power-ai.html}} estimates that inference will account for roughly two-thirds of total AI compute in 2026, up from about one-third in 2023. Consequently, efficient LLM inference is not merely a problem of computational optimization, but a system-level challenge that requires carefully balancing responsiveness, throughput, and resource costs.

Key-value (KV) caching \citep{pope2023efficiently,Kwon2023PagedAttention} has become a standard technique for accelerating LLM inference, largely due to the extensive use of attention mechanisms \citep{vaswani2017attention}. Roughly speaking, the KV cache stores the previously computed key and value representations of past tokens in GPU memory to avoid their repeated computation when generating each new token. Thus, KV caching trades increased memory usage for reduced computation and achieves substantially accelerated inference. However, this acceleration comes at the cost of substantial GPU memory consumption and intensive memory access, making memory capacity and bandwidth increasingly important. This becomes a severe issue for inference requests with long input contexts or long generated sequences, which are increasingly common in agentic and reasoning workloads, making memory a critical bottleneck in modern LLM inference \citep{nvidia2026rubin}.

In this paper, we study LLM inference with KV-cache memory constraints through an operations lens. We abstract and model the inference scheduling task for LLMs as a multiclass queue. In this queueing system model, inference requests arrive over time with heterogeneous workloads, characterized by their prompt and output lengths, and a GPU server processes them according to prescribed dynamic scheduling decisions. The workload each request brings to the system consists of two sequential phases that must be executed in order: a prefill phase followed by an autoregressive decode phase.  The scheduler dynamically determines which requests are admitted and processed at each GPU time step. Our primary objective is to minimize the request latency, or, equivalently, the cumulative number of unfinished requests in the system.

A crucial feature of LLM serving systems, and hence our queueing model, is the KV-cache constraint on GPU memory capacity and bandwidth. As a request progresses through the prefill and decode stages, its KV cache in GPU memory continues to grow until the request completes, while each processing step must also access the KV states accumulated up to that point. These memory requirements can significantly affect the number of requests that can be processed concurrently and the time required for each processing step. In particular, GPU-step processing time depends on both the amount of prefill computation and the KV-cache I/O  volume, while the set of concurrently processed requests must satisfy a hard memory-capacity constraint.

To develop an efficient dynamic scheduling policy for this LLM inference queue, we consider a fluid asymptotic regime in which the token processing rate and request arrival rate scale proportionally, while the total memory capacity remains fixed.  This choice of regime directly captures the key memory bottleneck in LLM inference systems with KV caching, which shifts the primary computational burden toward memory constraints and I/O costs.

The fluid approximation yields a tractable yet expressive limiting system with a clear scheduling structure. We show that, under the prescribed limiting regime, the memory-constrained LLM inference system can be related to a classical multiclass queueing model with homogeneous servers. However, the effective service rates in the fluid system must account for prefill computation, KV-cache I/O, and the effect of memory capacity constraints on service availability. This motivates us to adapt the \(c\mu\) priority rule \citep{Buyukkoc1985Cmu,VanMieghem1995Generalized,Mandelbaum2004Generalized}, which minimizes cumulative holding cost while maximizing full-run output-token throughput in the fluid model. We establish that the \(c\mu\)-rule remains optimal for this fluid system under the two-phase workload structure and the memory constraint in the LLM serving system. Furthermore, to bridge the fluid control problem and the prelimit LLM inference system, we construct a dynamic scheduling policy inspired by the $c\mu$-rule. We show that the proposed policy converges to the $c\mu$-rule in the fluid limit and that the prelimit systems under this policy converge to the proposed fluid system controlled by the $c\mu$-rule.

Concretely, the paper makes the following main contributions.
\begin{itemize}[leftmargin=*]
    \item \textbf{LLM system modeling.}
    We develop a queueing model of LLM inference that explicitly accounts for KV-cache memory capacity constraints and the associated memory I/O costs, with processing times calibrated using real GPU measurements to capture their dependence on the composition of the scheduled workload.

    \item \textbf{Identifying the asymptotic regime.}
    We identify a fluid regime in which processing and arrival rates increase together while memory capacity remains fixed, preserving the system's memory bottleneck.
    The resulting model reveals a classical single-server scheduling structure with LLM-specific behaviors such as two-phase workload and memory constraints. 

    \item \textbf{Asymptotic optimality of the $c\mu$-rule.} We derive the effective service rate $\mu$, which jointly captures prefill computation, KV-cache I/O, and memory capacity constraints.
    We then prove that the $c\mu$-rule minimizes cumulative holding cost and maximizes full-run throughput for finite workloads in the fluid limit. To bridge the fluid and finite-system models, we construct an implementable $c\mu$-inspired policy and show that the resulting dynamics converge to the optimal fluid limit.

    \item \textbf{Numerical and real-system validation.}
    We evaluate practical implementations of the proposed $c\mu$-rule using semi-synthetic workloads and real workload data on a real LLM serving engine. In both settings, our policy demonstrates substantial reductions in average latency relative to the default scheduler while maintaining high throughput.
\end{itemize}
Together, these results bridge classical queueing theory and priority-based scheduling policies with implementable, effective scheduling rules for modern LLM inference systems.

\subsection{Related Work}

Recent advances in LLM serving systems have enabled finer-grained control over batching and GPU memory management. \citet{Yu2022Orca} developed iteration-level scheduling, allowing the batch to be reconstructed after each model iteration. \citet{Kwon2023PagedAttention} introduced PagedAttention, which stores KV caches in blocks to reduce GPU memory waste and enable multiple requests to share cached data. \citet{Agrawal2024SarathiServe} split long prefill computations into smaller chunks and co-batched them with decode work. Together, these advances provide the infrastructure for fine-grained control over batching and GPU memory management, serving as both the foundation for our model and the practical testbed for our experiments.

To analyze LLM serving systems analytically, several pioneering studies have employed queueing models.  \citet{Dai2026Throughput} characterize stable throughput under a per-batch token-budget constraint and establish the throughput optimality of work-conserving policies, while \citet{Nie2026Stability} develop a fine-grained queueing model that tracks KV-cache growth and establish stability and overload conditions for continuous batching under the restrictive assumption of a constant batch processing time. Furthermore, 
\citet{Lin2026Heterogeneous} study prefill--decode contention in multi-GPU settings and develop a scheduling policy that achieves asymptotically optimal revenue as the number of GPUs grows, without explicitly modeling GPU memory-capacity constraints. Recently, \citet{yao2026multi} study KV-cache reuse in multi-turn LLM conversations under a least-recently-used (LRU) eviction policy and establish a mean-field limit for the cache hit ratio as the conversation arrival rate and memory capacity scale proportionally.

The most relevant literature is  \citet{Ao2026WAIT} and \citet{Ao2026Congestion}. For requests with known output lengths, \citet{Ao2026WAIT} propose WAIT to maximize throughput under a memory constraint. The policy waits for enough requests of the same class to accumulate, then processes them together. However, setting these thresholds requires arrival-rate information, and their extension to time-varying arrivals imposes additional load and class-mix conditions. Moreover, obtaining valid parameters for WAIT requires certain load and memory conditions to be satisfied, which may not hold in practical settings. We find that, although theoretically optimal, WAIT's performance can deteriorate in real-GPU implementations due to unnecessary idling. This work models batch processing time only through KV-cache I/O, without considering the extra computation time required to process tokens during the prefill phase.  \citet{Ao2026Congestion} further show that admitting requests whenever memory is available can cause recurring evictions and throughput loss, and mitigate this via rate-limited admission and request mixing, especially in multi-GPU deployments.

Besides queueing-based modeling, other studies approach this problem from an online linear programming perspective: \citet{Jaillet2025Online} develop a hindsight benchmark and the MC-SF algorithm, while \citet{Kong2026Geometry} propose the Smallest Volume First policy based on cumulative KV-cache usage. For offline scheduling, prior work develops algorithms and evaluates their total completion time against a clairvoyant benchmark. Specifically, \citet{Wang2025Variable} study the setting with known output lengths, \citet{Chen2025Adaptive} consider interval predictions of output lengths, and \citet{Feng2026General,Feng2026Competitive} develop algorithms that do not rely on output-length predictions.

LLM inference under KV-cache memory constraints is closely related to stochastic bin packing, where randomly arriving requests must be served subject to resource capacity constraints. Early studies establish stability conditions for First-Fit and Best-Fit scheduling \citep{coffman2001bandwidth,gamarnik2004stochastic}. In infinite-server models, greedy policies with safety stocks \citep{stolyar2013large}, greedy randomized placement \citep{stolyar2015asymptotic}, and randomized Best-Fit \citep{ghaderi2014asymptotic} asymptotically minimize the number of occupied servers. \citet{stolyar2013infinite} establishes asymptotic optimality of greedy allocation for a convex approximation to the server-count objective. For a fixed number of servers, \citet{maguluri2012stochastic} show that Best-Fit need not be throughput optimal and develop MaxWeight-based scheduling policies with throughput guarantees. \citet{gupta2015lagrangian} develop Lagrangian-based packing algorithms for permanent and departing requests, including a model with configuration-dependent service slowdowns. 

Despite similar constraints and objectives, LLM inference scheduling differs from classical stochastic bin-packing problems in several ways. First, while stochastic bin-packing models assume that each request has a fixed size upon arrival, the KV-cache usage of an LLM request grows as generation progresses. Second, LLM inference involves two distinct phases, prefill and decode, with different memory dynamics. Third, whereas service times in classical stochastic bin-packing models are typically independent of current memory occupancy, processing times in LLM inference depend on memory usage due to KV-cache I/O costs.

The objective of minimizing cumulative holding cost connects our problem to the literature on $c\mu$ priority scheduling. \citet{Buyukkoc1985Cmu} establish the optimality of the $c\mu$ rule for a preemptive single-server queue with linear holding costs, geometric service times, and a general arrival process. \citet{VanMieghem1995Generalized} extends this framework to convex delay costs and establishes the asymptotic optimality of a generalized $c\mu$ rule in heavy traffic. \citet{Mandelbaum2004Generalized} further extend the analysis to heterogeneous servers and prove heavy-traffic asymptotic optimality under complete resource pooling. 
\citet{Zychlinski2023Resources} study queues in which jobs may require multiple servers simultaneously and establish the asymptotic optimality of an idle-avoid \(c\mu/m\) policy in a many-server regime. Their model is related to ours through resource constraints, but differs in three key aspects: in LLM inference, a job's resource usage (KV cache) grows during service, service rates depend  on the types and numbers of jobs processed together, and memory capacity remains \(M=O(1)\), so each job's resource requirement is non-negligible relative to total capacity.

The remainder of the paper is organized as follows. In Section \ref{sec:LLM-system}, we introduce the LLM inference systems. In Section \ref{sec:queueing-scaling}, we
present our queueing model and its fluid control counterpart for LLM inference. 
In Section \ref{sec:fluid-optimal-control}, we establish the optimality of the $c\mu$ rule for the fluid model. In Section \ref{sec:prelimit-convergence}, we construct a feasible policy for the discrete system and prove that its scaled dynamics converge to those under the optimal fluid policy. Section \ref{sec:numerical} presents real GPU experiments that validate our theoretical results. Finally, we conclude the paper and discuss limitations and directions for future work in Section \ref{sec:conclusion}.

\section{LLM Inference Systems}
\label{sec:LLM-system}
We consider an online LLM inference system in which user requests arrive stochastically over time and wait in a queue for service on a single GPU. Each arriving request contains an input prompt and asks the model to generate an output sequence. A scheduler determines which requests are processed by the GPU over time.

\begin{figure}[ht]
\centering
\includegraphics[width=\textwidth]{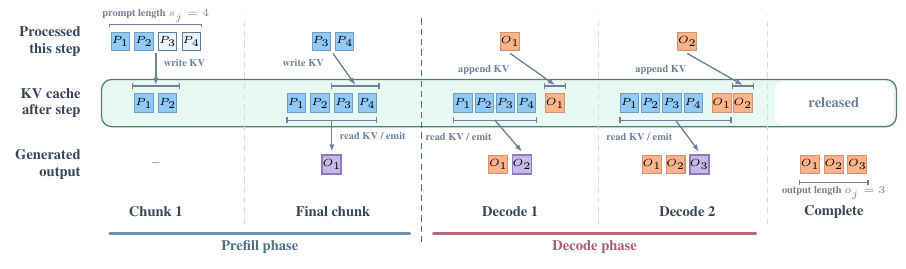}
\caption{Lifecycle of a request with prompt length $s_j=4$, output length $o_j=3$, and two-token prefill chunks. The first prefill step processes $P_1,P_2$ and writes their keys and values to the KV cache without emitting an output token. The final prefill step processes $P_3,P_4$, writes their keys and values to the cache, and emits $O_1$. Each decode step takes the latest output token as input, appends its keys and values to the cache, and emits the next token. After the request completes, its KV cache is released.}
\label{fig:llm-inference-lifecycle}
\end{figure}

The processing of each request consists of a prefill phase followed by an autoregressive decode phase. During the prefill phase, the model processes the prompt tokens and constructs their attention key--value (KV) states. With chunked prefill, the input prompt is divided into several smaller chunks that are processed sequentially and the first output token is generated when the final chunk is processed. During the decode phase, the model generates the remaining output tokens one at a time. Each decode step takes the latest output token as input, appends its KV state, and generates the next output token. The request's KV cache occupies GPU memory throughout execution and is released upon completion. Figure~\ref{fig:llm-inference-lifecycle} summarizes this lifecycle.

Modern inference engines use a GPU step as the smallest unit of scheduling. At each GPU step, a batch is formed from requests that may be at different stages of execution. After the step completes, the scheduler updates request states and may re-form the next batch from eligible waiting and running requests; requests may therefore enter or leave the scheduled batch between GPU steps. This mechanism is commonly called continuous batching~\citep{Yu2022Orca}. Modern engines may also use mixed prefill--decode batching, in which a GPU step processes prefill work from some requests and decode work from others~\citep{Agrawal2024SarathiServe}.

\begin{figure}[t]
\centering
\includegraphics[width=0.8\textwidth]{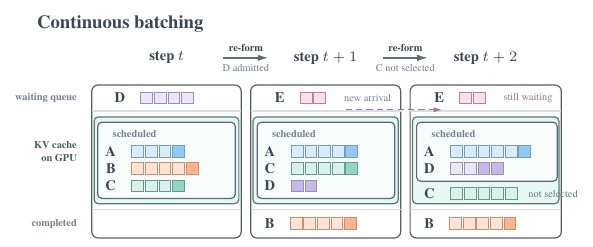}
\caption{Continuous batching across GPU-step boundaries. During a step, the batch remains fixed and the scheduled requests advance by one step. At the next boundary, completed requests leave, and the scheduler may admit waiting requests or leave resident requests unscheduled.}
\label{fig:continuous-batching}
\end{figure}

Every processed token contributes KV states, so a request's KV-cache occupancy grows during prefill and then by one token at each decode step. The running requests share a finite GPU KV-cache capacity, which we measure in tokens. Figure~\ref{fig:continuous-batching} illustrates how the scheduled batch changes at GPU-step boundaries.

Our queueing model represents request progress through its prefill and decode stages. Requests are partitioned into $J$ classes. A class-$j$ request has known prompt length $s_j$ and output length $o_j$. The GPU has KV-cache capacity $M$ tokens. Each prompt is partitioned into chunks of size $c_{\mathrm{chunk}}$. Accordingly, we define the number of prefill chunks of a class-$j$ request by
$
L_j=\left\lceil\frac{s_j}{c_{\mathrm{chunk}}}\right\rceil.
$ 
When a class-$j$ request has $k$ prefill chunks remaining, the next chunk contains
\begin{equation*}
\ell_{j,k}:=\min\{c_{\mathrm{chunk}},s_j-(L_j-k)c_{\mathrm{chunk}}\},\qquad k=1,\ldots,L_j,
\end{equation*}
prompt tokens.
Therefore, a class-$j$ request has $L_j$ sequential prefill stages followed by $o_j-1$ sequential decode stages, with the final prefill stage producing the first output token. If $o_j=1$, the request leaves the system immediately after its final prefill stage.

Let $m(t)$ denote the maximum KV-cache usage during step $t$.  It consists of the KV states already occupying GPU memory and those newly written during the step.  The memory constraint is
\begin{equation*}
m(t)\leq M,\qquad t=0,1,\dots.
\end{equation*}

For the batch executed in step $t$, let $n_{\mathrm{pf}}(t)$ be the total number of prompt tokens in the prefill stage to be processed in this batch, and let $n_{\mathrm{kv}}(t)$ be the batch's total KV-cache I/O volume in token units. In our model, this volume counts the existing KV states accessed by the batch and the new KV states written during the step.  

\begin{figure}[t]
\centering
\begin{minipage}[t]{0.49\textwidth}
\centering
\includegraphics[width=\linewidth]{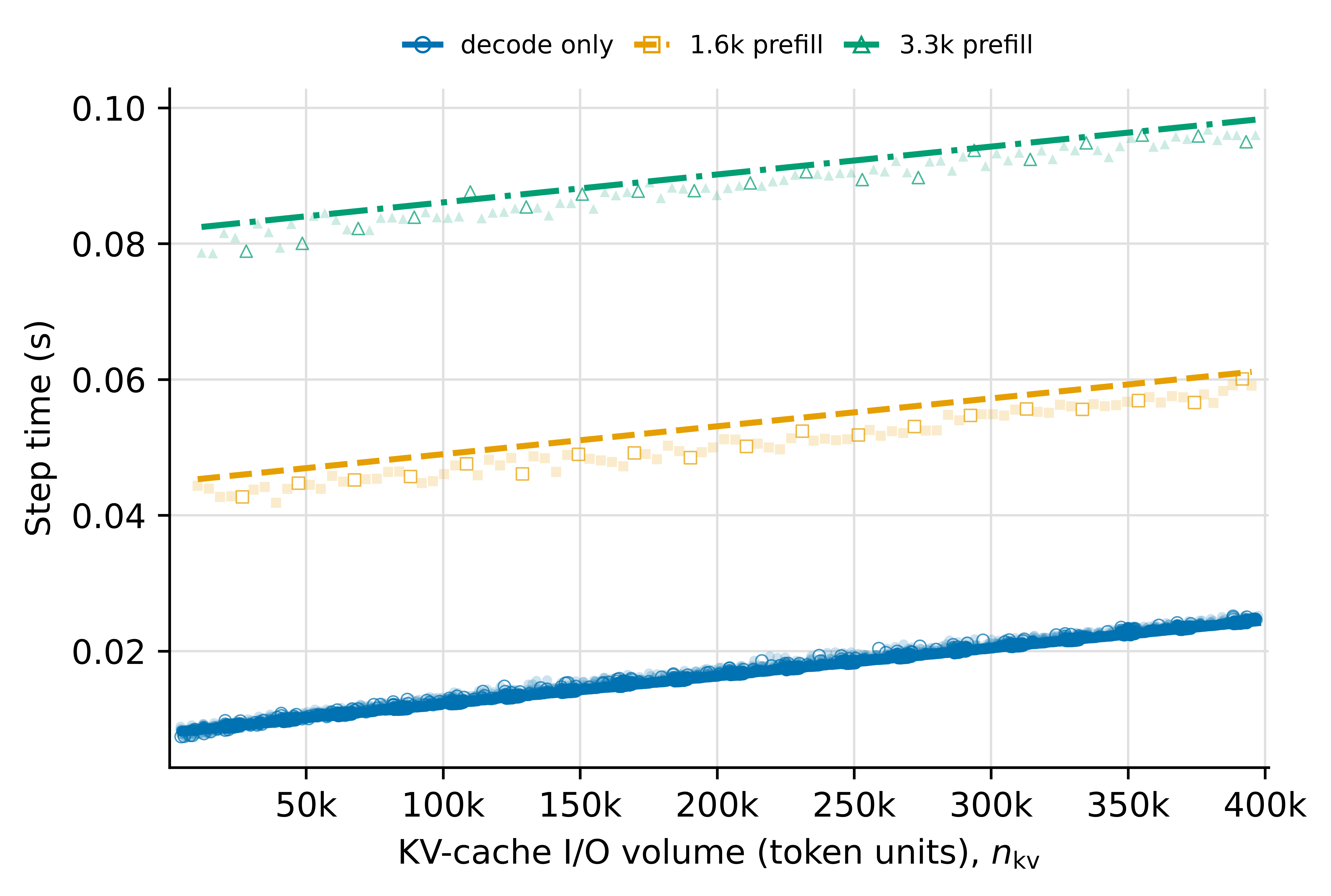}
\end{minipage}\hfill
\begin{minipage}[t]{0.49\textwidth}
\centering
\includegraphics[width=\linewidth]{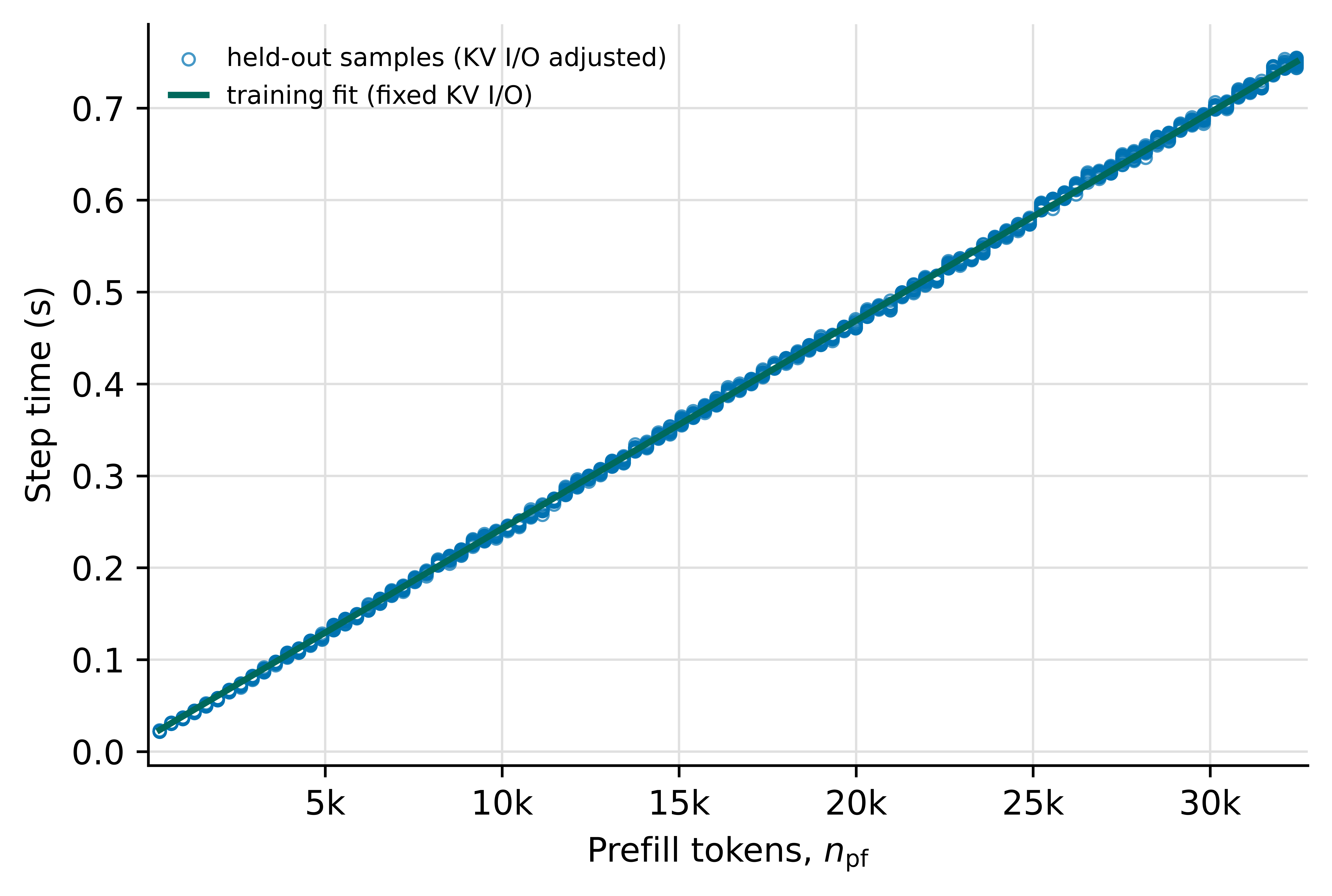}
\end{minipage}
\caption{Conditional fits of the GPU step-time model, estimated on the training sample. Left: KV-cache I/O effect at fixed prefill levels. Right: prefill effect after adjusting held-out samples to fixed KV-cache I/O.}
\label{fig:step-time-validation}
\end{figure}

By empirically testing different mixed batches on an NVIDIA H800, we observe that the step time (batch processing time) varies approximately linearly with both KV-cache I/O volume and the number of prompt tokens processed.
Figure~\ref{fig:step-time-validation} shows these linear relationships.

Motivated by these observations, we fit the step time of step $t$ as
\begin{equation}\label{eq:step-time-model}
\widehat\Delta(t):=\widehat a_{\mathrm{pf}}n_{\mathrm{pf}}(t)+\widehat a_{\mathrm{kv}}n_{\mathrm{kv}}(t)+\widehat b,
\end{equation}
where $\widehat a_{\mathrm{pf}}$ is the fitted per-prompt-token processing-time coefficient, $\widehat a_{\mathrm{kv}}$ is the fitted per-token KV-cache I/O coefficient, and $\widehat b$ is the estimated fixed per-step overhead. Table~\ref{tab:step-time-fit} reports the coefficient estimates and fit statistics. 

\begin{table}[t]
\centering
\caption{Batch processing time estimates in an NVIDIA H800. Out-of-sample statistics use held-out existing-KV strata.}
\label{tab:step-time-fit}
\begin{minipage}[t]{0.36\textwidth}
\vspace{0pt}
\centering
\begin{tabular}{@{}lc@{}}
\toprule
Parameter & Estimate \\
\midrule
$\widehat a_{\mathrm{pf}}$ (s/token) & $2.264\times10^{-5}$ \\
$\widehat a_{\mathrm{kv}}$ (s/token) & $4.108\times10^{-8}$ \\
$\widehat b$ (s) & $7.794\times10^{-3}$ \\
\bottomrule
\end{tabular}
\end{minipage}\hfill
\begin{minipage}[t]{0.52\textwidth}
\vspace{0pt}
\centering
\begin{tabular}{@{}lcc@{}}
\toprule
Metric & In sample & Out of sample \\
\midrule
Observations & $16{,}205$ & $3{,}987$ \\
$R^2$ & $0.999902$ & $0.999905$ \\
RMSE (ms) & $2.294$ & $2.268$ \\
\bottomrule
\end{tabular}
\end{minipage}
\end{table}

\section{System Modeling and Fluid Regime}
\label{sec:queueing-scaling}
In this section, we formulate a sequence of prelimit queueing models for the LLM inference system and derive the corresponding fluid control problem.  We first specify the system state, queueing dynamics, control decisions, and memory constraints. We then apply fluid scaling and formulate the corresponding fluid control problem. 
Finally, we show that the fluid control problem also exhibits a state-space collapse phenomenon commonly observed in the queueing literature \citep{Bramson1998StateSpaceCollapse,Williams1998DiffusionApproximations}.

\subsection{Queueing Model}
Consider a sequence of queueing systems indexed by $r=1,2,\dots$, while the fluid model is obtained as $r\to\infty$. The arrival rate of requests in each class and the GPU processing speed both scale with $r$, while the KV-cache capacity $M$ remains fixed. This regime represents the memory bottleneck created by the KV caching techniques that trade memory capacity and I/O cost for computation speed in modern LLM inference systems. Concretely, for the step-time model in \eqref{eq:step-time-model}, we assume that the positive coefficients of the $r$th system satisfy
\begin{equation}\label{eq:processing_time_limit}
r a_{\mathrm{pf}}^r\to a_{\mathrm{pf}}>0,\quad r a_{\mathrm{kv}}^r\to a_{\mathrm{kv}}>0,\quad r b^r\to b>0.
\end{equation}
Thus, the one-step processing time of fixed batch is of order $1/r$. Accordingly, define the step time
\begin{equation*}
\begin{aligned}
\Delta^r(t)&:=a_{\mathrm{pf}}^r n_{\mathrm{pf}}^r(t)+a_{\mathrm{kv}}^r n_{\mathrm{kv}}^r(t)+b^r,
\end{aligned}
\end{equation*}
where \(n_{\mathrm{pf}}^r(t)\) is the number of prefill tokens processed in step \(t\), and \(n_{\mathrm{kv}}^r(t)\) is the batch's KV-cache I/O volume in token units, as formally defined in \eqref{eq:def:npfkv}.

Throughout the paper, we index GPU processing steps by the integer $t = 0, 1, 2, \dots$ and denote clock time by $\tau \geq 0$, since processing steps may have different durations. Thus, the clock time at the beginning of step \(t\) is given by
$$
f^r(t):=\sum_{u=0}^{t-1}\Delta^r(u). 
$$

At the beginning of GPU step $t$, let $p_{j,k}^{r}(t)$ be the number of class-$j$ requests with $k$ prefill chunks waiting to be processed and let $d_{j,k}^{r}(t)$ be the number of class-$j$ requests with $k$ decode tokens remaining. The system state records the number of requests at each prefill and decode stage:
\begin{equation*}
z_t^r:=\bigl(\{p_{j,k}^r(t)\}_{1\leq j\leq J,\,1\leq k\leq L_j},\{d_{j,k}^r(t)\}_{1\leq j\leq J,\,1\leq k\leq o_j-1}\bigr),
\end{equation*}
where \(p_{j,L_j}^{r}(0)\) represents the class-\(j\) requests already present in the system at time zero.

At each GPU step $t$, a policy $\pi^r$ determines how many requests to process at each prefill and decode stage. Specifically, the policy schedules one prefill chunk for $\Delta\mathcal U_{j,k}^r(t)$ class-$j$ requests with $k$ prefill chunks remaining, and one decode step for $\Delta\mathcal V_{j,k}^r(t)$ class-$j$ requests with $k$ output tokens remaining. Since in practice, each request usually consumes only a small fraction of the total memory capacity \(M\) (at most $2\%$ in our datasets), we allow fractional requests and treat queue lengths and per-step controls as nonnegative real-valued quantities. The number of requests processed at each step cannot exceed the corresponding queue length:
\[
0\leq\Delta\mathcal U_{j,k}^r(t)\leq p_{j,k}^r(t)
\quad \text{and}\quad
0\leq\Delta\mathcal V_{j,k}^r(t)\leq d_{j,k}^r(t).
\]
To facilitate the fluid scaling and the analysis of the limiting dynamics in Section \ref{subsec:fluid}, we introduce the following cumulative stage-completion processes:
\[
\mathcal U_{j,k}^r(t):=\sum_{s=0}^{t-1}\Delta\mathcal U_{j,k}^r(s)\quad \text{and}\quad
\mathcal V_{j,k}^r(t):=\sum_{s=0}^{t-1}\Delta\mathcal V_{j,k}^r(s),
\]
with $\mathcal{U}_{j,k}^r(0)=0$ and $\mathcal{V}_{j,k}^r(0)=0$. Since the increments are recovered by taking one-step differences, the cumulative processes provide an equivalent representation of the controls. We therefore write \(\pi^r=(\mathcal U^r,\mathcal V^r)\).

The selected batch must also satisfy the KV-cache capacity constraint.
The maximum KV-cache usage during step $t$ is
\begin{equation*}
\begin{aligned}
m^r(t)&=\sum_{j=1}^J\Biggl(\sum_{k=1}^{L_j}c_{\mathrm{chunk}}(L_j-k)p_{j,k}^r(t)+\sum_{k=1}^{o_j-1}(s_j+o_j-k-1)d_{j,k}^r(t)\\
&\hphantom{=\sum_{j=1}^J\Biggl(}{}+\sum_{k=1}^{L_j}\ell_{j,k}\Delta\mathcal U_{j,k}^r(t)+\sum_{k=1}^{o_j-1}\Delta\mathcal V_{j,k}^r(t)\Biggr).
\end{aligned}
\end{equation*}
The first two sums count the KV states already in the memory before step $t$ begins, while the last two count the new KV states for which the selected batch reserves slots. Thus, the memory usage should satisfy the constraint $m^r(t)\leq M$ at all time steps.

For the selected batch, let $n_{\mathrm{pf}}^r(t)$ be the number of prefill tokens processed in step $t$, and let $n_{\mathrm{kv}}^r(t)$ be its KV-cache I/O volume in token units:
\begin{equation}
\begin{aligned}
n_{\mathrm{pf}}^r(t)&=\sum_{j=1}^J\sum_{k=1}^{L_j}\ell_{j,k}\Delta\mathcal U_{j,k}^r(t),\\
n_{\mathrm{kv}}^r(t)&=\sum_{j=1}^J\Biggl(\sum_{k=1}^{L_j}[c_{\mathrm{chunk}}(L_j-k)+\ell_{j,k}]\Delta\mathcal U_{j,k}^r(t)+\sum_{k=1}^{o_j-1}(s_j+o_j-k)\Delta\mathcal V_{j,k}^r(t)\Biggr).
\end{aligned}
\label{eq:def:npfkv}
\end{equation}
Thus, each selected prefill stage contributes its existing KV states and the states written for its current chunk, while each selected decode stage contributes its existing KV states and one new write. 

The arrival processes are modeled as nonhomogeneous Poisson processes, with arrival rates scaled by a factor of $r$ to match the processing-time scaling. Specifically, let $\lambda_j^r(\tau)=r\lambda_j(\tau)$ denote the class-$j$ arrival rate at clock time $\tau$ of the $r$th system. For each $r\ge 1$, let $N_j^r$, $j=1,\dots,J$, be independent unit-rate Poisson processes. The number of class-$j$ requests arriving during step $t$ is represented as 
\begin{equation}
\eta_j^r(t):=N_j^r\!\left(r\int_{0}^{f^r(t+1)}\lambda_j(s)\,ds\right) - N_j^r\!\left(r\int_{0}^{f^r(t)}\lambda_j(s)\,ds\right)
\label{eq:prelimit-arrival-increments}
\end{equation}
Therefore, including the class-\(j\) requests already present in the system at time zero, the cumulative class-\(j\) arrival process is
\begin{equation}
A_j^r(t):=p_{j,L_j}^r(0)+\sum_{u=0}^{t-1}\eta_j^r(u) = p_{j,L_j}^r(0) + N_j^r\!\left( r\int_0^{f^r(t)}\lambda_j(s)\,ds \right),\qquad t\in\mathbb Z_+.
\label{eq:A-def}
\end{equation}

These Poisson arrivals are the only sources of randomness in the prelimit systems. Let $\cF_\tau^r$ be the information contained in the arrival process before clock time $\tau$; i.e.
\[
\cF_\tau^r = \sigma\crbk{N_j^r\crbk{r\int_0^u\lambda_j(s)ds}:j=1,\ds,J; u\leq \tau}
\]
Under the admissibility condition in Definition \ref{def:admissible-prelimit-policy}, we will establish that $\set{f^r(t):t\ge 0}$ is a sequence of $\cF_\tau^r$ stopping times. Thus, we can define the history at each decision time as $\cG^r_t:=\cF^r_{f^r(t)}$.

Requests arriving during step $t$ join the waiting queue at step $t+1$ and can be scheduled from that step onward. The queueing dynamics are
\begin{equation}
{\medmuskip=1mu\thickmuskip=2mu \begin{aligned}
p_{j,L_j}^{r}(t+1)&=p_{j,L_j}^{r}(t)+\eta_j^r(t)-\Delta\mathcal U_{j,L_j}^{r}(t),\\
p_{j,k}^{r}(t+1)&=p_{j,k}^{r}(t)+\Delta\mathcal U_{j,k+1}^{r}(t)-\Delta\mathcal U_{j,k}^{r}(t), & k=1,\dots,L_j-1;\\
d_{j,o_j-1}^{r}(t+1)&=d_{j,o_j-1}^{r}(t)+\Delta\mathcal U_{j,1}^{r}(t)-\Delta\mathcal V_{j,o_j-1}^{r}(t),\\
d_{j,k}^{r}(t+1)&=d_{j,k}^{r}(t)+\Delta\mathcal V_{j,k+1}^{r}(t)-\Delta\mathcal V_{j,k}^{r}(t), & k=1,\dots,o_j-2,
\end{aligned}}
\label{eq:prelimit-dynamics}
\end{equation}
when $o_j\geq2$. If $o_j=1$, the request completes immediately after its final prefill stage. The first equation states that the class-$j$ waiting-queue length at step $t+1$ equals its value at step $t$, plus arrivals during step $t$, minus the admission increment $\Delta\mathcal U_{j,L_j}^r(t)$. The remaining three equations have the same inflow-minus-outflow form: completions from the preceding stage enter the state, and completions of the current stage leave it.

We now specify the admissibility requirements for prelimit policies.
\begin{definition}[Admissible prelimit policy]
\label{def:admissible-prelimit-policy}
Fix an admissible initial state $z^r_0$. A prelimit policy $\pi^r=(\mathcal U^r,\mathcal V^r)$ is admissible if the following conditions hold at every GPU step $t$.
\begin{enumerate}
\item $\pi^r$ is non-anticipative; i.e. $\Delta \mathcal U_{j,k}^r(t)$ and $\Delta \mathcal V_{j,k}^r(t)$ are $\cG^r_t$-measureable for all indices $(j,k)$ ranging over  prefill and decode stages.
\item The KV-cache capacity constraint holds: $m^r(t)\leq M$.
\item For every class $j=1,\dots,J$, the control increments satisfy
\begin{equation*}
0\leq\Delta\mathcal U_{j,k}^r(t)\leq p_{j,k}^r(t)\quad\text{and}\quad 0\leq\Delta\mathcal V_{j,k}^r(t)\leq d_{j,k}^r(t).
\end{equation*}
\end{enumerate}
\end{definition}

In Appendix \ref{sec:prelimit-nonanticipativity}, we note that the non-anticipativity condition and $\cG_t^r$ are well-defined by showing that $\set{f^r(t):t\ge 0}$ is a sequence of $\cF_\tau^r$ stopping times.

Let $c_j\geq0$ be the holding-cost rate for each class-$j$ request. The normalized expected cumulative holding cost under a prelimit scheduling policy $\pi^r$ is
\begin{equation}\label{eq:prelimit_obj}
\mathcal J^r(\pi^r;z_0^r):=\frac1r\mathbb E\Biggl[\sum_{t=0}^{\infty}\Delta^r(t) \sum_{j=1}^Jc_j\biggl(\sum_{k=1}^{L_j}p_{j,k}^r(t)+\sum_{k=1}^{o_j-1}d_{j,k}^r(t)\biggr)\,\Biggm|\,z_0^r\Biggr].
\end{equation}
The prelimit control problem is to minimize $\mathcal J^r(\pi^r;z_0^r)$ over all feasible scheduling policies.

We assume that each \(\lambda_j\) is nonnegative and locally bounded, with \(\lambda_j(\tau)=0\) for all \(\tau>H\), where \(H<\infty\) is fixed.
\subsection{Fluid Scaling}
\label{subsec:fluid}
Within an $O(1)$ unit of clock time in the $r$th system, since the arrival rates scale as $O(r)$ and each GPU step takes $O(1/r)$ units of clock time, cumulative arrivals, cumulative stage completions, and waiting-queue lengths can all be of order $O(r)$. This suggests the following scaling, which we use to restore stability and obtain a sensible limiting system. For $\xi\in\mathbb R_+$, define
\begin{equation*}
\begin{aligned}
\bar p_{j,k}^{r}(\xi) &:= \frac{1}{r}p_{j,k}^{r}(\lfloor r\xi\rfloor), & \bar d_{j,k}^{r}(\xi) &:= \frac{1}{r}d_{j,k}^{r}(\lfloor r\xi\rfloor),\\
\bar{\mathcal U}_{j,k}^{r}(\xi) &:= \frac{1}{r}\mathcal U_{j,k}^{r}(\lfloor r\xi\rfloor), & \bar{\mathcal V}_{j,k}^{r}(\xi) &:= \frac{1}{r}\mathcal V_{j,k}^{r}(\lfloor r\xi\rfloor),\\
\bar A_j^r(\xi) &:= \frac{1}{r}A_j^r(\lfloor r\xi\rfloor), & \bar f^r(\xi) &:= f^r(\lfloor r\xi\rfloor).
\end{aligned}
\end{equation*}
Here $\xi=t/r$ is fluid-scaled step time. When limits of these scaled processes exist, we omit the bars and the index $r$ and write the limiting processes as $p$, $d$, $\mathcal U$, $\mathcal V$, and $f$. The formal convergence from the prelimit system to the fluid limit is established in Section~\ref{sec:prelimit-convergence}.

A crucial phenomenon as a consequence of this fluid regime and the memory constraint stems from the following observation. For $k=1,\ldots,L_j-1$, each request counted by $p_{j,k}^r$ occupies $c_{\mathrm{chunk}}(L_j-k)$ KV-cache slots. Likewise, for $k=1,\ldots,o_j-1$, each request counted by $d_{j,k}^r$ occupies $s_j+o_j-k-1$ KV-cache slots. Since the total capacity $M$ is fixed independently of $r$, the memory constraint bounds each of these unscaled states by a constant independent of $r$. Therefore,
\begin{equation}
\bar p_{j,k}^{r}(\xi) \leq\frac{M}{r c_{\mathrm{chunk}}(L_j-k)} \rightarrow0\quad\text{and}\quad \bar d_{j,k}^{r}(\xi) \leq\frac{M}{r(s_j+o_j-k-1)} \rightarrow0.
\label{eq:vanishing-intermediate-states}
\end{equation}
Thus, every fluid-scaled intermediate-stage state vanishes as $r\to\infty$. The fluid-scaled waiting-queue length $\bar p_{j,L_j}^r$ need not vanish because waiting requests do not occupy GPU memory. This observation naturally leads to the phenomenon of stage-level state-space collapse, yielding a more tractable limiting system for which we can characterize an optimal policy. We provide a detailed characterization of this
phenomenon in Section~\ref{subsec:stage-level-collapse}. 
As a result, on every bounded horizon, the normalized holding-cost contribution of requests in intermediate stages vanishes. Thus, although $\mathcal J^r$ counts all requests in the system, its fluid counterpart contains only the waiting-queue lengths $p_{j,L_j}$.

The cumulative stage completions need not vanish with these states. A fixed fluid-time interval $[\xi_1,\xi_2]$ contains $O(r)$ GPU steps. Since each step can complete $O(1)$ requests at a given stage, the cumulative number of stage completions over the interval can be $O(r)$. After division by $r$, the corresponding cumulative flow can therefore remain $O(1)$.

\subsection{Fluid Control Problem}
We now formulate the limiting fluid control problem driven by the arrival rates $\lambda_j$. Fluid evolution is indexed by step time $\xi\in\mathbb{R}_+$, the continuous analogue of the discrete GPU-step index, whereas $\tau\in\mathbb{R}_+$ still represents clock time. Under a policy $\pi$, the real-time clock $f^\pi$ maps step time $\xi$ to clock time $\tau=f^\pi(\xi)$.

At step time $\xi$, $p_{j,k}(\xi)\in\mathbb{R}_+$ is the fluid-scaled number of class-$j$ requests with $k$ prefill chunks remaining, and $d_{j,k}(\xi)\in\mathbb{R}_+$ is the fluid-scaled number of class-$j$ requests with $k$ decode tokens remaining. The system state is,
\begin{equation*}
z_\xi=\bigl(\{p_{j,k}(\xi)\}_{1\leq j\leq J,\,1\leq k\leq L_j},\{d_{j,k}(\xi)\}_{1\leq j\leq J,\,1\leq k\leq o_j-1}\bigr).
\end{equation*}
For the fluid system, \eqref{eq:vanishing-intermediate-states} implies that the memory constraint restricts admissible initial states to
\begin{equation*}
\begin{aligned}
&p_{j,L_j}(0)\geq0,\\
&p_{j,k}(0)=0,\quad k=1,\dots,L_j-1,\\
&d_{j,k}(0)=0,\quad k=1,\dots,o_j-1,
\end{aligned}
\end{equation*}
for every class $j$.

Prefill and decode schedules are represented by the respective cumulative control processes. For $k=1,\dots,L_j$, $\mathcal U_{j,k}(\xi)$ is the fluid-scaled cumulative number of class-$j$ requests that have completed the transition from $k$ to $k-1$ remaining prefill chunks by step time $\xi$. For $k=1,\dots,o_j-1$, $\mathcal V_{j,k}(\xi)$ is defined analogously for decode tokens. 

For cumulative control processes $\pi=(\mathcal U,\mathcal V)$, the real-time clock is
\begin{equation*}
\begin{aligned}
&f^\pi(\xi)=a_{\mathrm{pf}}\sum_{j=1}^J\sum_{k=1}^{L_j}\ell_{j,k}\mathcal U_{j,k}(\xi)\\
&\quad+a_{\mathrm{kv}}\sum_{j=1}^J\Biggl(\sum_{k=1}^{L_j}[c_{\mathrm{chunk}}(L_j-k)+\ell_{j,k}]\mathcal U_{j,k}(\xi)+\sum_{k=1}^{o_j-1}(s_j+o_j-k)\mathcal V_{j,k}(\xi)\Biggr)+b\xi.
\end{aligned}
\end{equation*}
The three clock terms accumulate the clock time due to prefill computation, KV-cache I/O, and the fixed per-step overhead, respectively.

Given this clock, the fluid state is determined by the cumulative form of the stagewise flow-balance equations, with arrivals evaluated up to clock time $f^\pi(\xi)$:
\begin{equation}
{\medmuskip=1mu\thickmuskip=2mu\scriptspace=0pt\thinmuskip=2mu
\begin{aligned}
p_{j,L_j}^{\pi}(\xi)&=p_{j,L_j}(0)+\int_0^{f^\pi(\xi)}\lambda_j(s)\,ds-\mathcal U_{j,L_j}(\xi),\\
p_{j,k}^{\pi}(\xi)&=p_{j,k}(0)+\mathcal U_{j,k+1}(\xi)-\mathcal U_{j,k}(\xi),&k=1,\dots,L_j-1;\\
d_{j,o_j-1}^{\pi}(\xi)&=d_{j,o_j-1}(0)+\mathcal U_{j,1}(\xi)-\mathcal V_{j,o_j-1}(\xi),\\
d_{j,k}^{\pi}(\xi)&=d_{j,k}(0)+\mathcal V_{j,k+1}(\xi)-\mathcal V_{j,k}(\xi), & k=1,\dots,o_j-2.
\end{aligned}}
\label{eq:fluid-stage-flow-balance}
\end{equation}
The equation for $d_{j,o_j-1}^{\pi}$ applies only when $o_j\geq2$. If $o_j=1$, the request completes immediately after its final prefill stage, and there are no decode states.
The fluid control admissibility conditions below inherit the structure of the prelimit model. 
\begin{definition}[Admissible fluid policy]
\label{def:admissible-fluid-policy}
Fix an admissible initial state $z_0$. A policy $\pi=(\mathcal U,\mathcal V)$ is admissible if the following conditions hold:
\begin{enumerate}
\item The cumulative control processes take values in $\mathbb{R}_+$, are continuous and nondecreasing, and start from zero; i.e. $\mathcal U_{j,k}(0)=0$ and $\mathcal V_{j,k}(0)=0$. 
\item For every $\xi_1>\xi_2\geq0$, the cumulative controls satisfy the memory constraint
\begin{equation}
{\medmuskip=1mu\thickmuskip=2mu\thinmuskip=2mu\scriptspace=0pt \begin{aligned}
&\sum_{j=1}^J\!\Biggl(\sum_{k=1}^{L_j}[c_{\mathrm{chunk}}(L_j-k)+\ell_{j,k}]\bigl(\mathcal U_{j,k}(\xi_1)-\mathcal U_{j,k}(\xi_2)\bigr)\\
&\qquad+\sum_{k=1}^{o_j-1}(s_j+o_j-k)\bigl(\mathcal V_{j,k}(\xi_1)-\mathcal V_{j,k}(\xi_2)\bigr)\Biggr)\leq M(\xi_1-\xi_2).
\end{aligned}}
\label{eqn:lips}
\end{equation}
\item For every $\xi\geq0$, the waiting-queue lengths are nonnegative and all intermediate-stage states are zero:
\begin{equation}
\begin{aligned}
p_{j,L_j}^{\pi}(\xi) &\ge 0,\qquad j=1,\dots,J,\\
p_{j,k}^{\pi}(\xi) &= 0,\qquad j=1,\dots,J,\quad k=1,\dots,L_j-1,\\
d_{j,k}^{\pi}(\xi) &= 0,\qquad j=1,\dots,J,\quad k=1,\dots,o_j-1.
\end{aligned}
\label{eq:fluid-state-constraints}
\end{equation}
\end{enumerate}
\end{definition}

Definition \ref{def:admissible-fluid-policy} imposes three sets of constraints that reflect the requirements for prelimit policies in Definition \ref{def:admissible-prelimit-policy}. First, the (deterministic) fluid cumulative control processes \(\mathcal U\) and \(\mathcal V\) must start from zero and be nondecreasing. Second, the per-step memory constraint \(m^r(t)\leq M\) in the prelimit model translates into the Lipschitz constraint \eqref{eqn:lips} in the fluid model. The left-hand side measures the cumulative KV-cache usage associated with work completed during \((\xi_2,\xi_1]\), which cannot exceed \(M(\xi_1-\xi_2)\), the memory capacity integrated over the interval. Third, the nonnegativity of the prelimit queues is preserved in the fluid limit, while the convergence in \eqref{eq:vanishing-intermediate-states} implies that all intermediate-stage fluid states must vanish.

\begin{remark}
\label{rem:fluid-clock-properties}
Since $b>0$ and the cumulative controls of an admissible policy are continuous and nondecreasing, $f^\pi$ is a continuous strictly increasing bijection from $\mathbb R_+$ onto itself and therefore admits a continuous inverse. Moreover, $f^\pi$ has bounded variation on every bounded interval and hence induces a Lebesgue--Stieltjes measure $df^\pi$. All integrals with respect to $df^\pi$ below are understood in this sense.\end{remark}

Since admissibility implies that all intermediate-stage states are zero, the fluid control problem considers the analogue of the prelimit objective \eqref{eq:prelimit_obj}, seeking an admissible policy $\pi$ that minimizes
\begin{equation}\label{eq:limit_obj}
\mathcal J(\pi;z_0) = \int_0^\infty \sum_{j=1}^{J}c_jp_{j,L_j}^{\pi}(\xi)\,df^\pi(\xi).
\end{equation}
The cumulative holding cost is measured in clock time, with $df^\pi$ converting increments in step time into elapsed real-time increments.

\subsection{Stage-Level Collapse}
\label{subsec:stage-level-collapse}
In this subsection, we focus on the stage-level state-space collapse phenomenon and use it to simplify the fluid control problem. 
For an admissible policy, substituting the zero intermediate-stage constraints~\eqref{eq:fluid-state-constraints} into the stagewise flow-balance equations \eqref{eq:fluid-stage-flow-balance} gives
\begin{equation*}
\begin{aligned}
0&=p_{j,k}^{\pi}(\xi)=p_{j,k}(0)+\mathcal U_{j,k+1}(\xi)-\mathcal U_{j,k}(\xi)=\mathcal U_{j,k+1}(\xi)-\mathcal U_{j,k}(\xi),\quad k=1,\dots,L_j-1,\\[1ex]
0&=d_{j,o_j-1}^{\pi}(\xi)=d_{j,o_j-1}(0)+\mathcal U_{j,1}(\xi)-\mathcal V_{j,o_j-1}(\xi)=\mathcal U_{j,1}(\xi)-\mathcal V_{j,o_j-1}(\xi),\quad o_j\geq2,\\[1ex]
0&=d_{j,k}^{\pi}(\xi)=d_{j,k}(0)+\mathcal V_{j,k+1}(\xi)-\mathcal V_{j,k}(\xi)=\mathcal V_{j,k+1}(\xi)-\mathcal V_{j,k}(\xi),\quad k=1,\dots,o_j-2.
\end{aligned}
\end{equation*}
Thus, for every class $j$, define its cumulative completion process by
\begin{equation}
\mathcal X_j^\pi(\xi):=\mathcal U_{j,L_j}(\xi)=\cdots=\mathcal U_{j,1}(\xi).
\label{eq:class-level-service}
\end{equation}
For classes with $o_j\geq2$, this process also equals the cumulative completions at every decode stage:
\begin{equation*}
\mathcal X_j^\pi(\xi)=\mathcal V_{j,o_j-1}(\xi)=\cdots=\mathcal V_{j,1}(\xi).
\end{equation*}
These identities exhibit a stage-level state-space collapse, a phenomenon well studied in the heavy-traffic and fluid limits of queueing systems; see \citet[Chapter~7]{xu2023drift}. The fluid state is therefore fully described by the class-level waiting queues, and for each class \(j\), the cumulative completion process \(\mathcal{X}_j^\pi\) represents all stage-level controls.

Based on the stage-level collapse, we here reformulate the constraints and objective of the fluid control problem in terms of the cumulative class-$j$ request completions $\mathcal X_j^\pi(\xi)$.
First, define the cumulative KV-cache usage of a class-$j$ request, summed over its prefill and decode stages, by
\begin{equation}
\tilde m_j:=\sum_{k=1}^{L_j}\left[c_{\mathrm{chunk}}(L_j-k)+\ell_{j,k}\right]+\sum_{k=1}^{o_j-1}(s_j+o_j-k).
\label{eq:cumulative-kv-occupancy}
\end{equation}
Throughout the remainder of the paper, assume $\tilde m_j>0$ for every class $j$.

Substituting \eqref{eq:class-level-service} into the stage-level model and using $\sum_{k=1}^{L_j}\ell_{j,k}=s_j$ and \eqref{eq:cumulative-kv-occupancy} yields
\begin{equation}
p_{j,L_j}^{\pi}(\xi)=p_{j,L_j}(0)+\int_0^{f^\pi(\xi)}\lambda_j(s)\,ds-\mathcal X_j^\pi(\xi)\geq0.
\label{eq:class-level-waiting-state}
\end{equation}
\begin{equation}
f^\pi(\xi)=\sum_{j=1}^J\left(a_{\mathrm{pf}}s_j+a_{\mathrm{kv}}\tilde m_j\right)\mathcal X_j^\pi(\xi)+b\xi.
\label{eq:class-level-clock}
\end{equation}
For $\xi_1>\xi_2\geq0$,
\begin{equation}
\sum_{j=1}^J\tilde m_j\left(\mathcal X_j^\pi(\xi_1)-\mathcal X_j^\pi(\xi_2)\right) \leq M(\xi_1-\xi_2).
\label{eq:class-level-memory}
\end{equation}
Conversely, any collection of nondecreasing processes
$\{\mathcal X_j^\pi\}_{j=1}^J$ satisfying
the initial conditions $\mathcal X_j^\pi(0)=0$,
the queue nonnegativity constraint
\eqref{eq:class-level-waiting-state},
and the memory constraint
\eqref{eq:class-level-memory},
with the clock defined by
\eqref{eq:class-level-clock},
uniquely determines an admissible fluid policy
$\pi=(\mathcal U,\mathcal V)$. Therefore, we henceforth study the fluid control problem in terms of $\{\mathcal X_j^\pi\}_{j=1}^J$.

\section{Optimal Control of the Fluid Model}
\label{sec:fluid-optimal-control}

In this section, we study the $c\mu$ priority policy for the fluid model. We first introduce the basic idea of the $c\mu$ rule and define the effective service rate $\mu_j$. We then describe the fluid $c\mu$ policy and establish that it minimizes the cumulative holding cost. Finally, we show that the policy also minimizes the completion time and, consequently, maximizes the full-run output-token throughput.

\subsection{Effective Service Rates and the \texorpdfstring{$c\mu$}{c-mu} Rule}

The fluid control objective \eqref{eq:limit_obj} can be viewed as a linear holding cost on the queue lengths of waiting requests. Each unit of waiting class-$j$ requests incurs a holding cost at rate $c_j$. Serving class-$j$ requests at rate $\mu_j$ therefore reduces the holding-cost rate by $c_j\mu_j$ per unit of time. This motivates the classical $c\mu$ rule, which prioritizes classes in decreasing order of $c_j\mu_j$. To apply this rule to our fluid model and establish its optimality, we first identify the effective service rates determined by the memory constraint and the step-time model.

Consider an interval during which only class $j$ is served. The memory constraint in \eqref{eq:class-level-memory} gives, in differential form, that $\tilde m_j\,d\mathcal X_j^\pi(\xi)\leq M\,d\xi.$
Hence the class-$j$ service rate in step time is at most $M/\tilde m_j$, with equality when the KV capacity is saturated. At this saturated service rate, the amount of class-$j$ service completed per unit of step time is $M/\tilde m_j$, while the corresponding elapsed clock time is
\begin{equation*}
\left(a_{\mathrm{pf}}s_j+a_{\mathrm{kv}}\tilde m_j\right)\frac{M}{\tilde m_j}+b.
\end{equation*}
This motivates the following definition.

\begin{definition}[Effective service rate]
The effective service rate of class $j$ is
\begin{equation}
\mu_j:=\frac{M/\tilde m_j}{(a_{\mathrm{pf}}s_j+a_{\mathrm{kv}}\tilde m_j)M/\tilde m_j+b}=\left[a_{\mathrm{pf}}s_j+\left(a_{\mathrm{kv}}+\frac bM\right)\tilde m_j\right]^{-1}.
\label{eq:effective_service_rate}
    \end{equation}
\end{definition}

This definition of the effective service rate is key to our results. In particular, in $\mu_j^{-1}$, the first term $a_{\mathrm{pf}}s_j$ corresponds to prefill computation, while the second term $\left(a_{\mathrm{kv}}+\frac{b}{M}\right)\tilde m_j$ captures the combined effects of KV-cache I/O and memory capacity.

The same formula \eqref{eq:effective_service_rate} also characterizes the allocation of real-time service among multiple classes. Substituting the memory constraint \eqref{eq:class-level-memory} into the equation for the real-time clock \eqref{eq:class-level-clock} gives
\begin{equation*}
\begin{aligned}
df^\pi(\xi)&=\sum_{j=1}^J\left(a_{\mathrm{pf}}s_j+a_{\mathrm{kv}}\tilde m_j\right)d\mathcal X_j^\pi(\xi) +b\,d\xi\\
&\geq\sum_{j=1}^J\left[a_{\mathrm{pf}}s_j+\left(a_{\mathrm{kv}}+\frac{b}{M}\right)\tilde m_j\right]d\mathcal X_j^\pi(\xi)=\sum_{j=1}^J\mu_j^{-1}\,d\mathcal X_j^\pi(\xi),
\end{aligned}
\end{equation*}
with equality whenever the KV capacity is saturated. Here and below, $d\mathcal X_j^\pi/df^\pi$ denotes the Radon--Nikodym derivative of the Lebesgue--Stieltjes measure induced by $\mathcal X_j^\pi$ with respect to that induced by $f^\pi$. Define
\begin{equation*}
u_j^\pi(\xi):=\mu_j^{-1}\left(\frac{d\mathcal X_j^\pi}{df^\pi}\right)(\xi).
\end{equation*}
Then
\begin{equation*}
\left(\frac{d\mathcal X_j^\pi}{df^\pi}\right)(\xi)=\mu_j u_j^\pi(\xi),\qquad \sum_{j=1}^J u_j^\pi(\xi)\leq1,
\end{equation*}
where the inequality is an equality whenever the KV capacity is saturated. Here, $u_j^\pi$ is the fraction of real-time service effort allocated to class $j$, and $\mu_j$ is its effective service rate. Thus, the reduced GPU model has the same service-effort representation as the classical multiclass single-server model.

We apply the $c\mu$ rule using these effective service rates and index the classes so that
\[
c_1\mu_1\geq c_2\mu_2\geq\cdots\geq c_J\mu_J.
\]
The policy allocates service in this priority order, serving each class at the maximum rate permitted by its waiting requests, new arrivals, and the memory constraint \eqref{eq:class-level-memory}, and passing any remaining capacity to lower-priority classes. For a more detailed characterization of the $c\mu$ rule and its generalizations, see for example \citet{Buyukkoc1985Cmu}.

\subsection{Holding-Cost Optimality}
We measure the GPU's service budget in units of clock time: an interval of length $\Delta\tau$ provides $\Delta\tau$ units of budget. Since $\mu_j$ is the effective service rate of class $j$, processing $x$ class-$j$ requests requires $x/\mu_j$ units of effective work. At each time, the fluid $c\mu$ policy considers classes in decreasing order of $c_j\mu_j$. For each class, it uses the available budget to serve waiting requests and new arrivals, then passes any remaining budget to the next class. To establish the holding-cost optimality of the $c\mu$ rule, we use cumulative completion processes to represent the policy and construct these processes recursively in priority order.

For this recursive construction, let $A_j(\tau)$ denote the cumulative class-$j$ arrivals by clock time $\tau$:
\begin{equation}
A_j(\tau):=p_{j,L_j}(0)+\int_0^\tau\lambda_j(s)\,ds.
\label{eq:cumulative-demand}
\end{equation}

\begin{definition}[Fluid $c\mu$ policy]
\label{def:fluid-cmu-policy}
Set $R_0^*(\tau):=\tau$ for $\tau\geq0$.
Recursively, for $j=1,\dots,J$, define
\begin{equation}
R_j^*(\tau):=\sup_{0\leq s\leq\tau}\left[R_{j-1}^*(s)-\mu_j^{-1}A_j(s)\right]^+.
\label{eq:target-residual}
\end{equation}
Define
\begin{equation}
{\medmuskip=1.5mu\thickmuskip=2mu
F(\tau):=\frac1M\sum_{j=1}^J\tilde m_j\mu_j\left(R_{j-1}^*(\tau)-R_j^*(\tau)\right)+\frac1bR_J^*(\tau),}
\label{eq:target-F}
\end{equation}
and let the policy clock be its generalized inverse:
\begin{equation}
f^*(\xi):=\inf\left\{\tau\geq0:F(\tau)\geq\xi\right\}.
\label{eq:target-clock}
\end{equation}
For each class $j=1,\dots,J$, set
\begin{equation}
\mathcal X_j^*(\xi) :=\mu_j\left(R_{j-1}^*(f^*(\xi))-R_j^*(f^*(\xi))\right),
\label{eq:target-admission}
\end{equation}
and define $\pi^*:=(\mathcal U^*,\mathcal V^*)$ by
\begin{equation*}
\mathcal U_{j,k}^{*}=\mathcal X_j^*,\quad k=1,\dots,L_j,\quad\text{and}\quad \mathcal V_{j,k}^{*}=\mathcal X_j^*,\quad k=1,\dots,o_j-1.
\end{equation*}
\end{definition}

\begin{remark}[Policy interpretation]
The basic idea is to use Skorokhod reflection to represent the remaining service budget passed to the next class and use the effective service rate $\mu_j$ to convert the budget used by class $j$ into its cumulative completion process. Specifically, $R_0^*(\tau)=\tau$ is the total service budget available by clock time $\tau$. After classes $1,\dots,j$ have been considered, $R_j^*(\tau)$ is the cumulative unused service budget passed to the lower-priority classes. Hence, $R_{j-1}^*(\tau)-R_j^*(\tau)$ is the service budget used by class $j$. Since these budgets are indexed by clock time $\tau$, we evaluate this difference at $\tau=f^*(\xi)$. Multiplying it by $\mu_j$ then gives $\mathcal X_j^*(\xi)$, the cumulative number of class-$j$ completions by step time $\xi$.
\end{remark}

\begin{theorem}[Fluid holding-cost optimality]
\label{thm:fluid-cmu-optimality}
The fluid $c\mu$ policy $\pi^*$ is admissible and minimizes the cumulative holding cost. Specifically, for every admissible policy $\pi$,
\begin{equation}
\mathcal J(\pi^*;z_0)\leq\mathcal J(\pi;z_0).
\label{eq:fluid-cmu-optimality}
\end{equation}
\end{theorem}

The proof is given in Appendix~\ref{sec:proof-fluid-cmu-optimality}. 

\begin{figure}[H]
\centering
\includegraphics{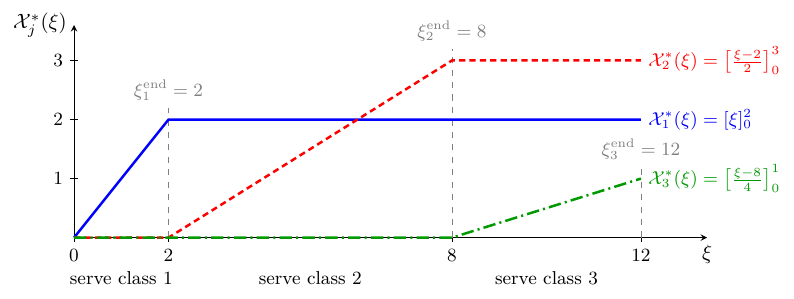}
\caption{Illustration of the fluid $c\mu$ policy for $J=3$, assuming classes are already indexed in priority order. Here $M=1$, $(p_{1,L_1}(0),p_{2,L_2}(0),p_{3,L_3}(0))=(2,3,1)$, and $(\tilde m_1,\tilde m_2,\tilde m_3)=(1,2,4)$.}
\label{fig:priority-policy-example-ordered}
\end{figure}

To better illustrate the priority rule, we consider the case with no new arrivals. We show in Appendix~\ref{sec:no-arrival-closed-forms} that the policy admits the following closed-form representation:
\begin{equation}
\mathcal X_j^*(\xi)=\left[\frac{M}{\tilde m_j}\left(\xi-\sum_{\ell=1}^{j-1}\frac{\tilde m_\ell p_{\ell,L_\ell}(0)}{M}\right)\right]_0^{p_{j,L_j}(0)}.
\label{eq:no-arrival-explicit-policy}
\end{equation}
Here $[x]_0^y:=\max\{0,\min\{x,y\}\}$. In~\eqref{eq:no-arrival-explicit-policy}, the slope $M/\tilde m_j$ represents the service rate of class $j$ when it occupies the full memory capacity. At this rate, completing its initial requests takes $\tilde m_j p_{j,L_j}(0)/M$ units of step time. The shift $\sum_{\ell=1}^{j-1}\tilde m_\ell p_{\ell,L_\ell}(0)/M$ therefore equals the total service time required by the preceding classes. The truncation keeps $\mathcal X_j^*(\xi)$ between $0$ and $p_{j,L_j}(0)$. Thus, the formula serves the classes sequentially in index order, which is the decreasing order of $c_j\mu_j$ prescribed by the $c\mu$ rule. Figure~\ref{fig:priority-policy-example-ordered} illustrates the $c\mu$ priority rule for three classes, with each class beginning service when the preceding class has been completed. For an illustration under general arrivals, see Appendix~\ref{sec:general-arrival-policy-illustration}.

\subsection{Completion Time and Throughput Optimality}
\label{subsec:fluid-throughput}

For an admissible policy $\pi$, the clock $f^\pi$ admits a continuous inverse, as noted in Remark~\ref{rem:fluid-clock-properties}. Therefore we define its average output-token throughput over the real-time interval $[0,\tau]$ by
\begin{equation}
{\medmuskip=1mu\thickmuskip=2mu
\begin{aligned}
\operatorname{TP}^\pi(\tau)&:=\frac{1}{\tau}\sum_{j=1}^J\left[\mathcal U_{j,1}\!\bigl((f^\pi)^{-1}(\tau)\bigr)+\sum_{k=1}^{o_j-1}\mathcal V_{j,k}\!\bigl((f^\pi)^{-1}(\tau)\bigr)\right]=\frac{1}{\tau}\sum_{j=1}^Jo_j\mathcal X_j^\pi\!\bigl((f^\pi)^{-1}(\tau)\bigr),\qquad \tau>0.
\end{aligned}}
\label{eq:fluid-throughput}
\end{equation}
The second equality follows from \eqref{eq:class-level-service}.

Because $\lambda_j$ vanishes after $H$, $A_j(H)$ is the total number of class-$j$ requests. For an admissible policy $\pi$ that eventually completes all requests, define its completion time as
\begin{equation}
\tau_{\mathrm{c}}^\pi:=\inf\left\{\tau\geq0:\mathcal X_j^\pi\!\left((f^\pi)^{-1}(\tau)\right)=A_j(H)\right\}.
\label{eq:fluid-completion-time}
\end{equation}
The equality must hold for all classes simultaneously.
If $\sum_{j=1}^Jo_jA_j(H)>0$, then $\tau_{\mathrm{c}}^\pi>0$ and
\begin{equation}
\operatorname{TP}^\pi(\tau_{\mathrm{c}}^\pi)=\frac{\sum_{j=1}^Jo_jA_j(H)}{\tau_{\mathrm{c}}^\pi}.
\label{eq:fluid-throughput-at-completion}
\end{equation}

\begin{theorem}[Fluid completion-time and throughput optimality]
\label{thm:fluid-throughput-optimality}
The $c\mu$ policy $\pi^*$ in Definition~\ref{def:fluid-cmu-policy} minimizes the time needed to complete all requests. Specifically, for every admissible policy $\pi$ that eventually completes all requests,
\begin{equation}
\tau_{\mathrm{c}}^{\pi^*}\leq\tau_{\mathrm{c}}^\pi.
\label{eq:fluid-completion-time-optimality}
\end{equation}
When $\sum_{j=1}^Jo_jA_j(H)>0$, completion-time optimality also implies that $\pi^*$ maximizes the full-run output-token throughput:
\begin{equation}
\operatorname{TP}^\pi(\tau_{\mathrm{c}}^\pi)\leq\operatorname{TP}^{\pi^*}(\tau_{\mathrm{c}}^{\pi^*}).
\label{eq:fluid-throughput-optimality}
\end{equation}
\end{theorem}
The proof is given in Appendix~\ref{sec:proof-fluid-throughput-optimality}.

\section{Fluid-Limit Convergence}
\label{sec:prelimit-convergence}

In this section, we demonstrate that the optimal fluid $c\mu$ policy provides operational guidance for the prelimit $r$th system. We construct a smoothed $c\mu$-inspired dynamic policy that is admissible at every finite scale $r$. At a high level, it follows the $c\mu$ priority order when forming admission plans and spreads each plan over a short window to satisfy the KV-cache capacity constraint. In particular, this policy does not require knowledge of the arrival rates $\lambda_j(\cdot)$. Then, we show that under this policy and fluid scaling, the prelimit system and control converge to the fluid system and control under the optimal $c\mu$ policy. 

To track the memory used by a request that runs continuously after admission, define the reserved KV-cache usage of a class-$j$ request at the start of the GPU step occurring $h$ steps after admission by
\begin{equation*}
{\renewcommand{\arraystretch}{1.2}\medmuskip=2mu\thickmuskip=2.5mu k_j^{\mathrm{kv}}(h):=\left\{\begin{array}{@{}l@{}}
\min\{c_{\mathrm{chunk}}(h+1),s_j\},\enspace 0\leq h\leq L_j-1,\\
s_j+h-L_j+1,\enspace L_j\leq h\leq L_j+o_j-2,\\
0,\enspace\text{otherwise}.
\end{array}\right.}
\end{equation*}
During prefill, this includes the KV slots reserved for the current chunk. During decode, it includes the prompt KV states, the KV states of previously generated tokens, and the slot reserved for the current input token. It becomes zero after the request completes. Thus, with
\begin{equation*}
H_{\mathrm{kv}}:=\max_{1\leq j\leq J}(L_j+o_j-2),
\end{equation*}
the cumulative KV-cache usage $\tilde m_j$ satisfies
\begin{equation}
\sum_{h=0}^{H_{\mathrm{kv}}}k_j^{\mathrm{kv}}(h)=\sum_{k=1}^{L_j}\left[c_{\mathrm{chunk}}(L_j-k)+\ell_{j,k}\right]+\sum_{k=1}^{o_j-1}(s_j+o_j-k)=\tilde m_j.
\label{eq:prelimit-mtilde-stage-identity}
\end{equation}

The KV-cache usage at step $t$ is then
\begin{equation}
m^{r}(t) = \sum_{j=1}^{J}\sum_{h=0}^{H_{\mathrm{kv}}} k_j^{\mathrm{kv}}(h)\Delta\mathcal U_{j,L_j}^{r}(t-h). \label{eq:prelimit-memory-convolution}
\end{equation}
For class $j$, the requests admitted at step $t-h$ contribute $k_j^{\mathrm{kv}}(h)\Delta\mathcal U_{j,L_j}^{r}(t-h)$ to the KV-cache usage at step $t$. Summing these contributions over all classes and lags gives $m^r(t)$. Using the priority order $c_1\mu_1\geq\cdots\geq c_J\mu_J$, we define the smoothed $c\mu$-inspired policy as follows.

\begin{definition}[Smoothed $c\mu$-inspired policy]
\label{def:moving-average-policy}

Choose integer-valued $w_r$ such that
\begin{equation*}
w_r\longrightarrow\infty, \qquad \frac{w_r}{r}\longrightarrow0, \qquad w_r>H_{\mathrm{kv}};
\end{equation*}
set $y_j^{r}(t)=\Delta\mathcal U_{j,L_j}^{r}(t)=0$ for $t<0$, and define the new class-$j$ planned-admission amount recursively in priority order by
\begin{equation}
{\medmuskip=2.5mu\thickmuskip=3mu y_j^r(t)=\min\Biggl\{\left[p_{j,L_j}^r(0)+\sum_{u<t}\eta_j^r(u)-\sum_{u<t}y_j^r(u)\right]^+,\tilde m_j^{-1}\left[\frac{w_rM}{w_r+H_{\mathrm{kv}}}-\sum_{\ell<j}\tilde m_\ell y_\ell^r(t)\right]^+\Biggr\}.}
\label{eq:planned-admission-policy}
\end{equation}
Construct the admission increment directly by the moving average
\begin{equation}
\Delta\mathcal U_{j,L_j}^{r}(t) \;=\; \frac{1}{w_r}\sum_{v=0}^{w_r-1}y_j^{r}(t-v). \label{eq:smoothed-admission}
\end{equation}
Recall that we allow fractional requests, so both the planned-admission amounts $y_j^r(t)$ and the admission increments $\Delta\mathcal U_{j,L_j}^r(t)$ may take noninteger values. 
For $k=1,\dots,L_j-1$, we have
\begin{equation}
\Delta\mathcal U_{j,k}^r(t)=\Delta\mathcal U_{j,L_j}^r\bigl(t-(L_j-k)\bigr).
\label{eq:prelimit-U-shift}
\end{equation}
For $k=1,\dots,o_j-1$, we have
\begin{equation}
\Delta\mathcal V_{j,k}^r(t)=\Delta\mathcal U_{j,L_j}^r\bigl(t-(L_j+o_j-1-k)\bigr).
\label{eq:prelimit-V-shift}
\end{equation}
\end{definition}

\begin{remark}
The choice $w_r=\max\{H_{\mathrm{kv}}+1,\lceil\sqrt r\rceil\}$
satisfies the required conditions. The condition $w_r\to\infty$ smooths the step-level admission plans, whereas $w_r/r\to0$ ensures that the smoothing window vanishes in fluid time. At step $t$, $y_j^r(t)$ specifies how many class-$j$ requests are assigned an admission plan. Equation~\eqref{eq:smoothed-admission} schedules $y_j^r(t)/w_r$ of this plan to enter the GPU in each of steps $t,\dots,t+w_r-1$, so the plan is carried out over a $w_r$-step window. Once admitted, each request runs continuously until completion; its subsequent stage-completion increments are therefore the time shifts in \eqref{eq:prelimit-U-shift}--\eqref{eq:prelimit-V-shift}. In Equation~\eqref{eq:planned-admission-policy}, the first term is the number of class-$j$ requests available at step $t$ that have not yet received an admission plan, while the second is the remaining per-step planning budget expressed in class-$j$ requests. The first term can equivalently be written using the waiting-queue length as
\begin{equation*}
\left[ p_{j,L_j}^{r}(t) -\frac{1}{w_r} \sum_{v=1}^{w_r-1}(w_r-v)y_j^{r}(t-v) \right]^+.
\end{equation*}
The subtracted sum is the portion of earlier admission plans that has not yet entered the GPU, so pending plans are not included when a new plan is formed.
\end{remark}
The following proposition establishes the admissibility of the policy in the prelimit queueing model.

\begin{proposition}[Admissibility]
\label{prop:pathwise-memory-safety}
Suppose that, initially, some requests may be waiting, but no requests are being processed on the GPU. Then, the policy in Definition~\ref{def:moving-average-policy} is admissible for the prelimit queueing model. 
\end{proposition}

Proposition~\ref{prop:pathwise-memory-safety} establishes policy admissibility at every finite scale (proof in Appendix~\ref{sec:proof-pathwise-memory-safety}). We next show that the smoothed $c\mu$-inspired policy converges to the fluid $c\mu$ policy under fluid scaling.

\begin{theorem}[Fluid-limit convergence]
\label{thm:moving-average-fluid-limit}
Suppose initially no requests are being processed on the GPU, and the initial waiting queue lengths satisfy
\begin{equation*}
\frac{1}{r}p_{j,L_j}^{r}(0)\longrightarrow p_{j,L_j}(0),\qquad j=1,\dots,J.
\end{equation*}
Consider prelimit systems \eqref{eq:prelimit-dynamics} controlled by the smoothed $c\mu$-inspired policy in Definition~\ref{def:moving-average-policy}. Then, for every $T<\infty$ and every class $j$, the prelimit cumulative control processes and system clock converge to the limiting $c\mu$-policy and the corresponding clock: 
\begingroup
\setlength{\abovedisplayskip}{3pt}
\setlength{\belowdisplayskip}{3pt}
\setlength{\abovedisplayshortskip}{0pt}
\setlength{\belowdisplayshortskip}{3pt}
{\setlength{\jot}{0pt}
\begin{align}
\sup_{0\leq\xi\leq T}
\left|\frac1r\mathcal U_{j,k}^r(\lfloor r\xi\rfloor)
-\mathcal X_j^*(\xi)\right|
&\xrightarrow{\mathbb P}0,
\quad k=1,\dots,L_j,
\label{eq:prelimit-U-convergence}\\
\sup_{0\leq\xi\leq T}
\left|\frac1r\mathcal V_{j,k}^r(\lfloor r\xi\rfloor)
-\mathcal X_j^*(\xi)\right|
&\xrightarrow{\mathbb P}0,
\quad k=1,\dots,o_j-1;
\label{eq:prelimit-V-convergence}\\
\sup_{0\leq\xi\leq T}\left|\bar f^r(\xi)-f^*(\xi)\right|&\xrightarrow{\mathbb P}0.\label{eq:prelimit-clock-convergence}
\end{align}
}
\endgroup
Moreover, the state processes converge to that of the limiting system under the optimal $c\mu$-policy:
\begingroup
\setlength{\abovedisplayskip}{3pt}
\setlength{\belowdisplayskip}{3pt}
\setlength{\abovedisplayshortskip}{0pt}
\setlength{\belowdisplayshortskip}{3pt}
{\setlength{\jot}{0pt}
\begin{align}
\sup_{0\leq\xi\leq T}\left|\frac1r p_{j,k}^r(\lfloor r\xi\rfloor)-p_{j,k}^{\pi^*}(\xi)\right|&\xrightarrow{\mathbb P}0,\quad k=1,\dots,L_j,\label{eq:prelimit-p-convergence}\\
\sup_{0\leq\xi\leq T}\left|\frac1r d_{j,k}^r(\lfloor r\xi\rfloor)-d_{j,k}^{\pi^*}(\xi)\right|&\xrightarrow{\mathbb P}0,\quad k=1,\dots,o_j-1.\label{eq:prelimit-d-convergence}
\end{align}
}
\endgroup
\end{theorem}

The proof of Theorem~\ref{thm:moving-average-fluid-limit} is deferred to Appendix~\ref{sec:proof-moving-average-fluid-limit}. The limit $(f^*,\mathcal X^*)$ is the fluid $c\mu$ policy in Definition~\ref{def:fluid-cmu-policy}. It is holding-cost optimal by Theorem~\ref{thm:fluid-cmu-optimality} and completion-time and throughput optimal by Theorem~\ref{thm:fluid-throughput-optimality}. Thus, the smoothing used at finite scale disappears under fluid scaling, and the real-time clock, system states, and the fluid-scaled cumulative stage-completion processes converge uniformly on compacts in probability to the optimal fluid trajectory. This convergence, together with the optimality of the limiting trajectory, yields the following corollary on the convergence of policy values. The proof is given in Appendix~\ref{sec:proof-moving-average-cost-convergence}.
\begin{corollary}[Holding-cost convergence]
\label{cor:moving-average-cost-convergence}
Under the assumptions of Theorem~\ref{thm:moving-average-fluid-limit}, the smoothed $c\mu$-inspired policy $\widehat\pi^r$ in Definition~\ref{def:moving-average-policy} satisfies
\begin{equation}
\mathcal J^r(\widehat\pi^r;z_0^r)\longrightarrow\mathcal J(\pi^*;z_0).
\label{eq:moving-average-cost-convergence}
\end{equation}
\end{corollary}
\section{Numerical Experiments}
\label{sec:numerical}

In this section, we evaluate a practical implementation of the $c\mu$ rule on a real LLM serving engine. We set $c=1$ to assign equal weight to all classes and use workloads derived from a real-world dataset~\citep{Wang2025KVCacheWild}. 
We conduct experiments on semi-synthetic workloads in Section~\ref{sec:generated-workloads}, followed by experiments in a realistic serving environment under known prompt length and unknown output lengths in Section \ref{sec:compressed-traces}.

\subsection{Experimental Setup}
\label{sec:numerical-datasets}

We run Llama-3.1-8B-Instruct~\citep{Grattafiori2024Llama3} using vLLM~0.34.0~\citep{Kwon2023PagedAttention} on a single NVIDIA H800 GPU with 364,736 usable KV-cache tokens and chunked prefill enabled. All reported times in this section include the scheduler's runtime, which accounts for around 4\% of total clock time.
Within each experimental group, all policies are run sequentially on the same GPU. We calibrate the step-time model separately for each GPU to estimate $a_{\mathrm{pf}}$, $a_{\mathrm{kv}}$, and $b$. Substituting these estimates and request $i$'s prompt and output lengths into the effective-service-rate formula~\eqref{eq:effective_service_rate} gives its estimated service rate $\widehat\mu_i$ for the $c\mu$ ordering.

We use the To-B, To-C, Thinking, and Coder traces from the Qwen-Bailian Anonymous Dataset~\citep{Wang2025KVCacheWild}. To-B covers business-facing API calls and task automation, while To-C covers consumer-facing interactive chat. Thinking focuses on reasoning-intensive chat, and Coder covers code generation and interactive programming. These anonymized production traces each cover two hours. Each record provides a relative arrival time and prompt and output lengths. We retain only requests whose combined prompt and output lengths do not exceed the model's 131,072-token context window. We adjust the arrival rates to achieve a traffic intensity of approximately 0.95. For replay, we configure the LLM engine to generate exactly the specified number of output tokens for each request. Table~\ref{tab:numerical-workloads} summarizes their lengths and eligible request counts.
\begin{table}[!ht]
\centering
\caption{Eligible requests in the four production traces. Prompt and output lengths are measured in tokens.}
\label{tab:numerical-workloads}
\begin{tabular*}{\textwidth}{@{\extracolsep{\fill}}lccc@{}}
\toprule
Trace & Requests & Mean prompt length & Mean output length \\
\midrule
To-B & 172,800 & 914.5 & 91.6 \\
To-C & 43,058 & 2,331.3 & 430.5 \\
Thinking & 10,812 & 4,762.9 & 3,650.8 \\
Coder & 43,011 & 5,747.7 & 808.3 \\
\bottomrule
\end{tabular*}
\end{table}

\begin{figure}[htbp]
\centering
\begin{minipage}[t]{.49\textwidth}
\centering
\includegraphics[width=\linewidth]{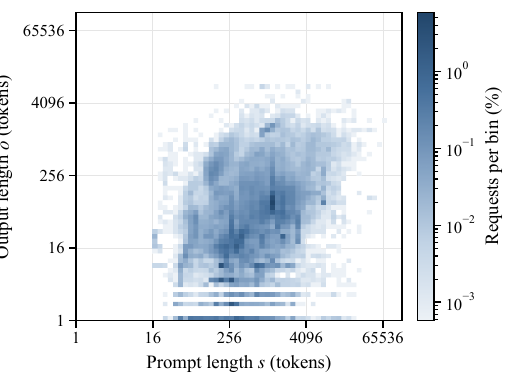}
\par\smallskip
{\small (a) To-B}
\end{minipage}\hfill
\begin{minipage}[t]{.49\textwidth}
\centering
\includegraphics[width=\linewidth]{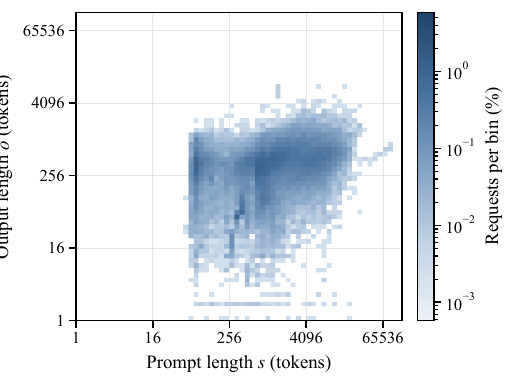}
\par\smallskip
{\small (b) To-C}
\end{minipage}

\medskip
\begin{minipage}[t]{.49\textwidth}
\centering
\includegraphics[width=\linewidth]{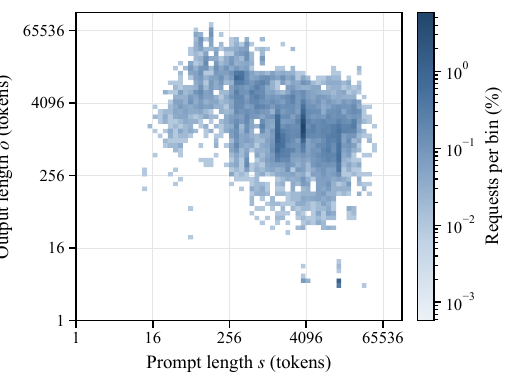}
\par\smallskip
{\small (c) Thinking}
\end{minipage}\hfill
\begin{minipage}[t]{.49\textwidth}
\centering
\includegraphics[width=\linewidth]{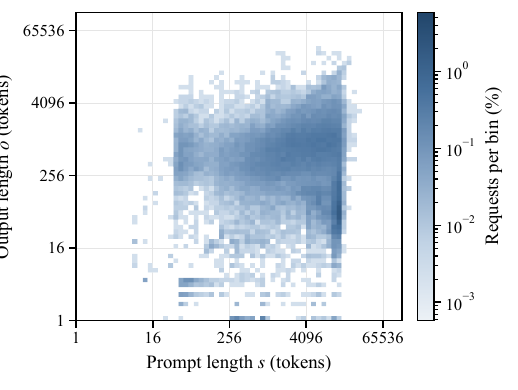}
\par\smallskip
{\small (d) Coder}
\end{minipage}
\caption{Joint prompt--output length distributions of eligible requests in (a) To-B, (b) To-C, (c) Thinking, and (d) Coder. Both axes use logarithmic token scales.}
\label{fig:production-lengths}
\end{figure}

Figure~\ref{fig:production-lengths} compares the joint prompt--output length distributions. On average, To-B has the shortest prompts and outputs, followed by To-C. Coder has the longest prompts and Thinking has the longest outputs. We use To-B to construct the semi-synthetic workloads in Section~\ref{sec:generated-workloads}. In Section~\ref{sec:compressed-traces}, we combine requests from all four traces into one workload and replay it after arrival-time rescaling.

\subsection{Policies and Performance Measures}
Based on the $c\mu$ priority rule, we implement a practical scheduling policy denoted by $c\mu$. We compare it with the default vLLM scheduler, MC-SF~\citep{Jaillet2025Online}, SVF~\citep{Kong2026Geometry}, and WAIT-frac, a fractional implementation of WAIT~\citep{Ao2026WAIT}. The $c\mu$, MC-SF, and SVF policies order waiting requests by priority. In Section~\ref{sec:generated-workloads}, all policies except the default vLLM scheduler use nonpreemptive implementations and assume that each request's output length is known. In Section~\ref{sec:compressed-traces}, we evaluate $c\mu$, MC-SF, and SVF using predicted output lengths and vLLM's default preemption mechanism, alongside the default vLLM scheduler.

The $c\mu$, MC-SF, and SVF policies check future KV-cache usage before admission, accounting for cache growth and memory released upon completion. A request is admitted only if the projected total KV-cache usage remains within $M$ at every future GPU step.

Using the future KV-cache feasibility check, the $c\mu$ policy scans waiting requests in decreasing $c_i\widehat\mu_i$ order and stops at the first infeasible request.

MC-SF orders waiting requests by increasing output length $o_i$, while SVF orders them by increasing cumulative KV-cache volume $s_io_i+o_i(o_i+1)/2$. Both use the same future KV-cache feasibility check. When $c_i=1$, $a_{\mathrm{pf}}=0$, and prefill is not chunked, the $c\mu$ rule reduces to a volume-based priority rule, as in SVF.

The vanilla version of WAIT described in \cite{Ao2026WAIT} uses positive integer service quotas for each class and stage, so the smallest positive quota is one request per GPU step. When the quota is one, a single waiting request is sufficient to make its class eligible for service in the next GPU step, eliminating the intended waiting mechanism for request accumulation. In our experiments, however, matching class-wise arrival rates may require service quotas below one request per GPU step. We therefore introduce WAIT-frac, which maintains fractional credits separately for each class and stage. Serving one request at a stage consumes one unit of credit; fractional remainders carry over to subsequent steps, subject to an expiry rule. Runtime KV-cache checks enforce memory feasibility. 
In WAIT-frac, the reference step time $T_{\mathrm{WAIT}}$ converts class-wise arrival rates into the credit increments for each GPU step. Under our GPU configuration, the original parameter formula for WAIT does not yield a positive value of $T_{\mathrm{WAIT}}$, because \citet{Ao2026WAIT} assumes a different batch service-time model.
We therefore calibrate $T_{\mathrm{WAIT}}$ automatically for each trace by matching it to the mean step time during peak load.

For a trace containing $N$ requests, let $\tau_i^{\mathrm{arr}}$ and $\tau_i^{\mathrm{comp}}$ be the arrival and completion times of request $i$. We define the total request latency as
\[
\operatorname{Latency}_{\mathrm{tot}}:=\sum_{i=1}^N\left(\tau_i^{\mathrm{comp}}-\tau_i^{\mathrm{arr}}\right).
\]
We set $c_i=1$ for every request, so $cQ(\tau)$ counts all arrived but unfinished requests, including those on the GPU. Our primary performance measure is the average request latency:
\begin{equation*}
\operatorname{Latency}_{\mathrm{avg}}:=\frac{1}{N}\operatorname{Latency}_{\mathrm{tot}}=\frac{1}{N}\int_0^\infty cQ(\tau)\,d\tau.
\end{equation*}

\begin{remark}
Request $i$ contributes one to $cQ(\tau)$ exactly when $\tau_i^{\mathrm{arr}}\leq\tau<\tau_i^{\mathrm{comp}}$. Summing these contributions and integrating gives
\[
\int_0^\infty cQ(\tau)\,d\tau=\sum_{i=1}^N\int_{\tau_i^{\mathrm{arr}}}^{\tau_i^{\mathrm{comp}}}1\,d\tau=\sum_{i=1}^N\left(\tau_i^{\mathrm{comp}}-\tau_i^{\mathrm{arr}}\right)=\operatorname{Latency}_{\mathrm{tot}}.
\]
\end{remark}

We also report the mean time to first token (TTFT), measured from arrival to the first output, and the 99th percentile of request latency (P99 latency). Clearance time is $\max_i\tau_i^{\mathrm{comp}}-\min_i\tau_i^{\mathrm{arr}}$. Output-token throughput is the total number of generated tokens divided by the time from experiment start to the last completion.

\subsection{Semi-Synthetic Workloads}
\label{sec:generated-workloads}

We estimate a Gaussian mixture model to cluster eligible To-B requests into 16 classes. For each class, we generate Poisson arrivals with time-varying rates based on the source trace and scaled by a common factor to achieve an overall traffic intensity of approximately 0.95. We then assume that all requests within each class have identical prompt and decode lengths, given by their respective cluster centers.
We repeat the experiment eight times, replaying the same generated trace for every policy within each repetition. 
Table~\ref{tab:generated-policy-comparison} summarizes the results across all eight repetitions. The $c\mu$ and SVF policies have nearly identical average latencies. The $c\mu$ policy achieves the lowest mean TTFT and P99 latency. All policies perform similarly in terms of clearance time and output token throughput.

\begin{table}[H]
\centering
\small
\setlength{\tabcolsep}{5pt}
\caption{Scheduling results on eight semi-synthetic To-B traces. Average latency and clearance are reported as the mean $\pm$ sample standard deviation across traces. The remaining columns report means, with P99 latency computed separately on each trace. Bold values are best in each column.}
\label{tab:generated-policy-comparison}
\begin{tabular*}{\textwidth}{@{\extracolsep{\fill}}lrrrrr@{}}
\toprule
Policy & \shortstack{Average\\latency (s)} & \shortstack{Mean\\TTFT (s)} & \shortstack{P99\\latency (s)} & \shortstack{Clearance\\(s)} & \shortstack{Output\\tokens/s} \\
\midrule
$c\mu$ & $\boldsymbol{8.412\pm0.573}$ & \textbf{2.463} & \textbf{67.480} & $3{,}625.4\pm13.2$ & 2,431.2 \\
SVF & $\boldsymbol{8.412\pm0.565}$ & 2.485 & 90.707 & $3{,}624.2\pm13.7$ & 2,432.0 \\
MC-SF & $9.503\pm0.612$ & 3.703 & 143.020 & $\boldsymbol{3{,}623.9\pm13.2}$ & \textbf{2,432.3} \\
WAIT-frac & $12.462\pm0.892$ & 6.819 & 91.595 & $3{,}682.6\pm22.4$ & 2,393.5 \\
vLLM & $29.356\pm5.068$ & 23.716 & 97.215 & $3{,}644.9\pm18.8$ & 2,418.2 \\
\bottomrule
\end{tabular*}
\end{table}

\subsection{Real Production Traces}
\label{sec:compressed-traces}

We combine all four production traces to evaluate scheduling with predicted output lengths. We train the predictors on 60\% of the data and test their accuracy on the remaining 40\%. For scheduling, we randomly sample 20\% of the test requests without replacement and replay them in arrival order.

To construct the \(c\mu\) rule when the output length is unknown, we estimate the mean effective service time,  using the request's input features. This can be achieved by estimating \(\tilde m\) \eqref{eq:cumulative-kv-occupancy}. For each request \(i\), let \(\widehat{O_i}\) and \(\widehat{O_i^2}\) denote the estimated first and second moments of its output length, respectively. We then compute the corresponding estimate \(\widehat{m_i}\) as follows:
\[
\begin{aligned}
\widehat {m_i}=\frac{c_{\mathrm{chunk}}L_i(L_i-1)}{2}+s_i\widehat{O_i}+\frac{\widehat{O_i^2}-\widehat{O_i}}{2}.
\end{aligned}
\]
The SVF and MC-SF policies follow the same estimation procedure.
\begin{figure}[!htbp]
\centering
\includegraphics[width=\textwidth]{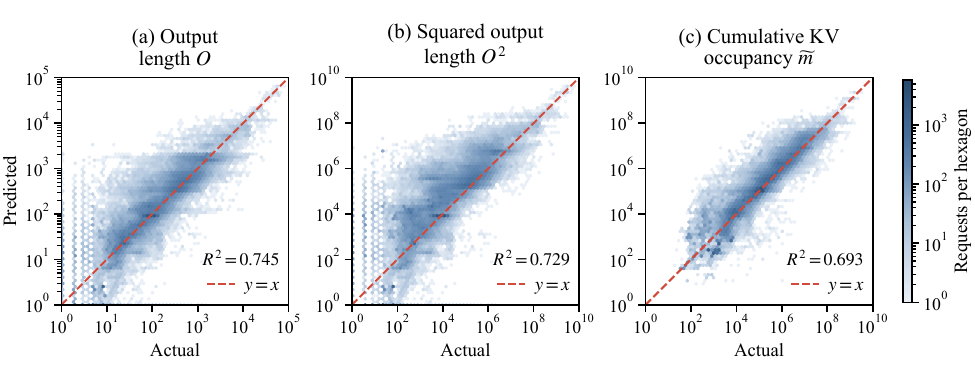}
\caption{Actual versus predicted values on the test set: (a) output length, (b) squared output length, and (c) cumulative KV-cache usage. Color indicates the number of requests per hexagon.}
\label{fig:prediction-accuracy}
\end{figure}

The out-of-sample coefficient of determination ($R^2$) for predictions of $O_i$, $O_i^2$, and $\tilde m_i$ are 0.745, 0.729, and 0.693, respectively.
Figure~\ref{fig:prediction-accuracy} compares actual and predicted values using logarithmic axes.  The prediction time is \(131.2 \pm 5.7\,\mu\mathrm{s}\), which is negligible compared with the average latency. Although output length predictions are imperfect, we show that the $c\mu$ rule still outperforms the baseline policies in terms of serving decisions.

For the scheduling experiments, we scale the interarrival times of the mixed test workload to target a traffic intensity of $\rho=0.95$. The $c\mu$ policy and SVF use both $\widehat{O_i}$ and $\widehat{O_i^2}$ to determine priorities, while MC-SF uses only $\widehat{O_i}$. All three policies retain the future KV-cache feasibility check, using $\lceil\widehat{O_i}\rceil$ as the predicted output length. We use vLLM's default preemption mechanism to avoid exceeding the KV-cache capacity.

The $c\mu$ policy achieves the lowest average latency, mean TTFT, and P99 latency. Specifically, its average latency is 7.7\% lower than that of SVF, 13.5\% lower than that of MC-SF, and 88.6\% lower than that of vLLM. These results demonstrate the robustness and the effectiveness of our proposed $c\mu$ policy across nearly all performance metrics.

\begin{table}[H]
\centering
\small
\setlength{\tabcolsep}{5pt}
\caption{Scheduling results on the mixed test workload. Clearance is the time from the first arrival until all requests complete. Bold values are best among the compared policies.}
\label{tab:real-policy-comparison}
\begin{tabular*}{\textwidth}{@{\extracolsep{\fill}}lrrrr@{}}
\toprule
Policy & \shortstack{Average\\latency (s)} & \shortstack{Mean\\TTFT (s)} & \shortstack{P99\\latency (s)} & \shortstack{Clearance\\(s)} \\
\midrule
$c\mu$ & \textbf{17.787} & \textbf{4.716} & \textbf{226.702} & 6,425.2 \\
SVF & 19.268 & 6.223 & 245.163 & 6,450.5 \\
MC-SF & 20.557 & 6.958 & 258.573 & 6,438.3 \\
vLLM & 156.035 & 143.388 & 533.507 & \textbf{6,211.0} \\
\bottomrule
\end{tabular*}
\end{table}

\FloatBarrier

\section{Conclusion and Future Work}
\label{sec:conclusion}
In this paper, we propose a simple modification of the $c\mu$ rule and establish its asymptotic optimality for LLM inference scheduling under memory constraints. Interesting directions for future research include extending the queueing model to multi-GPU settings and agentic workloads, as well as characterizing scheduling sensitivity to output-length prediction errors and improving policy robustness to these errors.

\bibliographystyle{plainnat}
\bibliography{references}

@book{Billingsley1999,
  author    = {Patrick Billingsley},
  title     = {Convergence of Probability Measures},
  edition   = {2},
  publisher = {Wiley},
  year      = {1999}
}

@book{JacodShiryaev2003,
  author    = {Jean Jacod and Albert N. Shiryaev},
  title     = {Limit Theorems for Stochastic Processes},
  edition   = {2},
  publisher = {Springer},
  year      = {2003},
  doi       = {10.1007/978-3-662-05265-5}
}

@book{Whitt2002,
  author    = {Ward Whitt},
  title     = {Stochastic-Process Limits: An Introduction to Stochastic-Process Limits and Their Application to Queues},
  publisher = {Springer},
  year      = {2002}
}

@misc{Whitt2002Supplement,
  author    = {Ward Whitt},
  title     = {Internet Supplement to {Stochastic-Process Limits: An Introduction to Stochastic-Process Limits and Their Application to Queues}},
  year      = {2002},
  url       = {https://www.columbia.edu/~ww2040/supplement.html}
}

@book{Rudin1976,
  author    = {Walter Rudin},
  title     = {Principles of Mathematical Analysis},
  edition   = {3},
  publisher = {McGraw--Hill},
  year      = {1976}
}

@inproceedings{Yu2022Orca,
  author    = {Gyeong-In Yu and Joo Seong Jeong and Geon-Woo Kim and Soojeong Kim and Byung-Gon Chun},
  title     = {Orca: A Distributed Serving System for {Transformer-Based} Generative Models},
  booktitle = {16th USENIX Symposium on Operating Systems Design and Implementation (OSDI 22)},
  pages     = {521--538},
  publisher = {USENIX Association},
  year      = {2022},
  url       = {https://www.usenix.org/conference/osdi22/presentation/yu}
}

@inproceedings{Kwon2023PagedAttention,
  author    = {Woosuk Kwon and Zhuohan Li and Siyuan Zhuang and Ying Sheng and Lianmin Zheng and Cody Hao Yu and Joseph E. Gonzalez and Hao Zhang and Ion Stoica},
  title     = {Efficient Memory Management for Large Language Model Serving with {PagedAttention}},
  booktitle = {Proceedings of the 29th Symposium on Operating Systems Principles},
  pages     = {611--626},
  publisher = {Association for Computing Machinery},
  year      = {2023},
  doi       = {10.1145/3600006.3613165},
  url       = {https://doi.org/10.1145/3600006.3613165}
}

@article{Grattafiori2024Llama3,
  author  = {Aaron Grattafiori and others},
  title   = {The {Llama 3} Herd of Models},
  journal = {arXiv preprint arXiv:2407.21783},
  year    = {2024},
  doi     = {10.48550/arXiv.2407.21783},
  url     = {https://arxiv.org/abs/2407.21783}
}

@inproceedings{Agrawal2024SarathiServe,
  author    = {Amey Agrawal and Nitin Kedia and Ashish Panwar and Jayashree Mohan and Nipun Kwatra and Bhargav Gulavani and Alexey Tumanov and Ramachandran Ramjee},
  title     = {Taming {Throughput-Latency} Tradeoff in {LLM} Inference with {Sarathi-Serve}},
  booktitle = {18th USENIX Symposium on Operating Systems Design and Implementation (OSDI 24)},
  pages     = {117--134},
  publisher = {USENIX Association},
  year      = {2024},
  url       = {https://www.usenix.org/conference/osdi24/presentation/agrawal}
}

@misc{Ao2026WAIT,
  author    = {Ruicheng Ao and Gan Luo and David Simchi-Levi and Xinshang Wang},
  title     = {Optimizing {LLM} Inference: {Fluid-Guided} Online Scheduling with Memory Constraints},
  publisher = {arXiv},
  year      = {2026},
  doi       = {10.48550/arXiv.2504.11320},
  url       = {https://arxiv.org/abs/2504.11320}
}

@misc{Dai2026Throughput,
  author    = {J. G. Dai and Tianze Deng and Yueying Li and Tianyi Peng},
  title     = {Throughput-Optimal Scheduling Algorithms for {LLM} Inference and {AI} Agents},
  publisher = {arXiv},
  year      = {2026},
  doi       = {10.48550/arXiv.2504.07347},
  url       = {https://arxiv.org/abs/2504.07347}
}

@inproceedings{Nie2026Stability,
  author    = {Chengyi Nie and Nian Si and Zijie Zhou},
  title     = {A Queueing-Theoretic Framework for Stability Analysis of {LLM} Inference with {KV} Cache Memory Constraints},
  booktitle = {Proceedings of the 43rd International Conference on Machine Learning},
  year      = {2026},
  url       = {https://arxiv.org/abs/2605.04595}
}

@misc{Lin2026Heterogeneous,
  author    = {Ruihan Lin and Zezhen Ding and Zean Han and Jiheng Zhang},
  title     = {Large-Scale {LLM} Inference with Heterogeneous Workloads: Prefill--Decode Contention and Asymptotically Optimal Control},
  publisher = {arXiv},
  year      = {2026},
  doi       = {10.48550/arXiv.2602.02987},
  url       = {https://arxiv.org/abs/2602.02987}
}

@misc{Jaillet2025Online,
  author    = {Patrick Jaillet and Jiashuo Jiang and Konstantina Mellou and Marco Molinaro and Chara Podimata and Zijie Zhou},
  title     = {Online Scheduling for {LLM} Inference with {KV} Cache Constraints},
  publisher = {arXiv},
  year      = {2025},
  doi       = {10.48550/arXiv.2502.07115},
  url       = {https://arxiv.org/abs/2502.07115}
}

@misc{Wang2025Variable,
  author    = {Meixuan Wang and Yinyu Ye and Zijie Zhou},
  title     = {{LLM} Serving Optimization with Variable Prefill and Decode Lengths},
  publisher = {arXiv},
  year      = {2025},
  doi       = {10.48550/arXiv.2508.06133},
  url       = {https://arxiv.org/abs/2508.06133}
}

@misc{Ao2026Congestion,
  author    = {Ruicheng Ao and Jing Dong and Gan Luo and David Simchi-Levi},
  title     = {Service-Induced Congestion in Memory-Constrained {LLM} Serving},
  publisher = {arXiv},
  year      = {2026},
  doi       = {10.48550/arXiv.2606.15555},
  url       = {https://arxiv.org/abs/2606.15555}
}

@misc{Kong2026Geometry,
  author    = {Li Kong and Qi Qi and Yinyu Ye and Zijie Zhou},
  title     = {Geometry-Aware Online Scheduling for {LLM} Serving: From Theoretical Bound to System Practice},
  publisher = {arXiv},
  year      = {2026},
  doi       = {10.48550/arXiv.2606.22327},
  url       = {https://arxiv.org/abs/2606.22327}
}

@misc{Chen2025Adaptive,
  author    = {Zixi Chen and Yinyu Ye and Zijie Zhou},
  title     = {Adaptively Robust {LLM} Inference Optimization under Prediction Uncertainty},
  publisher = {arXiv},
  year      = {2025},
  doi       = {10.48550/arXiv.2508.14544},
  url       = {https://arxiv.org/abs/2508.14544}
}

@misc{Feng2026Competitive,
  author    = {Yiding Feng and Zonghan Yang and Yuhao Zhang},
  title     = {Competitive Non-Clairvoyant {KV-Cache} Scheduling for {LLM} Inference},
  publisher = {arXiv},
  year      = {2026},
  doi       = {10.48550/arXiv.2601.22996},
  url       = {https://arxiv.org/abs/2601.22996}
}

@misc{Feng2026General,
  author    = {Yiding Feng and Siyu Liu and Zonghan Yang and Yuhao Zhang},
  title     = {General Non-Clairvoyant {KV-Cache} Scheduling via Regime-Aware Routing},
  publisher = {arXiv},
  year      = {2026},
  doi       = {10.48550/arXiv.2607.09248},
  url       = {https://arxiv.org/abs/2607.09248}
}

@inproceedings{Wang2025KVCacheWild,
  author    = {Jiahao Wang and Jinbo Han and Xingda Wei and Sijie Shen and Dingyan Zhang and Chenguang Fang and Rong Chen and Wenyuan Yu and Haibo Chen},
  title     = {{KVCache} Cache in the Wild: Characterizing and Optimizing {KVCache} Cache at a Large Cloud Provider},
  booktitle = {2025 USENIX Annual Technical Conference (USENIX ATC 25)},
  year      = {2025},
  publisher = {USENIX Association},
  url       = {https://www.usenix.org/conference/atc25/presentation/wang-jiahao}
}

@article{coffman2001bandwidth,
  title={Bandwidth packing},
  author={Coffman, EG and Stolyar, Alexander L},
  journal={Algorithmica},
  volume={29},
  number={1},
  pages={70--88},
  year={2001},
  publisher={Springer}
}

@article{stolyar2015asymptotic,
  title={Asymptotic optimality of a greedy randomized algorithm in a large-scale service system with general packing constraints},
  author={Stolyar, Alexander L and Zhong, Yuan},
  journal={Queueing Systems},
  volume={79},
  number={2},
  pages={117--143},
  year={2015},
  publisher={Springer}
}

@article{stolyar2013large,
  title={A large-scale service system with packing constraints: Minimizing the number of occupied servers},
  author={Stolyar, Alexander L and Zhong, Yuan},
  journal={ACM SIGMETRICS Performance Evaluation Review},
  volume={41},
  number={1},
  pages={41--52},
  year={2013},
  publisher={ACM New York, NY, USA}
}

@article{gamarnik2004stochastic,
  title={Stochastic bandwidth packing process: stability conditions via Lyapunov function technique},
  author={Gamarnik, David},
  journal={Queueing systems},
  volume={48},
  number={3},
  pages={339--363},
  year={2004},
  publisher={Springer}
}

@article{stolyar2013infinite,
  title={An infinite server system with general packing constraints},
  author={Stolyar, Alexander L},
  journal={Operations Research},
  volume={61},
  number={5},
  pages={1200--1217},
  year={2013},
  publisher={INFORMS}
}

@article{ghaderi2014asymptotic,
  title={Asymptotic optimality of BestFit for stochastic bin packing},
  author={Ghaderi, Javad and Zhong, Yuan and Srikant, Rayadurgam},
  journal={ACM SIGMETRICS Performance Evaluation Review},
  volume={42},
  number={2},
  pages={64--66},
  year={2014},
  publisher={ACM New York, NY, USA}
}

@article{gupta2015lagrangian,
  title={Lagrangian-based online stochastic bin packing},
  author={Gupta, Varun and Radovanovic, Ana},
  journal={ACM SIGMETRICS Performance Evaluation Review},
  volume={43},
  number={1},
  pages={467--468},
  year={2015},
  publisher={ACM New York, NY, USA}
}

@inproceedings{maguluri2012stochastic,
  title={Stochastic models of load balancing and scheduling in cloud computing clusters},
  author={Maguluri, Siva Theja and Srikant, Rayadurgam and Ying, Lei},
  booktitle={2012 Proceedings IEEE Infocom},
  pages={702--710},
  year={2012},
  organization={IEEE}
}

@article{huang2025orlm,
  title={Orlm: A customizable framework in training large models for automated optimization modeling},
  author={Huang, Chenyu and Tang, Zhengyang and Hu, Shixi and Jiang, Ruoqing and Zheng, Xin and Ge, Dongdong and Wang, Benyou and Wang, Zizhuo},
  journal={Operations Research},
  year={2025},
  publisher={INFORMS}
}

@article{yao2026multi,
  title={Multi-Turn LLM Conversations under the Least-Recently-Used Policy: Mean-Field Asymptotics and Hit Ratio Approximation},
  author={Yao, Heyuan and Gao, Chutong and Lyu, Yuan and Grosof, Izzy and Simchi-Levi, David},
  journal={arXiv preprint arXiv:2609.02027},
  year={2026}
}

@article{brown2020language,
  title={Language models are few-shot learners},
  author={Brown, Tom and Mann, Benjamin and Ryder, Nick and Subbiah, Melanie and Kaplan, Jared D and Dhariwal, Prafulla and Neelakantan, Arvind and Shyam, Pranav and Sastry, Girish and Askell, Amanda and others},
  journal={Advances in neural information processing systems},
  volume={33},
  pages={1877--1901},
  year={2020}
}

@article{openai2023gpt,
  title={{GPT}-4 technical report. arxiv 2303.08774},
  author={OpenAI},
  journal={View in Article},
  volume={2},
  number={5},
  year={2023}
}

@misc{chatgpt2023,
  author = {OpenAI},
  title = {{ChatGPT}},
  howpublished = {\url{https://chat.openai.com}},
  year = {2019}
}

@misc{githubcopilot2023,
  author = {{GitHub}},
  title = {{GitHub} Copilot},
  howpublished = {\url{https://github.com/features/copilot}},
  year = {2021}
}

@misc{chatgpt_users_2023,
  author       = {GilPress},
  title        = {{ChatGPT} Users Statistics},
  howpublished = {WhatsTheBigData},
  year         = {2024},
  url          = {https://whatsthebigdata.com/chatgpt-users/#google_vignette}
}

@misc{gordon_2024,
  author       = {Gordon, Cindy},
  title        = {{ChatGPT} and {Generative AI} Innovations Are Creating Sustainability Havoc},
  howpublished = {Forbes},
  year         = {2024},
  url          = {https://www.forbes.com/sites/cindygordon/2024/03/12/chatgpt-and-generative-ai-innovations-are-creating-sustainability-havoc/}
}

@article{xu2023drift,
  title={Drift method: from stochastic networks to machine learning},
  author={Xu, Kuang},
  journal={URL: https://web. stanford. edu/kuangxu/papers/driftmethod. pdf. Last visited on},
  volume={3},
  number={09},
  year={2023}
}

@misc{nvidia2026rubin,
  author       = {Eduardo Alvarez and Vishal Mehta and Farshad Ghodsian},
  title        = {Inside NVIDIA Rubin GPU Architecture: Powering the Era of Agentic AI},
  year         = {2026},
  month        = jul,
  howpublished = {NVIDIA Technical Blog},
  url          = {https://developer.nvidia.com/blog/inside-nvidia-rubin-gpu-architecture-powering-the-era-of-agentic-ai/},
  note         = {Accessed: 2026-09-20}
}

@article{pope2023efficiently,
  title={Efficiently scaling transformer inference},
  author={Pope, Reiner and Douglas, Sholto and Chowdhery, Aakanksha and Devlin, Jacob and Bradbury, James and Heek, Jonathan and Xiao, Kefan and Agrawal, Shivani and Dean, Jeff},
  journal={Proceedings of Machine Learning and Systems},
  volume={5},
  pages={606--624},
  year={2023}
}

@article{singh2025openai,
  title={Openai gpt-5 system card},
  author={Singh, Aaditya and Fry, Adam and Perelman, Adam and Tart, Adam and Ganesh, Adi and El-Kishky, Ahmed and McLaughlin, Aidan and Low, Aiden and Ostrow, AJ and Ananthram, Akhila and others},
  journal={arXiv preprint arXiv:2601.03267},
  year={2025}
}

@article{vaswani2017attention,
  title={Attention is all you need},
  author={Vaswani, Ashish and Shazeer, Noam and Parmar, Niki and Uszkoreit, Jakob and Jones, Llion and Gomez, Aidan N and Kaiser, {\L}ukasz and Polosukhin, Illia},
  journal={Advances in neural information processing systems},
  volume={30},
  year={2017}
}

@article{Buyukkoc1985Cmu,
  author  = {Buyukkoc, Cagatay and Varaiya, Pravin and Walrand, Jean},
  title   = {The {$c\mu$} Rule Revisited},
  journal = {Advances in Applied Probability},
  volume  = {17},
  number  = {1},
  pages   = {237--238},
  year    = {1985},
  doi     = {10.2307/1427064},
  url     = {https://doi.org/10.2307/1427064}
}

@article{VanMieghem1995Generalized,
  author  = {Van Mieghem, Jan A.},
  title   = {Dynamic Scheduling with Convex Delay Costs: The Generalized {$c\mu$} Rule},
  journal = {The Annals of Applied Probability},
  volume  = {5},
  number  = {3},
  pages   = {809--833},
  year    = {1995},
  doi     = {10.1214/aoap/1177004706},
  url     = {https://doi.org/10.1214/aoap/1177004706}
}

@article{Mandelbaum2004Generalized,
  author  = {Mandelbaum, Avishai and Stolyar, Alexander L.},
  title   = {Scheduling Flexible Servers with Convex Delay Costs: Heavy-Traffic Optimality of the Generalized {$c\mu$}-Rule},
  journal = {Operations Research},
  volume  = {52},
  number  = {6},
  pages   = {836--855},
  year    = {2004},
  doi     = {10.1287/opre.1040.0152},
  url     = {https://doi.org/10.1287/opre.1040.0152}
}

@article{Zychlinski2023Resources,
  author  = {Zychlinski, Noa and Chan, Carri W. and Dong, Jing},
  title   = {Managing Queues with Different Resource Requirements},
  journal = {Operations Research},
  volume  = {71},
  number  = {4},
  pages   = {1387--1413},
  year    = {2023},
  doi     = {10.1287/opre.2022.2284},
  url     = {https://doi.org/10.1287/opre.2022.2284}
}

@article{Bramson1998StateSpaceCollapse,
  author  = {Bramson, Maury},
  title   = {State Space Collapse with Application to Heavy Traffic Limits for Multiclass Queueing Networks},
  journal = {Queueing Systems},
  volume  = {30},
  number  = {1--2},
  pages   = {89--148},
  year    = {1998},
  doi     = {10.1023/A:1019160803783},
  url     = {https://doi.org/10.1023/A:1019160803783}
}

@article{Williams1998DiffusionApproximations,
  author  = {Williams, Ruth J.},
  title   = {Diffusion Approximations for Open Multiclass Queueing Networks: Sufficient Conditions Involving State Space Collapse},
  journal = {Queueing Systems},
  volume  = {30},
  number  = {1--2},
  pages   = {27--88},
  year    = {1998},
  doi     = {10.1023/A:1019108819713},
  url     = {https://doi.org/10.1023/A:1019108819713}
}

\clearpage
\appendix
\section*{Appendices}
\section{Non-Anticipativity of Prelimit Policies}
\label{sec:prelimit-nonanticipativity}

Fix $r\geq1$ and an initial state $z_0^r$. The arrival processes are the only source of randomness in the system. Define their filtration by
\[
\cF_\tau^r:=\sigma\!\left(N_j^r\!\left(r\int_0^s\lambda_j(u)\,du\right):j=1,\dots,J,\ 0\leq s\leq\tau\right),\qquad \tau\geq0,
\]
with the usual augmentation. A prelimit policy is non-anticipative if, at every GPU step $t$, all control increments $\Delta\mathcal U_{j,k}^r(t)$ and $\Delta\mathcal V_{j,k}^r(t)$ are measurable with respect to $\cG_t^r:=\cF_{f^r(t)}^r$.

We show by induction that $f^r(t)$ is a stopping time with respect to $\{\cF_\tau^r:\tau\geq0\}$ for every $t\in\mathbb Z_+$. First, $f^r(0)=0$ is a stopping time.

Suppose that $f^r(n)$ is a stopping time for every $n=0,\dots,t$. By the non-anticipativity condition, the control increments $\Delta\mathcal U_{j,k}^r(t)$ and $\Delta\mathcal V_{j,k}^r(t)$, over their respective prefill and decode stages, are $\cG_t^r$-measurable. The step duration $\Delta^r(t)$ is a function of these increments, so $f^r(t+1)=f^r(t)+\Delta^r(t)$ is also $\cG_t^r$-measurable.

For every $\tau\geq0$, the event $\{\Delta^r(t)\leq\tau-f^r(t)\}$ belongs to $\cG_t^r=\cF_{f^r(t)}^r$. By the definition of $\cF_{f^r(t)}^r$, intersecting this event with $\{f^r(t)\leq\tau\}$ gives an event in $\cF_\tau^r$. Since $\Delta^r(t)\geq0$, we obtain
\[
\{f^r(t+1)\leq\tau\}=\{f^r(t)\leq\tau,\ \Delta^r(t)\leq\tau-f^r(t)\}\in\cF_\tau^r.
\]
Thus, $f^r(t+1)$ is a stopping time, completing the induction. Consequently, $\cG_t^r=\cF_{f^r(t)}^r$ is well-defined for every $t$.

\section{Illustration of the Fluid \texorpdfstring{$c\mu$}{c-mu} Policy under General Arrivals}
\label{sec:general-arrival-policy-illustration}
\begin{figure}[htbp]
\centering
\includegraphics{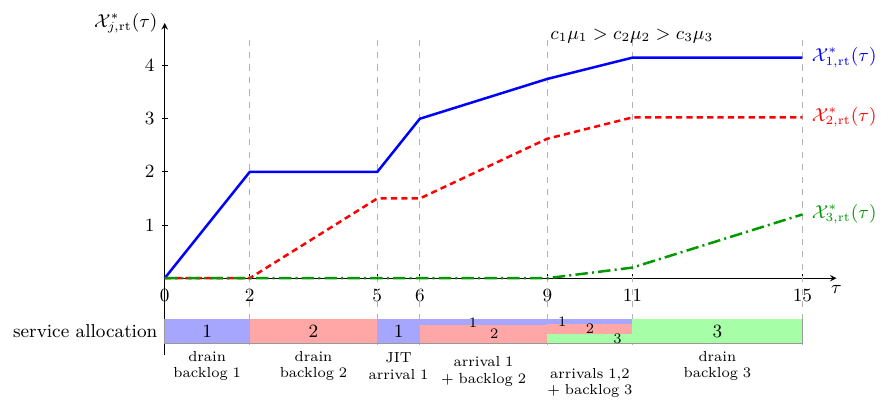}
\caption{Real-time representation of the optimal $c\mu$ policy for three priority classes under general arrivals. The curves show the cumulative class-level completion processes, and the lower bar shows the corresponding service-effort allocation.}
\label{fig:general-arrival-cmu-policy}
\end{figure}

Figure~\ref{fig:general-arrival-cmu-policy} illustrates the policy for three classes with $c_1\mu_1>c_2\mu_2>c_3\mu_3$. Here, $\mathcal X_{j,\mathrm{rt}}^*(\tau)$ denotes the cumulative number of class-$j$ completions by clock time $\tau$. During $[0,2]$, the class-1 backlog receives all available service, while lower-priority classes wait. During $[6,9]$, class-1 requests are served as they arrive, using only part of the service budget; the remainder is used to serve the class-2 backlog. During $[9,11]$, the policy keeps up with arrivals to both classes 1 and 2 and allocates the remaining budget to class 3.

\section{Proof of Theorem~\ref{thm:fluid-cmu-optimality}}
\label{sec:proof-fluid-cmu-optimality}

We first introduce the workload processes used to derive a lower bound on the cumulative holding cost.

\begin{definition}[Workload processes]
For an admissible policy $\pi$ and each priority prefix $\{1,\dots,\ell\}$, define the prefix workload at step time $\xi$ by
\begin{equation}
W_\ell^\pi(\xi):=\sum_{j=1}^\ell\mu_j^{-1}p_{j,L_j}^{\pi}(\xi), \qquad \ell=1,\dots,J. \label{eq:prefix-workload-W}
\end{equation}
In clock time $\tau$, define the net-input process $g_\ell$ and the corresponding reflected workload process $\varphi_\ell$ by
\begin{equation*}
g_\ell(\tau):=\sum_{j=1}^{\ell}\mu_j^{-1}A_j(\tau)-\tau,\qquad \varphi_\ell(\tau):=g_\ell(\tau)-\inf_{0\le s\le\tau}(g_\ell(s)\wedge0).
\end{equation*}
\end{definition}

\begin{remark}[Workload interpretation]
The queue length $p_{j,L_j}^{\pi}$ records the number of waiting class-$j$ requests. We use the effective service rate $\mu_j$ to convert this quantity into real-time units, so $W_\ell^\pi$ measures the total effective workload of waiting requests in the first $\ell$ classes. Applying the same conversion to cumulative arrivals, $g_\ell$ describes the net workload input to these classes. Since unused service capacity cannot be saved for future arrivals, the Skorokhod reflection $\varphi_\ell$ gives the remaining workload when the full service capacity is devoted to these classes whenever they have work waiting. This yields the smallest remaining workload permitted by the available capacity, and hence a lower bound on $W_\ell^\pi$ at the corresponding clock time.
\end{remark}

\begin{proof}
First, we establish a lower bound on $\mathcal J(\pi;z_0)$. To obtain this bound, we express the integrand in $\mathcal J(\pi;z_0)$ as a nonnegative weighted sum of the prefix workloads. Setting $c_{J+1}\mu_{J+1}:=0$, we obtain
\begin{equation}
\begin{aligned}
\sum_{\ell=1}^J \left(c_\ell\mu_\ell - c_{\ell+1}\mu_{\ell+1}\right)W_\ell^\pi(\xi)&=\sum_{\ell=1}^J \left(c_\ell\mu_\ell - c_{\ell+1}\mu_{\ell+1}\right)\sum_{j=1}^\ell\mu_j^{-1}p_{j,L_j}^{\pi}(\xi)\\
&=\sum_{j=1}^{J}\mu_j^{-1}p_{j,L_j}^{\pi}(\xi) \sum_{\ell=j}^J\left(c_\ell\mu_\ell - c_{\ell+1}\mu_{\ell+1}\right)\\
&=\sum_{j=1}^{J}c_jp_{j,L_j}^{\pi}(\xi).
\end{aligned}
\label{eq:telescoping-cost-W}
\end{equation}
It remains to show that $\varphi_\ell$ provides a lower bound on the prefix workload $W_\ell^\pi$. Using the queue balance $p_{j,L_j}^{\pi}(\xi)=A_j(f^\pi(\xi))-\mathcal X_j^\pi(\xi)$, we first write
\begin{equation}
W_\ell^\pi(\xi)=\sum_{j=1}^{\ell}\mu_j^{-1}A_j(f^\pi(\xi))-\sum_{j=1}^{\ell}\mu_j^{-1}\mathcal X_j^\pi(\xi). \label{eq:prefix-W-flow-balance}
\end{equation}

We next establish two inequalities for $\{\mathcal X_j^\pi\}_{j=1}^{J}$. The memory constraint in \eqref{eq:class-level-memory} gives
\begin{equation*}
d\xi\geq\frac{1}{M}\sum_{j=1}^{J}\tilde m_j\,d\mathcal X_j^\pi(\xi).
\end{equation*}
Together with the definitions of $f^\pi$ and $\mu_j$, this yields
\begin{equation}
\begin{aligned}
df^\pi(\xi)&=\sum_{j=1}^{J}(a_{\mathrm{pf}}s_j+a_{\mathrm{kv}}\tilde m_j)\,d\mathcal X_j^\pi(\xi)+b\,d\xi\\
&\geq\sum_{j=1}^{J}\left[a_{\mathrm{pf}}s_j+\left(a_{\mathrm{kv}}+\frac{b}{M}\right)\tilde m_j\right]d\mathcal X_j^\pi(\xi)\\
&=\sum_{j=1}^{J}\mu_j^{-1}\,d\mathcal X_j^\pi(\xi).
\end{aligned}
\label{eq:clock-mu-dX-bound}
\end{equation}
Since $\mathcal X_j^\pi$ is nondecreasing and $\mathcal X_j^\pi(0)=f^\pi(0)=0$, integrating over $[0,\xi]$ and $[\xi^\prime,\xi]$, respectively, gives, for $0\leq\xi^\prime\leq\xi$,
\begin{equation}
\sum_{j=1}^{\ell}\mu_j^{-1}\mathcal X_j^\pi(\xi)\leq f^\pi(\xi),\quad\text{and}\quad\sum_{j=1}^{\ell}\mu_j^{-1}\bigl[\mathcal X_j^\pi(\xi)-\mathcal X_j^\pi(\xi^\prime)\bigr]\leq f^\pi(\xi)-f^\pi(\xi^\prime). \label{eq:partial-work-bound}
\end{equation}

Substituting the first inequality into \eqref{eq:prefix-W-flow-balance}, we obtain
\begin{equation}
W_\ell^\pi(\xi)\geq\sum_{j=1}^{\ell}\mu_j^{-1}A_j(f^\pi(\xi))-f^\pi(\xi)=g_\ell(f^\pi(\xi)). \label{eq:W-lower-g}
\end{equation}
Substituting the second inequality and using $W_\ell^\pi(\xi^\prime)\geq0$, we obtain
\begin{equation}
\begin{aligned}
W_\ell^\pi(\xi)&\geq\sum_{j=1}^{\ell}\mu_j^{-1}\bigl[A_j(f^\pi(\xi))-\mathcal X_j^\pi(\xi^\prime)\bigr]-f^\pi(\xi)+f^\pi(\xi^\prime)\\
&=W_\ell^\pi(\xi^\prime)+g_\ell(f^\pi(\xi))-g_\ell(f^\pi(\xi^\prime))\\
&\geq g_\ell(f^\pi(\xi))-g_\ell(f^\pi(\xi^\prime)).
\end{aligned}
\label{eq:W-lower-g-increment}
\end{equation}

Combining these two workload bounds and taking the supremum over $0\leq\xi^\prime\leq\xi$, we obtain
\begin{equation}
\begin{aligned}
W_\ell^\pi(\xi)&\geq g_\ell(f^\pi(\xi))-\inf_{0\leq\xi^\prime\leq\xi}\bigl(g_\ell(f^\pi(\xi^\prime))\wedge0\bigr)\\
&=g_\ell(f^\pi(\xi))-\inf_{0\leq s\leq f^\pi(\xi)}\bigl(g_\ell(s)\wedge0\bigr)\\
&=\varphi_\ell(f^\pi(\xi)),
\end{aligned}
\label{eq:W-lower-varphi}
\end{equation}
where the change of variable $s=f^\pi(\xi^\prime)$ uses the continuity and monotonicity of $f^\pi$, together with $f^\pi(0)=0$.
The clock satisfies $f^\pi(\xi)\geq b\xi\to\infty$ because $b>0$. Thus, changing variables to $\tau=f^\pi(\xi)$ preserves the integration limits $0$ and $\infty$. Finally, using \eqref{eq:telescoping-cost-W}, \eqref{eq:W-lower-varphi}, and $c_\ell\mu_\ell-c_{\ell+1}\mu_{\ell+1}\geq 0$, we obtain
\begin{equation*}
\begin{aligned}
\mathcal J(\pi;z_0)&= \int_0^\infty \sum_{j=1}^J c_jp_{j,L_j}^{\pi}(\xi)\,df^\pi(\xi)\\
&= \int_0^\infty\sum_{\ell=1}^J \left(c_\ell\mu_\ell - c_{\ell+1}\mu_{\ell+1}\right)W_\ell^\pi(\xi)\,df^\pi(\xi)\\
&\geq \int_0^\infty\sum_{\ell=1}^J \left(c_\ell\mu_\ell - c_{\ell+1}\mu_{\ell+1}\right)\varphi_\ell(f^\pi(\xi))\,df^\pi(\xi)\\
&= \int_{f^\pi(0)}^{\lim_{\xi\rightarrow\infty}f^\pi(\xi)} \sum_{\ell=1}^J \left(c_\ell\mu_\ell - c_{\ell+1}\mu_{\ell+1}\right)\varphi_\ell(\tau)\,d\tau\\
&= \int_0^\infty \sum_{\ell=1}^J \left(c_\ell\mu_\ell - c_{\ell+1}\mu_{\ell+1}\right)\varphi_\ell(\tau)\,d\tau.
\end{aligned}
\end{equation*}

Next, we prove that $\pi^*$ is admissible by verifying the three conditions in Definition~\ref{def:admissible-fluid-policy}.

We first show that $F$ is a continuous strictly increasing bijection from $[0,\infty)$ onto itself, so that $f^*=F^{-1}$ is well defined.
Since $R_0^*$ and $A_j$ are continuous and the operations $[\cdot]^+$ and $\sup_{0\leq s\leq\tau}$ preserve continuity, each $R_j^*$ defined recursively by \eqref{eq:target-residual} is continuous. Each $R_j^*$ is also nondecreasing because it is defined using $\sup_{0\leq s\leq\tau}$. The recursion also gives $R_j^*(0)=0$. Therefore, by the definition of $F$ in \eqref{eq:target-F}, $F$ is continuous with $F(0)=0$.

To establish strict monotonicity and unboundedness of $F$, we use the following lemma.
\begin{lemma}\label{lem:R-increment}
Let $0\leq\tau_1\leq\tau_2$. Then, for $j=1,\dots,J$,
\begin{equation*}
R_j^*(\tau_2) - R_j^*(\tau_1) \leq R_{j-1}^*(\tau_2) - R_{j-1}^*(\tau_1).
\end{equation*}
\end{lemma}
For $0\leq\tau_1<\tau_2$, by the definition of $F$ in \eqref{eq:target-F}, we have
\begin{equation*}
{\medmuskip=1.25mu\thickmuskip=1.75mu\thinmuskip=1.25mu\scriptspace=0pt F(\tau_2)-F(\tau_1)=\sum_{j=1}^{J}\frac{\tilde m_j\mu_j}{M}\bigl[R_{j-1}^*(\tau_2)-R_j^*(\tau_2)-R_{j-1}^*(\tau_1)+R_j^*(\tau_1)\bigr]+\frac{1}{b}\bigl[R_J^*(\tau_2)-R_J^*(\tau_1)\bigr].}
\end{equation*}
By Lemma~\ref{lem:R-increment} and the monotonicity of $R_j^*$, all increments on the right-hand side are nonnegative. Since $b>0$, $M>0$, $\tilde m_j>0$, and $\mu_j>0$, we can set
\begin{equation*}
\kappa:=\min\left\{\frac{1}{b},\frac{\tilde m_1\mu_1}{M},\dots,\frac{\tilde m_J\mu_J}{M}\right\}>0.
\end{equation*}
We can therefore bound each coefficient below by $\kappa$. By telescoping and using $R_0^*(\tau)=\tau$, we obtain
\begin{equation*}
{\medmuskip=1.25mu\thickmuskip=1.75mu\thinmuskip=1.25mu\scriptspace=0pt
\begin{aligned}
F(\tau_2)-F(\tau_1)&\geq\kappa\Biggl(\sum_{j=1}^{J}\bigl[R_{j-1}^*(\tau_2)-R_j^*(\tau_2)-R_{j-1}^*(\tau_1)+R_j^*(\tau_1)\bigr]+R_J^*(\tau_2)-R_J^*(\tau_1)\Biggr)\\
&=\kappa\bigl[R_0^*(\tau_2)-R_0^*(\tau_1)\bigr]=\kappa(\tau_2-\tau_1)>0.
\end{aligned}}
\end{equation*}
This proves that $F$ is strictly increasing. Taking $\tau_1=0$ and using $F(0)=0$, we also obtain $F(\tau)\geq\kappa\tau\to\infty$. Together with continuity, these properties show that $F$ maps $[0,\infty)$ bijectively onto itself. Consequently, $f^*$ in \eqref{eq:target-clock} is well defined and satisfies $F(f^*(\xi))=\xi$.

We next verify that each $\mathcal X_j^*$ is continuous and nondecreasing, with $\mathcal X_j^*(0)=0$. Since $f^*$ and all $R_j^*$ are continuous, $\mathcal X_j^*$ defined in \eqref{eq:target-admission} is continuous. Since $f^*(0)=0$ and $R_j^*(0)=0$ for all $j$, \eqref{eq:target-admission} also gives $\mathcal X_j^*(0)=0$. For $0\leq\xi_1\leq\xi_2$, using $f^*(\xi_1)\leq f^*(\xi_2)$ and Lemma~\ref{lem:R-increment}, we obtain
\begin{equation*}
\mathcal X_j^*(\xi_2)-\mathcal X_j^*(\xi_1) =\mu_j\Bigl(R_{j-1}^*\bigl(f^*(\xi_2)\bigr)-R_{j-1}^*\bigl(f^*(\xi_1)\bigr)\Bigr)-\mu_j\Bigl(R_j^*\bigl(f^*(\xi_2)\bigr)-R_j^*\bigl(f^*(\xi_1)\bigr)\Bigr)\geq0.
\end{equation*}
We next verify the state constraints. For the intermediate-stage states, substituting $\mathcal U_{j,k}^*=\mathcal V_{j,k}^*=\mathcal X_j^*$ from Definition~\ref{def:fluid-cmu-policy} into \eqref{eq:fluid-stage-flow-balance}, we obtain
\begin{equation*}
p_{j,k}^{\pi^*}(\xi)=p_{j,k}(0)+\mathcal X_j^*(\xi)-\mathcal X_j^*(\xi)=0,\quad\text{and}\quad d_{j,k}^{\pi^*}(\xi)=d_{j,k}(0)+\mathcal X_j^*(\xi)-\mathcal X_j^*(\xi)=0.
\end{equation*}
It remains to show that the waiting-queue lengths satisfy $p_{j,L_j}^{\pi^*}(\xi)\geq 0$.
\begin{equation*}
{\medmuskip=1.5mu\thickmuskip=2mu p_{j,L_j}^{\pi^*}(\xi)=A_j(f^*(\xi))-\mathcal X_j^*(\xi)=\mu_j\left(\mu_j^{-1}A_j(f^*(\xi))-R_{j-1}^*(f^*(\xi))+R_j^*(f^*(\xi))\right).}
\end{equation*}
Since
\begin{equation*}
\begin{aligned}
R_j^*(f^*(\xi)) &= \sup_{0\leq s\leq f^*(\xi)} \left[R_{j-1}^*(s)-\mu_j^{-1}A_j(s)\right]^+\\
&\geq \left[R_{j-1}^*(f^*(\xi))-\mu_j^{-1}A_j(f^*(\xi))\right]^+\\
&\geq R_{j-1}^*(f^*(\xi))-\mu_j^{-1}A_j(f^*(\xi)),
\end{aligned}
\end{equation*}
we have
$p_{j,L_j}^{\pi^*}(\xi)\geq 0$.
We next verify the memory constraint \eqref{eq:class-level-memory}. Using the identity $F(f^*(\xi))=\xi$, for each $\xi\geq 0$,
\begin{equation}
\begin{aligned}
\xi=F(f^*(\xi))&=\frac{1}{M}\sum_{j=1}^{J}\tilde m_j\mu_j \left(R_{j-1}^*(f^*(\xi))-R_j^*(f^*(\xi))\right)+\frac{1}{b}R_J^*(f^*(\xi))\\
&=\frac{1}{M}\sum_{j=1}^{J}\tilde m_j\mathcal X_j^*(\xi) +\frac{1}{b}R_J^*(f^*(\xi)).
\end{aligned}
\label{eq:target-memory-identity}
\end{equation}
Thus, for $0\leq \xi_1\leq \xi_2$,
\begin{equation*}
\xi_2-\xi_1 =\frac{1}{M}\sum_{j=1}^{J}\tilde m_j \left(\mathcal X_j^*(\xi_2)-\mathcal X_j^*(\xi_1)\right)+\frac{1}{b} \left(R_J^*(f^*(\xi_2))-R_J^*(f^*(\xi_1))\right).
\end{equation*}
Equivalently,
\begin{equation*}
\sum_{j=1}^{J}\tilde m_j\left(\mathcal X_j^*(\xi_2)-\mathcal X_j^*(\xi_1)\right) =M(\xi_2-\xi_1)-\frac{M}{b}\left(R_J^*(f^*(\xi_2))-R_J^*(f^*(\xi_1))\right).
\end{equation*}
Since $f^*$ and $R_J^*$ are nondecreasing, $R_J^*(f^*(\xi_2))-R_J^*(f^*(\xi_1))\geq 0$.
Therefore, for $0\leq\xi_1\leq\xi_2$,
\begin{equation*}
\sum_{j=1}^{J}\tilde m_j \left(\mathcal X_j^*(\xi_2)-\mathcal X_j^*(\xi_1)\right) \leq M(\xi_2-\xi_1),
\end{equation*}
which is exactly the memory constraint.

Finally, we verify that $f^*$ satisfies the clock equation \eqref{eq:class-level-clock}. By telescoping and the definition of $\mathcal X_j^*$ in \eqref{eq:target-admission}, we obtain
\[
\begin{array}{@{}r@{\;}c@{\;}l@{}}
\displaystyle f^*(\xi)&=&\displaystyle\sum_{j=1}^{J}\left(R_{j-1}^*(f^*(\xi))-R_j^*(f^*(\xi))\right)+R_J^*(f^*(\xi))\\
&=&\displaystyle\sum_{j=1}^{J}\mu_j^{-1}\mathcal X_j^*(\xi)+R_J^*(f^*(\xi))\\
&\overset{\mathrm{(a)}}{=}&\displaystyle\sum_{j=1}^{J}\mu_j^{-1}\mathcal X_j^*(\xi)+b\xi-\frac{b}{M}\sum_{j=1}^{J}\tilde m_j\mathcal X_j^*(\xi)\\
&\overset{\mathrm{(b)}}{=}&\displaystyle\sum_{j=1}^{J}\left(a_{\mathrm{pf}}s_j+a_{\mathrm{kv}}\tilde m_j\right)\mathcal X_j^*(\xi)+b\xi,
\end{array}
\]
where (a) follows from \eqref{eq:target-memory-identity} and (b) uses $\mu_j^{-1}=a_{\mathrm{pf}}s_j+\left(a_{\mathrm{kv}}+\frac{b}{M}\right)\tilde m_j$.
Therefore, $\{\mathcal X_j^*\}_{j=1}^{J}$ satisfies all the required conditions. By Section~\ref{subsec:stage-level-collapse}, the corresponding policy $\pi^*=(\mathcal U^*,\mathcal V^*)$ is admissible.

It remains to prove that $\pi^*$ attains the lower bound. The following lemma expresses $R_j^*$ in terms of the net-input process $g_j$.

\begin{lemma}\label{lem:R-reflection}
For $j=1,\dots,J$ and $\tau\geq 0$,
\begin{equation*}
R_j^*(\tau) = - \inf_{0\leq s\leq \tau}(g_j(s)\wedge 0).
\end{equation*}
\end{lemma}

We now prove $W_j^{\pi^*}(\xi)=\varphi_j(f^*(\xi))$. Starting from \eqref{eq:prefix-W-flow-balance} and using Lemma~\ref{lem:R-reflection}, we obtain
\begin{equation}
{\medmuskip=2mu\thickmuskip=2.5mu\thinmuskip=2mu\scriptspace=0pt
\begin{aligned}
W_j^{\pi^*}(\xi)&=\sum_{r=1}^j\mu_r^{-1}A_r(f^*(\xi))-\sum_{r=1}^j\mu_r^{-1}\mathcal X_r^*(\xi)\\
&=\sum_{r=1}^j\Bigl(\mu_r^{-1}A_r(f^*(\xi))-R_{r-1}^*(f^*(\xi))+R_r^*(f^*(\xi))\Bigr)\\
&= \sum_{r=1}^j\mu_r^{-1}A_r(f^*(\xi)) -f^*(\xi) +R_j^*(f^*(\xi))\\
&= g_j(f^*(\xi)) - \inf_{0\leq s\leq f^*(\xi)}(g_j(s)\wedge 0) = \varphi_j(f^*(\xi)).
\end{aligned}}
\label{eq:prefix-workload-attainment}
\end{equation}

Finally, we compute the cost under this policy:
\begin{equation}
\begin{aligned}
\mathcal J(\pi^*;z_0)&= \int_0^\infty \sum_{j=1}^J c_jp_{j,L_j}^{\pi^*}(\xi)\,df^*(\xi)\\
&= \int_0^\infty\sum_{j=1}^J \left(c_j\mu_j - c_{j+1}\mu_{j+1}\right) W_j^{\pi^*}(\xi)\,df^*(\xi)\\
&= \int_0^\infty\sum_{j=1}^J \left(c_j\mu_j - c_{j+1}\mu_{j+1}\right) \varphi_j(f^*(\xi))\,df^*(\xi)\\
&= \int_0^\infty\sum_{j=1}^J \left(c_j\mu_j - c_{j+1}\mu_{j+1}\right) \varphi_j(\tau)\,d\tau.
\end{aligned}
\label{eq:fluid-optimal-value}
\end{equation}
The policy therefore attains the lower bound derived above, proving Theorem~\ref{thm:fluid-cmu-optimality}.
\end{proof}

\subsection{Proof of Lemma~\ref{lem:R-increment}}
\begin{proof}
By the definition of $R_j^*$ and the monotonicity of $R_{j-1}^*$, for every $s\in[0,\tau_1]$,
\begin{equation}
\begin{array}{@{}r@{\;}c@{\;}l@{}}
\displaystyle\left[R_{j-1}^*(s)-\mu_j^{-1}A_j(s)\right]^+&\leq&\displaystyle\sup_{0\leq u\leq\tau_1}\left[R_{j-1}^*(u)-\mu_j^{-1}A_j(u)\right]^+\\
&=&\displaystyle R_j^*(\tau_1)\leq R_j^*(\tau_1)+\left[R_{j-1}^*(\tau_2)-R_{j-1}^*(\tau_1)\right].
\end{array}
\label{eq:R-pointwise-bound}
\end{equation}
For $s\in(\tau_1,\tau_2]$, since $R_{j-1}^*$ and $A_j$ are nondecreasing and $\mu_j>0$, we have
{\medmuskip=1mu\thickmuskip=1.5mu\thinmuskip=1mu\scriptspace=0pt
\begin{equation*}
\left[R_{j-1}^*(s)-\mu_j^{-1}A_j(s)\right]^+\leq\left[R_{j-1}^*(\tau_2)-\mu_j^{-1}A_j(\tau_1)\right]^+.
\end{equation*}}
Using $[x+y]^+\leq x+[y]^+$ for $x\geq0$, we obtain
\begin{equation*}
{\medmuskip=1mu\thickmuskip=2mu\thinmuskip=2mu\scriptspace=0pt \begin{aligned}
\left[R_{j-1}^*(\tau_2)-\mu_j^{-1}A_j(\tau_1)\right]^+&= \left[R_{j-1}^*(\tau_2)-R_{j-1}^*(\tau_1)+R_{j-1}^*(\tau_1)-\mu_j^{-1}A_j(\tau_1)\right]^+\\
&\leq \bigl(R_{j-1}^*(\tau_2)-R_{j-1}^*(\tau_1)\bigr)+\left[R_{j-1}^*(\tau_1)-\mu_j^{-1}A_j(\tau_1)\right]^+\\
&\leq R_j^*(\tau_1)+\bigl(R_{j-1}^*(\tau_2)-R_{j-1}^*(\tau_1)\bigr),
\end{aligned}}
\end{equation*}
where the last inequality follows from \eqref{eq:R-pointwise-bound}.
Therefore, for $s\in[0,\tau_2]$,
\begin{equation*}
\left[R_{j-1}^*(s)-\mu_j^{-1}A_j(s)\right]^+\leq R_j^*(\tau_1)+\left[R_{j-1}^*(\tau_2)-R_{j-1}^*(\tau_1)\right].
\end{equation*}
Taking the supremum of both sides over $s\in[0,\tau_2]$, we obtain
\begin{equation*}
R_j^*(\tau_2)=\sup_{0\leq s\leq\tau_2}\left[R_{j-1}^*(s)-\mu_j^{-1}A_j(s)\right]^+\leq R_j^*(\tau_1)+\left[R_{j-1}^*(\tau_2)-R_{j-1}^*(\tau_1)\right].
\end{equation*}
Therefore, $R_j^*(\tau_2)-R_j^*(\tau_1)\leq R_{j-1}^*(\tau_2)-R_{j-1}^*(\tau_1)$.
\end{proof}

\subsection{Proof of Lemma~\ref{lem:R-reflection}}
\begin{proof}
We first prove, by induction on $j$, that
\begin{equation*}
R_j^*(\tau) = \sup_{0\leq u\leq \tau} \left[ u-\sum_{r=1}^{j}\mu_r^{-1}A_r(u) \right]^+,
\end{equation*}
where the sum is understood as $0$ when $j=0$. For $j=0$, $R_0^*(\tau)=\sup_{0\leq u\leq\tau}[u]^+$.
Suppose the claim holds for $j-1$. Then
\begin{equation*}
{\medmuskip=2mu\thickmuskip=2.5mu R_j^*(\tau)= \sup_{0\leq s\leq\tau}\left[R_{j-1}^*(s)-\mu_j^{-1}A_j(s)\right]^+ = \sup_{0\leq s\leq\tau}\left[\sup_{0\leq u\leq s}\left[u-\sum_{r=1}^{j-1}\mu_r^{-1}A_r(u)\right]^+-\mu_j^{-1}A_j(s)\right]^+.}
\end{equation*}
For $0\leq u\leq s\leq\tau$, since $A_j$ is nonnegative and nondecreasing, applying the identity $\bigl[[x]^+-a\bigr]^+=[x-a]^+$ for $a\geq0$ gives
\begin{equation*}
\left[\left[u-\sum_{r=1}^{j-1}\mu_r^{-1}A_r(u)\right]^+-\mu_j^{-1}A_j(s)\right]^+\leq \left[u-\sum_{r=1}^{j}\mu_r^{-1}A_r(u)\right]^+.
\end{equation*}
Taking the supremum over $0\leq u\leq s\leq\tau$, we obtain
\begin{equation*}
R_j^*(\tau) \leq \sup_{0\leq u\leq \tau} \left[ u-\sum_{r=1}^{j}\mu_r^{-1}A_r(u) \right]^+.
\end{equation*}
Conversely, for any $0\leq u\leq\tau$, by the definition of $R_j^*$,
\begin{equation*}
R_j^*(\tau) \geq \left[ R_{j-1}^*(u)-\mu_j^{-1}A_j(u) \right]^+ \geq \left[ \left[ u-\sum_{r=1}^{j-1}\mu_r^{-1}A_r(u) \right]^+ -\mu_j^{-1}A_j(u) \right]^+ = \left[ u-\sum_{r=1}^{j}\mu_r^{-1}A_r(u) \right]^+.
\end{equation*}
Taking the supremum over $u\in[0,\tau]$ gives the reverse inequality. Therefore,
\begin{equation*}
R_j^*(\tau)=\sup_{0\leq u\leq\tau}\left[u-\sum_{r=1}^{j}\mu_r^{-1}A_r(u)\right]^+ =\sup_{0\leq u\leq\tau}\left[-g_j(u)\right]^+ =-\inf_{0\leq u\leq\tau}\left(g_j(u)\wedge0\right).
\end{equation*}
\end{proof}

\section{Proof of Theorem~\ref{thm:fluid-throughput-optimality}}
\label{sec:proof-fluid-throughput-optimality}

\begin{proof}
We first show that the $c\mu$ policy $\pi^*$ completes all requests in finite time. By \eqref{eq:prefix-workload-W}, the total workload is
\begin{equation*}
W_J^{\pi^*}(\xi)=\sum_{j=1}^J\mu_j^{-1}p_{j,L_j}^{\pi^*}(\xi).
\end{equation*}
Since $\mu_j>0$ and $p_{j,L_j}^{\pi^*}(\xi)\geq0$, $W_J^{\pi^*}(\xi)=0$ if and only if $p_{j,L_j}^{\pi^*}(\xi)=0$ for all $j$. For any clock time $\tau\geq H$, all requests have arrived, so $W_J^{\pi^*}((f^*)^{-1}(\tau))=0$ if and only if all requests have been completed by time $\tau$.

Let $\overline\tau:=H+\sum_{j=1}^J\mu_j^{-1}A_j(H)$. We will show that $W_J^{\pi^*}((f^*)^{-1}(\overline\tau))=0$. Recall that, by \eqref{eq:prefix-workload-attainment},
\begin{equation*}
W_j^{\pi^*}(\xi)=\varphi_j(f^*(\xi)),\qquad j=1,\dots,J,\quad \xi\geq0.
\end{equation*}
By this identity, it suffices to prove that $\varphi_J(\overline\tau)=0$. We do so by determining the minimum of $g_J$ over $[0,\overline\tau]$. For $0\leq s\leq H$,
\begin{equation*}
g_J(s)=\sum_{j=1}^J\mu_j^{-1}A_j(s)-s\geq-s\geq-H.
\end{equation*}
For $H\leq s\leq\overline\tau$, the arrival truncation and the definition of $\overline\tau$ give
\begin{equation*}
g_J(s)=\sum_{j=1}^J\mu_j^{-1}A_j(H)-s =\overline\tau-H-s\geq-H,
\end{equation*}
with equality at $s=\overline\tau$. Thus, $\inf_{0\leq s\leq\overline\tau}\left(g_J(s)\wedge0\right)=-H$. By the definition of $\varphi_J$, we obtain
\begin{equation*}
\varphi_J(\overline\tau)=g_J(\overline\tau)-\inf_{0\leq s\leq\overline\tau}\left(g_J(s)\wedge0\right) =(-H)-(-H)=0.
\end{equation*}
Equation~\eqref{eq:prefix-workload-attainment} now gives $W_J^{\pi^*}((f^*)^{-1}(\overline\tau))=0$. Since $\overline\tau\geq H$, all requests are completed by this time, so $\tau_{\mathrm{c}}^{\pi^*}\leq\overline\tau<\infty$.

Next, we prove that if an arbitrary admissible policy $\pi$ completes all requests at $\tau_{\mathrm{c}}^\pi$, then $\pi^*$ also completes all requests by that time. For such a policy $\pi$, we have $W_J^\pi\!\left((f^\pi)^{-1}(\tau_{\mathrm{c}}^\pi)\right)=0$.
Recall that, by \eqref{eq:W-lower-varphi},
\begin{equation*}
W_J^\pi\!\left((f^\pi)^{-1}(\tau)\right)\geq\varphi_J(\tau),\qquad\tau\geq0.
\end{equation*}
Together with \eqref{eq:prefix-workload-attainment}, this gives
\begin{equation*}
0=W_J^\pi\!\left((f^\pi)^{-1}(\tau_{\mathrm{c}}^\pi)\right)\geq\varphi_J(\tau_{\mathrm{c}}^\pi)=W_J^{\pi^*}\!\left((f^*)^{-1}(\tau_{\mathrm{c}}^\pi)\right)\geq0.
\end{equation*}
Since all requests have arrived by $\tau_{\mathrm{c}}^\pi$, zero workload implies that $\pi^*$ has completed all requests by $\tau_{\mathrm{c}}^\pi$. Thus, $\tau_{\mathrm{c}}^{\pi^*}\leq\tau_{\mathrm{c}}^\pi$, proving \eqref{eq:fluid-completion-time-optimality}.

Finally, we prove the throughput claim. Every policy that completes all requests produces the same total number $\sum_{j=1}^Jo_jA_j(H)$ of output tokens. Therefore, when this number is positive, \eqref{eq:fluid-throughput-at-completion} and \eqref{eq:fluid-completion-time-optimality} give
\begin{equation*}
\operatorname{TP}^\pi(\tau_c^\pi)=\frac{\sum_{j=1}^Jo_jA_j(H)}{\tau_c^\pi} \leq\frac{\sum_{j=1}^Jo_jA_j(H)}{\tau_c^{\pi^*}} =\operatorname{TP}^{\pi^*}(\tau_c^{\pi^*}),
\end{equation*}
which proves \eqref{eq:fluid-throughput-optimality}.
\end{proof}

\section{Derivation of the No-Arrival Closed Forms}
\label{sec:no-arrival-closed-forms}

Let $B_j:=p_{j,L_j}(0)$ for $j=1,\dots,J$. Since $\lambda_j(\cdot)\equiv0$, we have $A_j(\tau)=B_j$. Substituting $A_j(s)=B_j$ into the definition of $R_j^*$ in \eqref{eq:target-residual} gives
\begin{equation}
R_j^*(\tau)=\left[\tau-\sum_{\ell=1}^{j}\mu_\ell^{-1}B_\ell\right]^+, \qquad j=0,1,\ldots,J. \label{eq:no-arrival-residual}
\end{equation}
For $\tau\in[\sum_{\ell=1}^{j-1}\mu_\ell^{-1}B_\ell,\sum_{\ell=1}^{j}\mu_\ell^{-1}B_\ell]$, we have $R_\ell^*(\tau)=0$ for $\ell=j,\dots,J$. Substituting \eqref{eq:no-arrival-residual} into the definition of $F$ in \eqref{eq:target-F} yields
\begin{equation*}
F(\tau)=\sum_{\ell=1}^{j-1}\frac{\tilde m_\ell B_\ell}{M}+\frac{\tilde m_j\mu_j}{M}\left(\tau-\sum_{\ell=1}^{j-1}\mu_\ell^{-1}B_\ell\right).
\end{equation*}
On this interval, substituting \eqref{eq:no-arrival-residual} into the definition of $\mathcal X_j^*$ in \eqref{eq:target-admission} gives
\begin{equation*}
\mathcal X_j^*(F(\tau))=\mu_j\left(R_{j-1}^*(\tau)-R_j^*(\tau)\right)=\mu_j\left(\tau-\sum_{\ell=1}^{j-1}\mu_\ell^{-1}B_\ell\right)=\frac{M}{\tilde m_j}\left(F(\tau)-\sum_{\ell=1}^{j-1}\frac{\tilde m_\ell B_\ell}{M}\right),
\end{equation*}
where the last equality uses the expression for $F(\tau)$ above. Since $F$ is continuous and strictly increasing, $\xi=F(\tau)$ ranges over $[\sum_{\ell=1}^{j-1}\frac{\tilde m_\ell B_\ell}{M},\sum_{\ell=1}^{j}\frac{\tilde m_\ell B_\ell}{M}]$ on the above $\tau$ interval. By \eqref{eq:no-arrival-residual}, $\mathcal X_j^*(\xi)=0$ before this step-time interval and $\mathcal X_j^*(\xi)=B_j$ after it. Hence, for all $\xi\geq0$,
\begin{equation*}
\mathcal X_j^*(\xi) = \left[ \frac{M}{\tilde m_j}\left(\xi-\sum_{\ell=1}^{j-1}\frac{\tilde m_\ell B_\ell}{M}\right) \right]_0^{B_j}.
\end{equation*}
Using $B_j=p_{j,L_j}(0)$ gives \eqref{eq:no-arrival-explicit-policy}.

\section{Proof of Proposition~\ref{prop:pathwise-memory-safety}}
\label{sec:proof-pathwise-memory-safety}

\begin{proof}
According to Definition~\ref{def:moving-average-policy}, $\sum_{u<t}\eta_j^r(u)$ and $\sum_{u<t}y_j^r(u)$ are determined by the history available at step $t$. The term $\sum_{\ell<j}\tilde m_\ell y_\ell^r(t)$ involves only higher-priority admission plans, which have already been computed using the same history. Thus, by recursion over the priority classes, $y_j^r(t)$ is $\cG_t^r$-measurable and hence non-anticipative. Since the moving average and time shifts in~\eqref{eq:smoothed-admission}--\eqref{eq:prelimit-V-shift} use only current and past values, all control increments $\Delta\mathcal U_{j,k}^r(t)$ and $\Delta\mathcal V_{j,k}^r(t)$ are also $\cG_t^r$-measurable, so the policy is non-anticipative.

We next prove by induction on $t$ that the cumulative planned admissions of class $j$ through step $t$ do not exceed its initial queue length plus cumulative arrivals before step $t$:
\begin{equation}
\sum_{u=0}^{t}y_j^r(u)\leq p_{j,L_j}^{r}(0)+\sum_{u=0}^{t-1}\eta_j^r(u), \qquad t\in\mathbb Z_+. \label{eq:prelimit-postplan-demand}
\end{equation}
For $t=0$, the planning rule \eqref{eq:planned-admission-policy} gives $0\leq y_j^r(0)\leq p_{j,L_j}^{r}(0)$. Suppose that \eqref{eq:prelimit-postplan-demand} holds at $t-1$. The induction hypothesis and the nonnegativity of arrivals give
\begin{equation*}
p_{j,L_j}^{r}(0) +\sum_{u=0}^{t-1}\eta_j^r(u) -\sum_{u=0}^{t-1}y_j^r(u)= \left[p_{j,L_j}^{r}(0)+\sum_{u=0}^{t-2}\eta_j^r(u)-\sum_{u=0}^{t-1}y_j^r(u)\right]+\eta_j^r(t-1)\geq \eta_j^r(t-1) \geq0.
\end{equation*}
Since this quantity is nonnegative, the first argument of the minimum in the planning rule \eqref{eq:planned-admission-policy} gives
\begin{equation*}
y_j^r(t)\leq\left[p_{j,L_j}^{r}(0)+\sum_{u=0}^{t-1}\eta_j^r(u)-\sum_{u=0}^{t-1}y_j^r(u)\right]^+=p_{j,L_j}^{r}(0)+\sum_{u=0}^{t-1}\eta_j^r(u)-\sum_{u=0}^{t-1}y_j^r(u).
\end{equation*}
Consequently,
\begin{equation*}
\begin{aligned}
\sum_{u=0}^{t}y_j^r(u)&= \sum_{u=0}^{t-1}y_j^r(u)+y_j^r(t)
\\&\leq \sum_{u=0}^{t-1}y_j^r(u)+\left[p_{j,L_j}^{r}(0)+\sum_{u=0}^{t-1}\eta_j^r(u)-\sum_{u=0}^{t-1}y_j^r(u)\right]= p_{j,L_j}^{r}(0)+\sum_{u=0}^{t-1}\eta_j^r(u).
\end{aligned}
\end{equation*}
This completes the induction.

Summing the admission increments in \eqref{eq:smoothed-admission} and setting $u=s-v$ gives
\begin{equation*}
\begin{aligned}
\mathcal U_{j,L_j}^{r}(t+1)&=\frac1{w_r}\sum_{v=0}^{w_r-1}\sum_{s=0}^{t}y_j^r(s-v)
=\frac1{w_r}\sum_{v=0}^{w_r-1}\sum_{u=0}^{t-v}y_j^r(u)\\
&\leq\left(\frac1{w_r}\sum_{v=0}^{w_r-1}1\right)\sum_{u=0}^{t}y_j^r(u)=\sum_{u=0}^{t}y_j^r(u)
\leq p_{j,L_j}^{r}(0)+\sum_{u=0}^{t-1}\eta_j^r(u),
\end{aligned}
\end{equation*}
where $y_j^r(u)=0$ for $u<0$ and sums with an upper limit below zero are empty. Subtracting $\mathcal U_{j,L_j}^{r}(t)$ from both sides of the preceding inequality and using the waiting-queue balance yields
\begin{equation*}
0\leq\Delta\mathcal U_{j,L_j}^{r}(t) \leq p_{j,L_j}^{r}(0)+\sum_{u=0}^{t-1}\eta_j^r(u)-\mathcal U_{j,L_j}^{r}(t) =p_{j,L_j}^{r}(t).
\end{equation*}

Next, we verify that the service increment at each intermediate stage does not exceed its queue length. For a prefill stage $k<L_j$, the queueing dynamics give
\begin{equation*}
p_{j,k}^r(u+1)-p_{j,k}^r(u)=\Delta\mathcal U_{j,k+1}^r(u)-\Delta\mathcal U_{j,k}^r(u).
\end{equation*}
Summing over $u=0,\dots,t-1$ and telescoping the left-hand side, we obtain
\begin{equation*}
p_{j,k}^r(t)-p_{j,k}^r(0)= \sum_{u=0}^{t-1}\left[p_{j,k}^r(u+1)-p_{j,k}^r(u)\right]= \sum_{u=0}^{t-1}\left[\Delta\mathcal U_{j,k+1}^r(u)-\Delta\mathcal U_{j,k}^r(u)\right].
\end{equation*}
Since the intermediate stages are initially empty, $p_{j,k}^r(0)=0$. Substituting the stage shifts in \eqref{eq:prelimit-U-shift} gives
\begin{equation*}
\begin{aligned}
p_{j,k}^r(t)&= \sum_{u=0}^{t-1}\Delta\mathcal U_{j,L_j}^r(u-(L_j-k-1)) -\sum_{u=0}^{t-1}\Delta\mathcal U_{j,L_j}^r(u-(L_j-k))\\
&= \Delta\mathcal U_{j,L_j}^r(t-(L_j-k))= \Delta\mathcal U_{j,k}^r(t).
\end{aligned}
\end{equation*}
The same calculation for the first decode stage and the remaining decode stages gives
\begin{equation*}{\medmuskip=1mu\thickmuskip=2mu\thinmuskip=2mu\scriptspace=0.25pt\begin{aligned}
&d_{j,o_j-1}^{r}(t)=\sum_{u=0}^{t-1}\left[\Delta\mathcal U_{j,1}^{r}(u)-\Delta\mathcal V_{j,o_j-1}^{r}(u)\right] =\Delta\mathcal U_{j,L_j}^{r}(t-L_j) =\Delta\mathcal V_{j,o_j-1}^{r}(t),\quad o_j\geq2,\\
&d_{j,k}^{r}(t)=\sum_{u=0}^{t-1}\left[\Delta\mathcal V_{j,k+1}^{r}(u)-\Delta\mathcal V_{j,k}^{r}(u)\right] =\Delta\mathcal U_{j,L_j}^{r}\bigl(t-(L_j+o_j-1-k)\bigr) =\Delta\mathcal V_{j,k}^{r}(t),\quad k=1,\dots,o_j-2.
\end{aligned}}\end{equation*}
Thus each intermediate stage serves exactly the requests available there. Together with the admission bound proved above, this verifies the queue-length bounds on all stage-completion increments and the nonnegativity of all queueing states.

We now prove that the memory constraint holds. Induction over the priority classes in the planning rule \eqref{eq:planned-admission-policy} gives, pathwise,
\begin{equation}
\sum_{j=1}^{J}\tilde m_jy_j^{r}(t) \leq\frac{w_r}{w_r+H_{\mathrm{kv}}}M. \label{eq:planning-budget}
\end{equation}
Substituting the moving-average admission formula in \eqref{eq:smoothed-admission} into the expression for $m^r(t)$ in \eqref{eq:prelimit-memory-convolution} gives
\begin{equation*}
m^{r}(t) =\sum_{j=1}^{J}\sum_{h=0}^{H_{\mathrm{kv}}} k_j^{\mathrm{kv}}(h)\Delta\mathcal U_{j,L_j}^{r}(t-h)=\frac{1}{w_r} \sum_{j=1}^{J}\sum_{h=0}^{H_{\mathrm{kv}}}\sum_{v=0}^{w_r-1} k_j^{\mathrm{kv}}(h)y_j^{r}(t-h-v).
\end{equation*}
Making the change of variable $u=h+v$, we obtain
\begin{equation*}
m^r(t)=\frac1{w_r}\sum_{j=1}^{J}\sum_{h=0}^{H_{\mathrm{kv}}}\sum_{u=h}^{h+w_r-1}k_j^{\mathrm{kv}}(h)y_j^r(t-u).
\end{equation*}
To exchange the sums over $h$ and $u$, note that the full range of $u$ is $0\leq u\leq w_r+H_{\mathrm{kv}}-1$. For a fixed $u$, the admissible values of $h$ are exactly those satisfying
$0\leq h\leq H_{\mathrm{kv}}$ and $0\leq u-h\leq w_r-1$.
After exchanging the order of summation, we factor $y_j^r(t-u)$ out of the sum over $h$, since it does not depend on $h$. This gives
\begin{equation*}
{\medmuskip=2mu\thickmuskip=2.5mu\thinmuskip=1.5mu \begin{aligned}
m^{r}(t)&=\frac{1}{w_r} \sum_{u=0}^{w_r+H_{\mathrm{kv}}-1}\sum_{j=1}^{J} y_j^{r}(t-u) \sum_{\substack{0\leq h\leq H_{\mathrm{kv}}\\
                  0\leq u-h\leq w_r-1}} k_j^{\mathrm{kv}}(h)\\
&\leq\frac{1}{w_r} \sum_{u=0}^{w_r+H_{\mathrm{kv}}-1}\sum_{j=1}^{J} \tilde m_jy_j^{r}(t-u)\\
&\leq \frac{w_r+H_{\mathrm{kv}}}{w_r} \frac{w_r}{w_r+H_{\mathrm{kv}}}M=M.
\end{aligned}}
\end{equation*}
This establishes the admissibility of the policy.

\end{proof}

\section{Proof of Theorem~\ref{thm:moving-average-fluid-limit}}
\label{sec:proof-moving-average-fluid-limit}

We first show that the moving average and the fixed stage delays produce vanishing errors under fluid scaling, so it suffices to analyze the cumulative planned-admission processes and the clock approximation constructed from them (Lemma~\ref{lem:prelimit-uoc-estimates}). We then establish joint $C$-tightness and show that every subsequential limit satisfies the fluid model's queue nonnegativity constraints, memory constraint, and clock equation (Lemma~\ref{lem:ctight-subsequential-limits}). The key step is to recover the priority rule in the limit: on any interval where a class's waiting-queue length is positive, its priority prefix uses the full capacity and the corresponding residual service budget remains constant. This yields the reflection representation that identifies every subsequential limit with the fluid $c\mu$ policy (Lemma~\ref{lem:residual-representation-uniqueness}).

We begin by defining the auxiliary processes used throughout the proof.

\begin{definition}[Auxiliary processes]
For $N\in\mathbb{Z}_+$, define the cumulative planned-admission process by
\begin{equation}
P_j^r(N):=\sum_{t=0}^{N-1}y_j^r(t), \qquad P_j^r(0)=0. \label{eq:P-def}
\end{equation}
For $\xi\geq0$, define its fluid-scaled version by
\begin{equation*}
\bar P_j^r(\xi):=\frac1rP_j^r(\lfloor r\xi\rfloor).
\end{equation*}
Let $n=\lfloor r\xi\rfloor$. Define the linear interpolation $\widehat P_j^r$ of $\bar P_j^r$ and use it to construct a clock approximation $\widehat f^{\,r}(\xi)$ of $\bar f^r(\xi)$ as follows:
\begingroup
\setlength{\abovedisplayskip}{3pt}
\setlength{\belowdisplayskip}{3pt}
\setlength{\abovedisplayshortskip}{0pt}
\setlength{\belowdisplayshortskip}{3pt}
{\setlength{\jot}{0pt}
\begin{align}
\widehat P_j^r(\xi) &:=\bar P_j^r(\xi)+\frac{r\xi-n}{r}y_j^r(n), \label{eq:P-linear}\\
\widehat f^{\,r}(\xi) &:=\sum_{j=1}^J(a_{\mathrm{pf}}s_j+a_{\mathrm{kv}}\tilde m_j)\widehat P_j^r(\xi)+b\xi. \label{eq:clock-proxy}
\end{align}
}
\endgroup

The continuous approximation of the fluid-scaled cumulative arrival process $\bar A_j^r$ is defined by
\begin{equation}
\widehat A_j^r(\xi):=A_j(\widehat f^{\,r}(\xi))=p_{j,L_j}(0)+\int_0^{\widehat f^{\,r}(\xi)}\lambda_j(s)\,ds.
\label{eq:A-proxy}
\end{equation}
\end{definition}

\begin{lemma}
\label{lem:prelimit-uoc-estimates}
Under the initial condition of Theorem~\ref{thm:moving-average-fluid-limit}, consider the smoothed $c\mu$-inspired policy in Definition~\ref{def:moving-average-policy}. Then the following statements hold.
\begin{enumerate}
  \item[\textup{(i)}] For every scale $r$, class $j$, and integer time $t$,
\begingroup
\setlength{\abovedisplayskip}{3pt}
\setlength{\belowdisplayskip}{3pt}
\setlength{\abovedisplayshortskip}{0pt}
\setlength{\belowdisplayshortskip}{3pt}
{\setlength{\jot}{0pt}
\begin{align}
0\leq y_j^r(t) &\leq \frac{M}{\tilde m_j}\frac{w_r}{w_r+H_{\mathrm{kv}}} \leq\frac{M}{\tilde m_j}, \label{eq:prop7-y-bound}\\
0\leq\Delta\mathcal U_{j,L_j}^{r}(t) &\leq \frac{M}{\tilde m_j}\frac{w_r}{w_r+H_{\mathrm{kv}}} \leq\frac{M}{\tilde m_j}. \label{eq:prop7-actual-increment-bound}
\end{align}
}
\endgroup
Moreover,
\begin{equation}
\sum_{j=1}^J\tilde m_jy_j^r(t) \leq\frac{w_rM}{w_r+H_{\mathrm{kv}}}. \label{eq:prop7-planning-budget}
\end{equation}
  \item[\textup{(ii)}] For every class $j$ and every $T<\infty$,
\begin{equation}
\sup_{0\leq\xi\leq T}\left| \bar P_j^r(\xi) -\frac{1}{r}\mathcal U_{j,L_j}^{r}(\lfloor r\xi\rfloor) \right|\!\leq\!\frac{M}{2\tilde m_j}\frac{w_r}{r} \to0.
\label{eq:prop7-planned-actual}
\end{equation}
  \item[\textup{(iii)}] For every class $j$ and every $T<\infty$, the following bounds hold for $k=1,\dots,L_j$ and $k=1,\dots,o_j-1$, respectively:
\begingroup
\medmuskip=2.5mu\thickmuskip=3mu
\begin{align}
\sup_{0\leq\xi\leq T}\left|\frac{\mathcal U_{j,k}^r(\lfloor r\xi\rfloor)}r-\frac{\mathcal U_{j,L_j}^r(\lfloor r\xi\rfloor)}r\right|&\leq\frac{H_{\mathrm{kv}}M}{r\tilde m_j}\longrightarrow0, \label{eq:prop7-fixed-shift-U}\\
\sup_{0\leq\xi\leq T}\left|\frac{\mathcal V_{j,k}^r(\lfloor r\xi\rfloor)}r-\frac{\mathcal U_{j,L_j}^r(\lfloor r\xi\rfloor)}r\right|&\leq\frac{H_{\mathrm{kv}}M}{r\tilde m_j}\longrightarrow0. \label{eq:prop7-fixed-shift-V}
\end{align}
\endgroup
  \item[\textup{(iv)}] For every $T<\infty$,
\begin{equation}
\sup_{0\leq\xi\leq T}\Biggl|\bar f^r(\xi) -\sum_{j=1}^J(a_{\mathrm{pf}}s_j+a_{\mathrm{kv}}\tilde m_j)\bar P_j^r(\xi)-b\xi\Biggr|\longrightarrow0.
\label{eq:prop7-clock-planned}
\end{equation}
  \item[\textup{(v)}] For every class $j$, every fixed $T<\infty$, and all $\xi_1,\xi_2\in[0,T]$,
\begin{equation}
\sup_{0\leq\xi\leq T}\left|\widehat P_j^r(\xi)-\bar P_j^r(\xi)\right|\leq\frac{M}{r\tilde m_j}\longrightarrow0, \qquad\text{and}\qquad \left|\widehat P_j^r(\xi_2)-\widehat P_j^r(\xi_1)\right|\leq\frac{M}{\tilde m_j}|\xi_2-\xi_1|.
\label{eq:prop7-P-proxy-bounds}
\end{equation}
Moreover, define
\begin{equation}
L_f:=b+M\sum_{j=1}^J \left(\frac{a_{\mathrm{pf}}s_j}{\tilde m_j}+a_{\mathrm{kv}}\right). \label{eq:prop7-Lf}
\end{equation}
Then $\widehat f^{\,r}$ is $L_f$-Lipschitz and
\begin{equation}
\sup_{0\leq\xi\leq T}\left| \bar f^r(\xi)-\widehat f^{\,r}(\xi) \right| \longrightarrow0. \label{eq:prop7-clock-proxy-close}
\end{equation}
  \item[\textup{(vi)}] For every class $j$ and every fixed $T<\infty$,
\begin{equation}
\sup_{0\leq\xi\leq T}\left| \bar A_j^r(\xi)-\widehat A_j^r(\xi) \right| \xrightarrow{\mathbb P}0. \label{eq:prop7-arrival-composition-close-vector}
\end{equation}
\end{enumerate}
\end{lemma}

\begin{remark}[Interpretation of the prelimit estimates]\leavevmode
\begin{enumerate}
\item[\textup{(i)}] Equations~\eqref{eq:prop7-y-bound} and \eqref{eq:prop7-actual-increment-bound} give memory-based upper bounds on both the class-$j$ planned-admission amount $y_j^r(t)$ and the admission increment $\Delta\mathcal U_{j,L_j}^r(t)$ at each step $t$. Equation~\eqref{eq:prop7-planning-budget} states that the plans for all classes share one memory-based planning budget.

\item[\textup{(ii)}] Equation~\eqref{eq:prop7-planned-actual} shows that the difference between the fluid-scaled cumulative planned admissions $\bar P_j^r(\xi)$ and the fluid-scaled cumulative admissions $\mathcal U_{j,L_j}^r(\lfloor r\xi\rfloor)/r$ vanishes uniformly on compact step-time intervals. Thus, spreading each plan over $w_r$ steps has a negligible effect on cumulative admissions at the fluid scale.

\item[\textup{(iii)}] The prefill-stage delay bound \eqref{eq:prop7-fixed-shift-U} and the decode-stage delay bound \eqref{eq:prop7-fixed-shift-V} show that the fixed delays between admission and later stages disappear under fluid scaling. The fluid-scaled cumulative stage-completion processes $\mathcal U_{j,k}^r(\lfloor r\xi\rfloor)/r$ and $\mathcal V_{j,k}^r(\lfloor r\xi\rfloor)/r$ are therefore asymptotically equivalent to the fluid-scaled cumulative admission process $\mathcal U_{j,L_j}^r(\lfloor r\xi\rfloor)/r$.

\item[\textup{(iv)}] The clock approximation in \eqref{eq:prop7-clock-planned} uses the fluid-scaled cumulative planned admissions $\bar P_j^r(\xi)$ to approximate the real-time clock $\bar f^r(\xi)$.

\item[\textup{(v)}] Equations~\eqref{eq:prop7-P-proxy-bounds} and \eqref{eq:prop7-clock-proxy-close} show that the linear interpolation $\widehat P_j^r(\xi)$ and the fluid-scaled cumulative planned-admission process $\bar P_j^r(\xi)$, as well as $\widehat f^{\,r}(\xi)$ and $\bar f^r(\xi)$, are asymptotically equivalent uniformly on compact intervals. Moreover, $\widehat P_j^r$ and $\widehat f^{\,r}$ are Lipschitz with constants $M/\tilde m_j$ and $L_f$, respectively, independent of $r$.

\item[\textup{(vi)}] Equation~\eqref{eq:prop7-arrival-composition-close-vector} shows that the difference between the fluid-scaled cumulative arrival process $\bar A_j^r$ and its continuous approximation $\widehat A_j^r$ vanishes uniformly on compact intervals in probability.
\end{enumerate}
\end{remark}

\begin{proof}

\noindent\textbf{(i)} By the policy definition,
\begin{equation*}
0\leq\tilde m_jy_j^r(t) \leq\left[\frac{w_rM}{w_r+H_{\mathrm{kv}}}-\sum_{\ell=1}^{j-1}\tilde m_\ell y_\ell^r(t)\right]^+ \leq\frac{w_rM}{w_r+H_{\mathrm{kv}}}.
\end{equation*}
Equation~\eqref{eq:smoothed-admission} expresses the admission increments as the moving average of the planned-admission amounts, which are zero at negative times. Hence,
\begin{equation*}
0\leq\Delta\mathcal U_{j,L_j}^{r}(t) =\frac1{w_r}\sum_{v=0}^{w_r-1}y_j^r(t-v) \leq\frac{M}{\tilde m_j}\frac{w_r}{w_r+H_{\mathrm{kv}}} \leq\frac{M}{\tilde m_j}.
\end{equation*}

The planning-budget bound follows by induction over the priority classes: starting from zero, each class uses at most the budget remaining after the preceding classes. Specifically, if the budget holds through class $j-1$, its remaining budget is nonnegative, and
\begin{equation*}
\begin{aligned}
\sum_{\ell=1}^{j}\tilde m_\ell y_\ell^r(t)&= \sum_{\ell=1}^{j-1}\tilde m_\ell y_\ell^r(t) +\tilde m_jy_j^r(t)\\
&\leq \sum_{\ell=1}^{j-1}\tilde m_\ell y_\ell^r(t) +\left[ \frac{w_rM}{w_r+H_{\mathrm{kv}}} -\sum_{\ell=1}^{j-1}\tilde m_\ell y_\ell^r(t) \right]^+\\
&= \frac{w_rM}{w_r+H_{\mathrm{kv}}}.
\end{aligned}
\end{equation*}
Taking $j=J$ proves \eqref{eq:prop7-planning-budget}.
\par\medskip
\noindent\textbf{(ii)} Fix $N\in\mathbb{Z}_+$. Changing variables $t=s-v$ in the inner sum, and using $y_j^r(t)=0$ for $t<0$, gives
\begin{equation*}
\mathcal U_{j,L_j}^{r}(N)= \sum_{s=0}^{N-1}\frac{1}{w_r} \sum_{v=0}^{w_r-1}y_j^r(s-v)= \frac{1}{w_r}\sum_{s=0}^{N-1}\sum_{t=s-w_r+1}^{s}y_j^r(t)= \frac{1}{w_r}\sum_{s=0}^{N-1}\sum_{t=\max\{0,s-w_r+1\}}^{s}y_j^r(t).
\end{equation*}
For each fixed $0\leq t\leq N-1$, the admissible values of $s$ satisfy $t\leq s\leq\min\{N-1,t+w_r-1\}$. Exchanging the order of summation therefore yields
\begin{equation*}
\mathcal U_{j,L_j}^{r}(N)= \frac{1}{w_r}\sum_{t=0}^{N-1}\sum_{s=t}^{\min\{N-1,t+w_r-1\}}y_j^r(t)= \frac{1}{w_r}\sum_{t=0}^{N-1} \min\{w_r,N-t\}\,y_j^r(t)= \sum_{t=0}^{N-1} \min\left\{1,\frac{N-t}{w_r}\right\}y_j^r(t).
\end{equation*}
Subtracting this expression from $P_j^r(N)$, and setting $v=N-t$ in the last line, gives
\begin{equation*}
\begin{aligned}
P_j^r(N) -\mathcal U_{j,L_j}^{r}(N)&= \sum_{t=0}^{N-1} \left[1-\min\left\{1,\frac{N-t}{w_r}\right\}\right]y_j^r(t)\\
&= \sum_{t=\max\{0,N-w_r+1\}}^{N-1} \left(1-\frac{N-t}{w_r}\right)y_j^r(t)\\
&= \frac{1}{w_r} \sum_{v=1}^{\min\{w_r-1,N\}}(w_r-v)y_j^r(N-v).
\end{aligned}
\end{equation*}
By \textup{(i)},
\begin{equation*}
0\leq P_j^r(N) -\mathcal U_{j,L_j}^{r}(N)\leq\frac{1}{w_r}\sum_{v=1}^{w_r-1}(w_r-v) \frac{M}{\tilde m_j}\frac{w_r}{w_r+H_{\mathrm{kv}}}=\frac{M}{\tilde m_j}\frac{w_r}{w_r+H_{\mathrm{kv}}} \frac{w_r-1}{2}.
\end{equation*}
Set $N=\lfloor r\xi\rfloor$, divide by $r$, and take the supremum over $0\leq\xi\leq T$.  Since $w_r/r\to0$, this gives \eqref{eq:prop7-planned-actual}.
\par\medskip
\noindent\textbf{(iii)} For any integers $N\geq0$ and $0\leq h\leq H_{\mathrm{kv}}$, \textup{(i)} gives
\begin{equation*}
0\leq\mathcal U_{j,L_j}^{r}(N)-\mathcal U_{j,L_j}^{r}([N-h]^+)
=\sum_{t=[N-h]^+}^{N-1}\Delta\mathcal U_{j,L_j}^{r}(t)
\leq\left[N-[N-h]^+\right]\frac{M}{\tilde m_j}
\leq h\frac{M}{\tilde m_j}\leq H_{\mathrm{kv}}\frac{M}{\tilde m_j}.
\end{equation*}
The stage-shift relations \eqref{eq:prelimit-U-shift} and \eqref{eq:prelimit-V-shift}, together with zero initial cumulative controls and the convention of zero increments at negative times, give
\begin{equation*}
\begin{aligned}
\mathcal U_{j,k}^r(N)&=\mathcal U_{j,L_j}^r\bigl([N-(L_j-k)]^+\bigr),\qquad 1\leq k\leq L_j,\\
\mathcal V_{j,k}^r(N)&=\mathcal U_{j,L_j}^r\bigl([N-(L_j+o_j-1-k)]^+\bigr),\qquad 1\leq k\leq o_j-1.
\end{aligned}
\end{equation*}
Apply the preceding bound with $h=L_j-k$ and $h=L_j+o_j-1-k$, respectively; both delays lie between zero and $H_{\mathrm{kv}}$. Setting $N=\lfloor r\xi\rfloor$, dividing by $r$, and taking the supremum over $0\leq\xi\leq T$ proves the two stated bounds.
\par\medskip
\noindent\textbf{(iv)} By the definitions of $\Delta^r$ and $f^r$, for $N=\lfloor r\xi\rfloor$,
\begin{equation*}
f^r(N)= \sum_{t=0}^{N-1}\Delta^r(t) = a_{\mathrm{pf}}^r\sum_{t=0}^{N-1}n_{\mathrm{pf}}^r(t) +a_{\mathrm{kv}}^r\sum_{t=0}^{N-1}n_{\mathrm{kv}}^r(t)+b^rN.
\end{equation*}
Substituting the expressions for $n_{\mathrm{pf}}^r(t)$ and $n_{\mathrm{kv}}^r(t)$ yields
\begin{equation*}
\begin{aligned}
f^r(N) &=a_{\mathrm{pf}}^r\sum_{j=1}^J\sum_{k=1}^{L_j} \ell_{j,k}\mathcal U_{j,k}^{r}(N) +a_{\mathrm{kv}}^r\sum_{j=1}^J\sum_{k=1}^{L_j} \left[c_{\mathrm{chunk}}(L_j-k)+\ell_{j,k}\right]\mathcal U_{j,k}^{r}(N)\\
&\quad+a_{\mathrm{kv}}^r\sum_{j=1}^J\sum_{k=1}^{o_j-1} (s_j+o_j-k)\mathcal V_{j,k}^{r}(N) +b^rN.
\end{aligned}
\end{equation*}
Combining \textup{(ii)} and \textup{(iii)} by the triangle inequality gives
\begin{equation*}
{\medmuskip=1mu\thickmuskip=2mu\thinmuskip=2mu\scriptspace=0.25pt\begin{aligned}
&\left|\mathcal U_{j,k}^{r}(N)-P_j^r(N)\right|\leq\left|\mathcal U_{j,k}^{r}(N)-\mathcal U_{j,L_j}^{r}(N)\right|+\left|\mathcal U_{j,L_j}^{r}(N)-P_j^r(N)\right| \leq\frac{M}{\tilde m_j}\left(H_{\mathrm{kv}}+\frac{w_r}{2}\right), \quad1 \leq k \leq L_j,\\
&\left|\mathcal V_{j,k}^{r}(N)-P_j^r(N)\right|\leq\left|\mathcal V_{j,k}^{r}(N)-\mathcal U_{j,L_j}^{r}(N)\right|+\left|\mathcal U_{j,L_j}^{r}(N)-P_j^r(N)\right| \leq\frac{M}{\tilde m_j}\left(H_{\mathrm{kv}}+\frac{w_r}{2}\right), \quad1 \leq k \leq o_j -1.
\end{aligned}}
\end{equation*}
Using $\sum_{k=1}^{L_j}\ell_{j,k}=s_j$ and the identity for $\tilde m_j$ in \eqref{eq:prelimit-mtilde-stage-identity}, we can write the clock difference as
\begin{equation*}
{\medmuskip=1mu\thickmuskip=2mu \begin{aligned}
&f^r(N)-\sum_{j=1}^J(a_{\mathrm{pf}}^rs_j+a_{\mathrm{kv}}^r\tilde m_j)P_j^r(N)-b^rN\\
&=\sum_{j=1}^J\Biggl(a_{\mathrm{pf}}^r\sum_{k=1}^{L_j}\ell_{j,k}\bigl[\mathcal U_{j,k}^r(N)-P_j^r(N)\bigr] +a_{\mathrm{kv}}^r\sum_{k=1}^{L_j}[c_{\mathrm{chunk}}(L_j-k)+\ell_{j,k}]\bigl[\mathcal U_{j,k}^r(N)-P_j^r(N)\bigr]\\
&\qquad{}+a_{\mathrm{kv}}^r\sum_{k=1}^{o_j-1}(s_j+o_j-k)\bigl[\mathcal V_{j,k}^r(N)-P_j^r(N)\bigr]\Biggr).
\end{aligned}}
\end{equation*}
Applying the two bounds above to these differences yields
\begin{equation*}
{\medmuskip=1mu\thickmuskip=2mu\thinmuskip=2mu\scriptspace=0pt \begin{aligned}
&\sup_{0\leq\xi\leq T}\Biggl|f^r(N)-\sum_{j=1}^J(a_{\mathrm{pf}}^rs_j+a_{\mathrm{kv}}^r\tilde m_j)P_j^r(N)-b^rN\Biggr|\\
&\leq M\Bigl(H_{\mathrm{kv}}+\frac{w_r}{2}\Bigr)\sum_{j=1}^J\frac1{\tilde m_j}\Biggl[a_{\mathrm{pf}}^r\sum_{k=1}^{L_j}\ell_{j,k} +a_{\mathrm{kv}}^r\Biggl(\sum_{k=1}^{L_j}[c_{\mathrm{chunk}}(L_j-k)+\ell_{j,k}] +\sum_{k=1}^{o_j-1}(s_j+o_j-k)\Biggr)\Biggr]\\
&=M\Bigl(H_{\mathrm{kv}}+\frac{w_r}{2}\Bigr)\sum_{j=1}^J\Bigl(\frac{a_{\mathrm{pf}}^rs_j}{\tilde m_j}+a_{\mathrm{kv}}^r\Bigr)\longrightarrow0.
\end{aligned}}
\end{equation*}
The convergence follows because $r a_{\mathrm{pf}}^r$ and $r a_{\mathrm{kv}}^r$ are bounded by the scaling assumptions and $(H_{\mathrm{kv}}+w_r)/r\to0$. Moreover, \textup{(i)} implies $0\leq \bar P_j^r(\xi)\leq MT/\tilde m_j$ for $0\leq\xi\leq T$. Since $r a_{\mathrm{pf}}^r\to a_{\mathrm{pf}}$, we have
\begin{equation*}
\sup_{0\leq\xi\leq T}\left| (r a_{\mathrm{pf}}^r-a_{\mathrm{pf}})\bar P_j^r(\xi) \right| \leq |r a_{\mathrm{pf}}^r-a_{\mathrm{pf}}|\frac{MT}{\tilde m_j} \longrightarrow0.
\end{equation*}
Similarly, $r a_{\mathrm{kv}}^r\to a_{\mathrm{kv}}$ gives
\begin{equation*}
\sup_{0\leq\xi\leq T}\left| (r a_{\mathrm{kv}}^r-a_{\mathrm{kv}})\bar P_j^r(\xi) \right| \leq |r a_{\mathrm{kv}}^r-a_{\mathrm{kv}}|\frac{MT}{\tilde m_j} \longrightarrow0.
\end{equation*}
For the overhead term, $N/r\leq T$ and $|N/r-\xi|\leq1/r$ imply
\begin{equation*}
\sup_{0\leq\xi\leq T}|b^rN-b\xi|=\sup_{0\leq\xi\leq T}\left| (r b^r-b)\frac{N}{r}+b\left(\frac{N}{r}-\xi\right) \right|\leq T|r b^r-b| +\frac{b}{r} \longrightarrow0,
\end{equation*}
where the convergence follows from $r b^r\to b$. Finally, combining the preceding estimates by the triangle inequality gives
\begin{equation*}
{\medmuskip=2.5mu\thickmuskip=3mu \begin{aligned}
&\sup_{0\leq\xi\leq T}\Biggl|\bar f^r(\xi)-\sum_{j=1}^J(a_{\mathrm{pf}}s_j+a_{\mathrm{kv}}\tilde m_j)\bar P_j^r(\xi)-b\xi\Biggr|\\
&\leq\sup_{0\leq\xi\leq T}\Biggl|\bar f^r(\xi)-\sum_{j=1}^J(a_{\mathrm{pf}}^rs_j+a_{\mathrm{kv}}^r\tilde m_j)P_j^r(N)-b^rN\Biggr|+\sum_{j=1}^Js_j\sup_{0\leq\xi\leq T}\left|(ra_{\mathrm{pf}}^r-a_{\mathrm{pf}})\bar P_j^r(\xi)\right|\\
&\quad+\sum_{j=1}^J\tilde m_j\sup_{0\leq\xi\leq T}\left|(ra_{\mathrm{kv}}^r-a_{\mathrm{kv}})\bar P_j^r(\xi)\right| +\sup_{0\leq\xi\leq T}|b^rN-b\xi|\longrightarrow0.
\end{aligned}}
\end{equation*}
This proves \eqref{eq:prop7-clock-planned}.
\par\medskip
\noindent\textbf{(v)} For $0\leq\xi\leq T$, the definition of $\widehat P_j^r$ in \eqref{eq:P-linear} and the one-step bound on $y_j^r$ in \textup{(i)} give
\begin{equation*}
0\leq\widehat P_j^r(\xi)-\bar P_j^r(\xi) =\frac{r\xi-\lfloor r\xi\rfloor}{r}y_j^r(\lfloor r\xi\rfloor) \leq\frac{M}{r\tilde m_j}.
\end{equation*}
This proves the first estimate in \eqref{eq:prop7-P-proxy-bounds}.  On each grid interval $[n/r,(n+1)/r]$, the slope of $\widehat P_j^r$ is $y_j^r(n)$, and
$0\leq y_j^r(n)\leq\frac{M}{\tilde m_j}$.
Consequently, for any $0\leq\xi_1\leq\xi_2\leq T$,
\begin{equation*}
0\leq\widehat P_j^r(\xi_2)-\widehat P_j^r(\xi_1) \leq\frac{M}{\tilde m_j}(\xi_2-\xi_1),
\end{equation*}
which proves the second estimate.

The definition of $\widehat f^{\,r}$ in \eqref{eq:clock-proxy} now gives, for $0\leq\xi_1\leq\xi_2\leq T$,
\begin{equation*}
{\medmuskip=1.5mu\thickmuskip=2mu
\begin{aligned}
\left|\widehat f^{\,r}(\xi_2)-\widehat f^{\,r}(\xi_1)\right|&\leq\sum_{j=1}^J(a_{\mathrm{pf}}s_j+a_{\mathrm{kv}}\tilde m_j) \left|\widehat P_j^r(\xi_2)-\widehat P_j^r(\xi_1)\right|+b|\xi_2-\xi_1|\\
&\leq\left[ b+M\sum_{j=1}^J \left(\frac{a_{\mathrm{pf}}s_j}{\tilde m_j}+a_{\mathrm{kv}}\right) \right]|\xi_2-\xi_1|=L_f|\xi_2-\xi_1|.
\end{aligned}}
\end{equation*}

Finally, the definition of $\widehat f^{\,r}$ in \eqref{eq:clock-proxy}, the clock approximation in \eqref{eq:prop7-clock-planned}, and the interpolation error bound \eqref{eq:prop7-P-proxy-bounds} give
\begin{equation*}
\begin{aligned}
\sup_{0\leq\xi\leq T}\left| \bar f^r(\xi)-\widehat f^{\,r}(\xi) \right|&\leq\sup_{0\leq\xi\leq T}\left| \bar f^r(\xi) -\sum_{j=1}^J(a_{\mathrm{pf}}s_j+a_{\mathrm{kv}}\tilde m_j) \bar P_j^r(\xi) -b\xi \right|\\
&\quad+\sum_{j=1}^J(a_{\mathrm{pf}}s_j+a_{\mathrm{kv}}\tilde m_j) \sup_{0\leq\xi\leq T}\left| \bar P_j^r(\xi)-\widehat P_j^r(\xi) \right|\to0.
\end{aligned}
\end{equation*}
\par\medskip
\noindent\textbf{(vi)} Since $0\leq\widehat f^{\,r}(\xi)\leq L_f T$ on $[0,T]$, the clock approximation in \eqref{eq:prop7-clock-proxy-close} implies that, for all sufficiently large $r$,
\begin{equation*}
0\leq \bar f^r(\xi)\leq L_f T+1, \qquad 0\leq\xi\leq T.
\end{equation*}
By the local boundedness of $\lambda_j$,
$\int_0^{L_f T+1}\lambda_j(s)\,ds<\infty$.
Apply the functional law of large numbers for unit-rate Poisson processes \citep[Corollary~3.2.1]{Whitt2002Supplement} on the cumulative-intensity interval
$\left[0,\int_0^{L_f T+1}\lambda_j(s)\,ds\right]$.
We obtain
\[
\sup_{0\leq u\leq\int_0^{L_f T+1}\lambda_j(s)\,ds} \left|\frac{N_j^r(ru)}{r}-u\right| \xrightarrow{\mathbb P}0.
\]

Also, the local boundedness of the arrival rates gives, for all $x,y\in[0,L_f T+1]$,
\begin{equation}
\left| \int_0^x\lambda_j(s)\,ds -\int_0^y\lambda_j(s)\,ds \right| \leq \left\|\lambda_j\right\|_{L^\infty([0,L_f T+1])}|x-y|.
\label{eq:prop7-arrival-integral-local-Lipschitz}
\end{equation}
Both clocks take values in $[0,L_f T+1]$ for all sufficiently large $r$. Using the definitions of $\bar A_j^r$ and $\widehat A_j^r$, for $0\leq\xi\leq T$ we have
\begin{equation*}
{\medmuskip=1mu\thickmuskip=2mu\scriptspace=0.25pt\thinmuskip=2mu
\begin{aligned}
\bar A_j^r(\xi)-\widehat A_j^r(\xi) &=\frac{p_{j,L_j}^r(0)}r+\frac1rN_j^r\!\left(r\int_0^{\bar f^r(\xi)}\lambda_j(s)\,ds\right) -p_{j,L_j}(0)-\int_0^{\widehat f^{\,r}(\xi)}\lambda_j(s)\,ds\\
&=\frac1rN_j^r\!\left(r\int_0^{\bar f^r(\xi)}\lambda_j(s)\,ds\right)-\int_0^{\bar f^r(\xi)}\lambda_j(s)\,ds +\frac{p_{j,L_j}^r(0)}r-p_{j,L_j}(0)\\
&\quad+\int_0^{\bar f^r(\xi)}\lambda_j(s)\,ds-\int_0^{\widehat f^{\,r}(\xi)}\lambda_j(s)\,ds.
\end{aligned}}
\end{equation*}
Therefore, the Lipschitz bound for $A_j$ in \eqref{eq:prop7-arrival-integral-local-Lipschitz} gives
\begin{equation*}
{\medmuskip=1mu\thickmuskip=2mu\thinmuskip=1mu\scriptspace=0pt \begin{aligned}
&\sup_{0\leq\xi\leq T}|\bar A_j^r(\xi)-\widehat A_j^r(\xi)|\\
&\leq\sup_{0\leq\xi\leq T}\left|\frac1rN_j^r\!\left(r\!\int_0^{\bar f^r(\xi)}\mkern-13mu\lambda_j(s)\,ds\right)-\int_0^{\bar f^r(\xi)}\mkern-13mu\lambda_j(s)\,ds\right| +\left|\frac{p_{j,L_j}^r(0)}r-p_{j,L_j}(0)\right|\\
&\quad+\sup_{0\leq\xi\leq T}\Biggl|\int_0^{\bar f^r(\xi)}\lambda_j(s)\,ds-\int_0^{\widehat f^{\,r}(\xi)}\lambda_j(s)\,ds\Biggr|\\
&\leq\sup_{0\leq u\leq\int_0^{L_fT+1}\lambda_j(s)\,ds}\left|\frac{N_j^r(ru)}r-u\right| +\left|\frac{p_{j,L_j}^r(0)}r-p_{j,L_j}(0)\right| +\|\lambda_j\|_{L^\infty([0,L_fT+1])}\sup_{0\leq\xi\leq T}|\bar f^r(\xi)-\widehat f^{\,r}(\xi)|\xrightarrow{\mathbb P}0.
\end{aligned}}
\end{equation*}
Here the first term converges to zero by the Poisson estimate above, the second by the initial-condition convergence in Theorem~\ref{thm:moving-average-fluid-limit}, and the third by \textup{(v)}. This proves \eqref{eq:prop7-arrival-composition-close-vector}.
\end{proof}

\begin{lemma}
\label{lem:ctight-subsequential-limits}
Under the initial condition of Theorem~\ref{thm:moving-average-fluid-limit}, consider the smoothed $c\mu$-inspired policy in Definition~\ref{def:moving-average-policy}. Then the following statements hold.
\begin{enumerate}
  \item[\textup{(i)}] \textup{\textbf{Joint $C$-tightness.}} For every fixed $T<\infty$, the sequence
$\set{\left( \widehat P_1^r,\dots,\widehat P_J^r, \widehat f^{\,r}, \bar A_1^r,\dots,\bar A_J^r \right):r\ge1}$ is $C$-tight in $D([0,T],\mathbb R^{2J+1})$ equipped with the Skorokhod $J_1$ topology; i.e., their joint laws are tight in this space, and every subsequential weak limit is in $C([0,T],\mathbb R^{2J+1})$ with probability one.  Taking $T=1,2,\dots$ yields local joint $C$-tightness on $[0,\infty)$.
  \item[\textup{(ii)}] \textup{\textbf{Skorokhod representation of subsequential limits.}} From every scale subsequence, one can extract a further subsequence and construct copies having the same joint laws on a new probability space such that there exist continuous paths $\mathcal X_1,\dots,\mathcal X_J,f$ for which, almost surely, uniformly on every compact interval,
\begin{equation}
\left(\widehat P_1^r,\dots,\widehat P_J^r,\widehat f^{\,r},\bar A_1^r,\dots,\bar A_J^r\right) \longrightarrow\left(\mathcal X_1,\dots,\mathcal X_J,f,A_1\circ f,\dots,A_J\circ f\right).
\label{eq:prop7-subsequence-uoc}
\end{equation}
  \item[\textup{(iii)}] \textup{\textbf{Structural relations of subsequential limits.}} Every subsequential limit in \textup{(ii)} satisfies $\mathcal X_j(0)=0$ almost surely for $j=1,\dots,J$.
The function $\mathcal X_j$ is nonnegative and nondecreasing on $[0,\infty)$. The processes $\mathcal X_1,\dots,\mathcal X_J$ and the real-time clock $f$ satisfy the clock equation
\begin{equation}
f(\xi) =\sum_{j=1}^J(a_{\mathrm{pf}}s_j+a_{\mathrm{kv}}\tilde m_j)\mathcal X_j(\xi)+b\xi, \qquad \xi\geq0, \label{eq:prop7-limit-clock}
\end{equation}
the queue nonnegativity constraints
\begin{equation}
A_j\bigl(f(\xi)\bigr)-\mathcal X_j(\xi)\geq0, \qquad j=1,\dots,J,\quad \xi\geq0, \label{eq:prop7-limit-demand}
\end{equation}
and, for $0\leq\xi_1\leq\xi_2$, the memory constraint
\begin{equation}
\sum_{j=1}^J\tilde m_j \bigl(\mathcal X_j(\xi_2)-\mathcal X_j(\xi_1)\bigr) \leq M(\xi_2-\xi_1). \label{eq:prop7-limit-capacity}
\end{equation}
\end{enumerate}
\end{lemma}

\begin{proof}
\leavevmode\newline
\noindent\textbf{(i) Joint $C$-tightness.}
Fix $T<\infty$. To apply the Arzel\`a--Ascoli theorem, we first verify that $\widehat P_j^r$, $\widehat f^{\,r}$, and $\widehat A_j^r$, $j=1,\dots,J$, are uniformly bounded and Lipschitz continuous on $[0,T]$, with bounds and Lipschitz constants independent of $r$.

Lemma~\ref{lem:prelimit-uoc-estimates}\textup{(v)} gives, for all $0\leq\xi,\xi_1,\xi_2\leq T$,
\begin{equation*}
\begin{aligned}
0\leq\widehat P_j^r(\xi) &\leq\frac{MT}{\tilde m_j}, & \left|\widehat P_j^r(\xi_2)-\widehat P_j^r(\xi_1)\right| &\leq\frac{M}{\tilde m_j}|\xi_2-\xi_1|,\\
0\leq\widehat f^{\,r}(\xi) &\leq L_f T, & \left|\widehat f^{\,r}(\xi_2)-\widehat f^{\,r}(\xi_1)\right| &\leq L_f|\xi_2-\xi_1|.
\end{aligned}
\end{equation*}
Combining these estimates with the Lipschitz bound for $A_j$ in \eqref{eq:prop7-arrival-integral-local-Lipschitz} and the definition of $\widehat A_j^r$ in \eqref{eq:A-proxy} gives
\begin{equation*}
\begin{aligned}
\left|\widehat A_j^r(\xi_2)-\widehat A_j^r(\xi_1)\right|&=\left|\int_0^{\widehat f^{\,r}(\xi_2)}\lambda_j(s)\,ds -\int_0^{\widehat f^{\,r}(\xi_1)}\lambda_j(s)\,ds\right|\\
&\leq\left\|\lambda_j\right\|_{L^\infty([0,L_f T+1])} \left|\widehat f^{\,r}(\xi_2)-\widehat f^{\,r}(\xi_1)\right|\\
&\leq\left\|\lambda_j\right\|_{L^\infty([0,L_f T+1])}L_f|\xi_2-\xi_1|,\\
0\leq\widehat A_j^r(\xi) &=p_{j,L_j}(0) +\int_0^{\widehat f^{\,r}(\xi)}\lambda_j(s)\,ds\\
&\leq p_{j,L_j}(0) +\left\|\lambda_j\right\|_{L^\infty([0,L_f T+1])}\widehat f^{\,r}(\xi)\\
&\leq p_{j,L_j}(0) +\left\|\lambda_j\right\|_{L^\infty([0,L_f T+1])}L_f T.
\end{aligned}
\end{equation*}
Using the Euclidean norm on $\mathbb R^{2J+1}$, define the constants
\begin{equation*}
\begin{aligned}
\rho_T&:=\Biggl[\sum_{j=1}^J\Bigl(\frac{MT}{\tilde m_j}\Bigr)^2+(L_fT)^2 +\sum_{j=1}^J\Bigl(p_{j,L_j}(0)+\|\lambda_j\|_{L^\infty([0,L_fT+1])}L_fT\Bigr)^2\Biggr]^{1/2},\\
L_T&:=\Biggl[\sum_{j=1}^J\Bigl(\frac M{\tilde m_j}\Bigr)^2+L_f^2 +\sum_{j=1}^J\Bigl(\|\lambda_j\|_{L^\infty([0,L_fT+1])}L_f\Bigr)^2\Biggr]^{1/2}.
\end{aligned}
\end{equation*}
Both constants are finite by the local boundedness of the arrival rates and $J<\infty$.  Therefore, all sample paths of the continuous approximation vector
\begin{equation*}
\widehat{\mathbf Z}^{\,r}:=\left( \widehat P_1^r,\dots,\widehat P_J^r, \widehat f^{\,r}, \widehat A_1^r,\dots,\widehat A_J^r \right)
\end{equation*}
belong to the deterministic set
\begin{equation}
{\medmuskip=2mu\thickmuskip=2.5mu \mathcal K_T:=\Biggl\{\mathbf h\in C([0,T],\mathbb R^{2J+1}):\,\|\mathbf h\|_{\infty,T}\leq\rho_T,\|\mathbf h(\xi_2)-\mathbf h(\xi_1)\|\leq L_T|\xi_2-\xi_1|,\enspace \forall\,\xi_1,\xi_2\in[0,T]\Biggr\}.}
\label{eq:prop7-proxy-compact-set}
\end{equation}
The two inequalities defining $\mathcal K_T$ are preserved under uniform limits, so $\mathcal K_T$ is closed under uniform convergence.  It is also uniformly bounded and equicontinuous.  By the Arzel\`a--Ascoli theorem \cite[Theorem~7.25]{Rudin1976}, it is relatively compact in the uniform topology and is therefore compact there.  Moreover,
\begin{equation}
\mathbb P\left(\widehat{\mathbf Z}^{\,r}\in \mathcal K_T\right)=1, \qquad r\geq1. \label{eq:prop7-proxy-compact-support}
\end{equation}
Thus the family of laws $\{\mathcal L(\widehat{\mathbf Z}^{\,r}):r\geq1\}$ is tight in $C([0,T],\mathbb R^{2J+1})$ under the uniform topology.  Since the inclusion into $D([0,T],\mathbb R^{2J+1})$ with the $J_1$ topology is continuous, $\{\widehat{\mathbf Z}^{\,r}\}$ is $C$-tight in $D([0,T],\mathbb R^{2J+1})$.

Set
$\mathbf Z^r:=\left( \widehat P_1^r,\dots,\widehat P_J^r, \widehat f^{\,r}, \bar A_1^r,\dots,\bar A_J^r \right)$.
The vector $\widehat{\mathbf Z}^{\,r}$ approximates $\mathbf Z^r$ by replacing only the cumulative arrival coordinates $\bar A_j^r$ with their continuous approximations $\widehat A_j^r$. Hence,
\begin{equation*}
\begin{aligned}
\|\mathbf Z^r-\widehat{\mathbf Z}^{\,r}\|_{\infty,T}&:=\sup_{0\leq\xi\leq T}\|\mathbf Z^r(\xi)-\widehat{\mathbf Z}^{\,r}(\xi)\|\\
&=\sup_{0\leq\xi\leq T} \left(\sum_{j=1}^J\left|\bar A_j^r(\xi)-\widehat A_j^r(\xi)\right|^2\right)^{1/2}\\
&\leq\left(\sum_{j=1}^J \left[\sup_{0\leq\xi\leq T}\left|\bar A_j^r(\xi)-\widehat A_j^r(\xi)\right|\right]^2\right)^{1/2} \xrightarrow{\mathbb P}0,
\end{aligned}
\end{equation*}
where the convergence follows from the approximation of $\bar A_j^r$ by $\widehat A_j^r$ in \eqref{eq:prop7-arrival-composition-close-vector} and $J<\infty$. Define $\Delta\mathbf Z^r:=\mathbf Z^r-\widehat{\mathbf Z}^{\,r}$. $\|\Delta\mathbf Z^r\|_{\infty,T}\xrightarrow{\mathbb P}0$ implies that $\Delta\mathbf Z^r$ converges weakly to the zero process and hence is tight. Since the zero process is continuous, $\{\Delta\mathbf Z^r\}$ is $C$-tight. Both $\{\widehat{\mathbf Z}^{\,r}\}$ and $\{\Delta\mathbf Z^r\}$ are therefore $C$-tight. Corollary~VI.3.33(a) of \citet{JacodShiryaev2003} gives that $\{\mathbf Z^r=\widehat{\mathbf Z}^{\,r}+\Delta\mathbf Z^r\}$ is jointly $C$-tight in $D([0,T],\mathbb R^{2J+1})$.

Finally, apply the preceding argument successively on $[0,m]$ for $m=1,2,\dots$, extract nested subsequences, and take a diagonal subsequence. This yields local joint $C$-tightness on $[0,\infty)$.
\par\medskip
\noindent\textbf{(ii) Skorokhod representation of subsequential limits.}
We prove the almost-sure convergence on compact intervals in \eqref{eq:prop7-subsequence-uoc} by first identifying the joint weak limit and then applying the Skorokhod representation theorem. Since $\{(\widehat P_1^r,\dots,\widehat P_J^r,\widehat f^{\,r})\}$ is jointly $C$-tight, every scale subsequence has a further subsequence, still indexed by $r$, along which
\begin{equation*}
\left(\widehat P_1^r,\dots,\widehat P_J^r,\widehat f^{\,r}\right)\Rightarrow\left(\mathcal X_1,\dots,\mathcal X_J,f\right)
\end{equation*}
in the local $J_1$ topology, with continuous limit paths. The clock limit $f$ is nonnegative because each $\widehat f^{\,r}$ is nonnegative. Each $A_j$ is locally Lipschitz, so the composition map $g\mapsto A_j\circ g$ is continuous under uniform convergence on compact intervals. Since $f$ is continuous, local $J_1$ convergence to $f$ implies uniform convergence on compact intervals. The continuous mapping theorem and the definition of $\widehat A_j^r$ in \eqref{eq:A-proxy} therefore give
\begin{equation*}
\widehat{\mathbf Z}^{\,r}\Rightarrow\left(\mathcal X_1,\dots,\mathcal X_J,f,A_1\circ f,\dots,A_J\circ f\right).
\end{equation*}
For every $T<\infty$, the error estimate in \textup{(i)} gives $\|\Delta\mathbf Z^r\|_{\infty,T}\xrightarrow{\mathbb P}0$. Applying Lemma~VI.3.31 of \citet{JacodShiryaev2003} to the joint vector and the error $\Delta\mathbf Z^r$ yields
\begin{equation*}
\mathbf Z^r=\widehat{\mathbf Z}^{\,r}+\Delta\mathbf Z^r \Rightarrow\left(\mathcal X_1,\dots,\mathcal X_J,f,A_1\circ f,\dots,A_J\circ f\right).
\end{equation*}

Since $D([0,\infty),\mathbb R^{2J+1})$ with the local $J_1$ topology is Polish, the Skorokhod representation theorem \cite[Theorem~6.7]{Billingsley1999} provides copies on a common probability space, preserving the joint law of the coordinates for each $r$, that converge almost surely in local $J_1$ to a copy of the displayed limit. Since that limit has continuous paths, the convergence is uniform on compact intervals \cite[p.~124]{Billingsley1999}. Retaining the original notation for these copies gives \eqref{eq:prop7-subsequence-uoc}.
\par\medskip
\noindent\textbf{(iii) Structural relations of subsequential limits.}
Work on a probability-one event where the uniform convergence in \eqref{eq:prop7-subsequence-uoc} and the prelimit path properties used below hold for the same-law copies.

We first verify the path properties and the clock equation. For every $r$, the path $\widehat P_j^r$ is nonnegative and nondecreasing and satisfies $\widehat P_j^r(0)=0$.  Uniform convergence on compact intervals in \eqref{eq:prop7-subsequence-uoc} therefore gives $\mathcal X_j(0)=0$.

Recall the definition of $\widehat f^{\,r}$ in \eqref{eq:clock-proxy}:
$\widehat f^{\,r}(\xi) =\sum_{j=1}^J(a_{\mathrm{pf}}s_j+a_{\mathrm{kv}}\tilde m_j) \widehat P_j^r(\xi)+b\xi.$
For every $T<\infty$, the uniform convergence in \eqref{eq:prop7-subsequence-uoc} yields
\begin{equation*}
\begin{aligned}
&\sup_{0\leq\xi\leq T}\Biggl|f(\xi)-\sum_{j=1}^J(a_{\mathrm{pf}}s_j+a_{\mathrm{kv}}\tilde m_j)\mathcal X_j(\xi)-b\xi\Biggr|\\
&\leq\sup_{0\leq\xi\leq T}|f(\xi)-\widehat f^{\,r}(\xi)| +\sum_{j=1}^J(a_{\mathrm{pf}}s_j+a_{\mathrm{kv}}\tilde m_j)\sup_{0\leq\xi\leq T}|\widehat P_j^r(\xi)-\mathcal X_j(\xi)|\to0.
\end{aligned}
\end{equation*}
Because $T$ is arbitrary, \eqref{eq:prop7-limit-clock} follows.

\medskip
\noindent\textbf{Queue nonnegativity constraints.}
For the original prelimit processes, the cumulative planned-admission bound \eqref{eq:prelimit-postplan-demand} gives
\begin{equation*}
\widehat P_j^r(n/r)=\frac{P_j^r(n)}{r}=\frac{1}{r}\sum_{u=0}^{n-1}y_j^r(u)\leq\frac{1}{r}\sum_{u=0}^{n}y_j^r(u) \leq\frac{1}{r}\left[p_{j,L_j}^{r}(0)+\sum_{u=0}^{n-1}\eta_j^r(u)\right]=\bar A_j^r(n/r).
\end{equation*}
Hence,
\begin{equation}
0\leq \bar A_j^r(n/r)-\widehat P_j^r(n/r), \qquad n\in\mathbb{Z}_+. \label{eq:prop7-unplanned-demand-nonnegative}
\end{equation}
This inequality is a joint path property of $(\bar A_j^r,\widehat P_j^r)$ and therefore also holds almost surely for their same-law copies.

Fix $T<\infty$.  Since $|\lfloor r\xi\rfloor/r-\xi|\leq1/r$, the uniform convergence in \eqref{eq:prop7-subsequence-uoc} and the continuity of $A_j\circ f$ and $\mathcal X_j$ imply, uniformly over $0\leq\xi\leq T$,
\begin{equation*}
\bar A_j^r\bigl(\lfloor r\xi\rfloor/r\bigr) -\widehat P_j^r\bigl(\lfloor r\xi\rfloor/r\bigr) \longrightarrow A_j\bigl(f(\xi)\bigr)-\mathcal X_j(\xi).
\end{equation*}
The left-hand side is nonnegative by \eqref{eq:prop7-unplanned-demand-nonnegative}; its limit is therefore nonnegative.  Since $T$ is arbitrary, \eqref{eq:prop7-limit-demand} follows.

\medskip
We next verify the memory constraint \eqref{eq:prop7-limit-capacity}. Summing the planning-budget bound \eqref{eq:prop7-planning-budget} from $n_1$ to $n_2-1$ gives, for integers $0\leq n_1\leq n_2$,
\begin{equation*}
\sum_{j=1}^J\tilde m_j(P_j^r(n_2)-P_j^r(n_1))=\sum_{j=1}^J\tilde m_j\sum_{t=n_1}^{n_2-1}y_j^r(t) =\sum_{t=n_1}^{n_2-1}\sum_{j=1}^J\tilde m_jy_j^r(t) \leq\frac{w_rM}{w_r+H_{\mathrm{kv}}}(n_2-n_1).
\end{equation*}
After division by $r$ and use of $\widehat P_j^r(n/r)=P_j^r(n)/r$, this becomes
\begin{equation*}
\sum_{j=1}^J\tilde m_j\left[\widehat P_j^r\left(\frac{n_2}{r}\right)-\widehat P_j^r\left(\frac{n_1}{r}\right)\right]
\leq\frac{w_rM}{w_r+H_{\mathrm{kv}}}\frac{n_2-n_1}{r}.
\end{equation*}
This is again a path property retained by the same-law copies.  Given $0\leq\xi_1\leq\xi_2$, take $n_i=\lfloor r\xi_i\rfloor$.  Then the uniform convergence in \eqref{eq:prop7-subsequence-uoc} and the continuity of $\mathcal X_j$ give
\begin{equation*}
\begin{gathered}
\widehat P_j^r(n_i/r)\longrightarrow \mathcal X_j(\xi_i), \quad i=1,2,
\text{while}\quad \frac{n_2-n_1}{r}\longrightarrow\xi_2-\xi_1, \quad \frac{w_r}{w_r+H_{\mathrm{kv}}}\longrightarrow1.
\end{gathered}
\end{equation*}
Passing to the limit in the preceding inequality proves \eqref{eq:prop7-limit-capacity}.
\end{proof}

\begin{lemma}[Residual representation and uniqueness]
\label{lem:residual-representation-uniqueness}
\leavevmode\newline
Under the initial condition of Theorem~\ref{thm:moving-average-fluid-limit}, every subsequential limit obtained in Lemma~\ref{lem:ctight-subsequential-limits} satisfies the following properties almost surely.
\begin{enumerate}
  \item[\textup{(i)}] The clock $f$ is a strictly increasing continuous bijection from $[0,\infty)$ onto itself.  Define
\begin{equation}
R_0(s):=s,\qquad R_j\bigl(f(\xi)\bigr):=f(\xi)-\sum_{\ell=1}^j\mu_\ell^{-1}\mathcal X_\ell(\xi),\enspace j=1,\dots,J.
\label{eq:prop7-candidate-residual}
\end{equation}
Then each $R_j$, $j=0,\dots,J$, is continuous, nonnegative, and nondecreasing, with $R_j(0)=0$.  Moreover, for $j=1,\dots,J$,
\begin{equation}
R_j(s) \geq R_{j-1}(s)-\mu_j^{-1}A_j(s), \qquad s\geq0. \label{eq:prop7-residual-obstacle}
\end{equation}
  \item[\textup{(ii)}] For every $j=1,\dots,J$ and $s\geq0$,
\begin{equation}
\int_0^s \left[ R_j(u)-R_{j-1}(u)+\mu_j^{-1}A_j(u) \right]\,dR_j(u)=0. \label{eq:prop7-residual-complementarity}
\end{equation}
  \item[\textup{(iii)}] The following identities hold:
\begin{equation}
R_j=R_j^*,\quad f=f^*,\quad\mathcal X_j=\mathcal X_j^*,
\label{eq:prop7-unique-identification}
\end{equation}
where $j=0,\dots,J$ for $R_j$ and $j=1,\dots,J$ for $\mathcal X_j$. Consequently, every subsequential limit is the same deterministic target path.
\end{enumerate}
\end{lemma}

\begin{proof}
\noindent\textbf{(i)}
Work on a probability-one event on which the conclusions of Lemma~\ref{lem:ctight-subsequential-limits} and all prelimit path properties used below hold for the same-law copies. For $0\leq\xi_1<\xi_2$, the clock equation \eqref{eq:prop7-limit-clock} and the monotonicity of the $\mathcal X_j$ give
\begin{equation*}
{\medmuskip=1mu\thickmuskip=2mu\scriptspace=0.25pt
f(\xi_2) -f(\xi_1)=\sum_{j=1}^J(a_{\mathrm{pf}} s_j+a_{\mathrm{kv}}\tilde m_j) \bigl(\mathcal X_j(\xi_2)-\mathcal X_j(\xi_1)\bigr)+b(\xi_2-\xi_1) \geq b(\xi_2-\xi_1) >0,
}
\end{equation*}
where the last inequality uses the standing assumption $b>0$.  Together with the continuity of $f$, $f(0)=0$, and $f(\xi)\geq b\xi\to\infty$, this shows that $f$ is a continuous bijection onto $[0,\infty)$ with a continuous inverse.

The clock equation \eqref{eq:prop7-limit-clock} and the definition of $R_j$ in \eqref{eq:prop7-candidate-residual} give
\begin{equation}
{\medmuskip=1mu\thickmuskip=2mu\scriptspace=0.25pt
\begin{aligned}
R_j\bigl(f(\xi)\bigr)&=\sum_{\ell=j+1}^J (a_{\mathrm{pf}} s_\ell+a_{\mathrm{kv}}\tilde m_\ell)\mathcal X_\ell(\xi)+b\xi -\frac bM\sum_{\ell=1}^j\tilde m_\ell \mathcal X_\ell(\xi)\\
&=\sum_{\ell=j+1}^J\left(a_{\mathrm{pf}}s_\ell+a_{\mathrm{kv}}\tilde m_\ell+\frac bM\tilde m_\ell\right)\mathcal X_\ell(\xi)+b\xi -\frac bM\sum_{\ell=j+1}^J\tilde m_\ell\mathcal X_\ell(\xi) -\frac bM\sum_{\ell=1}^j\tilde m_\ell\mathcal X_\ell(\xi)\\
&=\sum_{\ell=j+1}^J\mu_\ell^{-1}\mathcal X_\ell(\xi)+\frac bM\left( M\xi-\sum_{\ell=1}^J\tilde m_\ell \mathcal X_\ell(\xi) \right).
\end{aligned}}
\label{eq:prop7-residual-monotone-form}
\end{equation}
By Lemma~\ref{lem:ctight-subsequential-limits}, each $\mathcal X_\ell$ is continuous, nonnegative, and nondecreasing, with $\mathcal X_\ell(0)=0$. The function
$\xi\longmapsto M\xi-\sum_{\ell=1}^J\tilde m_\ell \mathcal X_\ell(\xi)$
is also continuous and equals zero at $\xi=0$. The memory constraint \eqref{eq:prop7-limit-capacity} shows that it is nonnegative and nondecreasing. Thus the representation \eqref{eq:prop7-residual-monotone-form} gives these same properties for $R_j\circ f$. Since $f^{-1}$ is continuous and increasing with $f^{-1}(0)=0$, they also hold for $R_j=(R_j\circ f)\circ f^{-1}$.

Finally, the definition of $R_j$ in \eqref{eq:prop7-candidate-residual} and the queue nonnegativity constraint \eqref{eq:prop7-limit-demand} give
\[
R_j\bigl(f(\xi)\bigr) -R_{j-1}\bigl(f(\xi)\bigr) +\mu_j^{-1}A_j\bigl(f(\xi)\bigr)=\mu_j^{-1} \bigl(A_j\bigl(f(\xi)\bigr)-\mathcal X_j(\xi)\bigr)\geq0.
\]
Setting $\xi=f^{-1}(s)$ proves \eqref{eq:prop7-residual-obstacle}.
\par\medskip
\noindent\textbf{(ii)}
We prove the complementarity relation \eqref{eq:prop7-residual-complementarity} by showing that $R_j$ is constant on intervals where the class-$j$ waiting-queue length is positive. At every GPU step $t$, the cumulative planned-admission bound \eqref{eq:prelimit-postplan-demand} gives
\begin{equation*}
p_{j,L_j}^r(0)+\sum_{u<t}\eta_j^r(u)-\sum_{u<t}y_j^r(u)\geq y_j^r(t)\geq0.
\end{equation*}
The planning-budget bound \eqref{eq:prop7-planning-budget} gives
\begin{equation*}
\frac{w_rM}{w_r+H_{\mathrm{kv}}}-\sum_{\ell=1}^{j-1}\tilde m_\ell y_\ell^r(t)\geq\sum_{\ell=j}^J\tilde m_\ell y_\ell^r(t)\geq0.
\end{equation*}
Thus the positive parts in the planning rule \eqref{eq:planned-admission-policy} can be removed:
\begin{equation*}
y_j^r(t)=\min\Biggl\{p_{j,L_j}^r(0)+\sum_{u<t}\eta_j^r(u)-\sum_{u<t}y_j^r(u),\tilde m_j^{-1}\left[\frac{w_rM}{w_r+H_{\mathrm{kv}}}-\sum_{\ell=1}^{j-1}\tilde m_\ell y_\ell^r(t)\right]\Biggr\}.
\end{equation*}
Applying the identity
$\bigl(a-\min\{a,b\}\bigr)\bigl(b-\min\{a,b\}\bigr)=0$
(valid for $a,b\geq0$) to the two arguments of the minimum
and multiplying by $\tilde m_j$ gives
\begin{equation}
\left[p_{j,L_j}^r(0)+\sum_{u<t}\eta_j^r(u)-\sum_{u\leq t}y_j^r(u)\right]\left[\frac{w_rM}{w_r+H_{\mathrm{kv}}}-\sum_{\ell=1}^j\tilde m_\ell y_\ell^r(t)\right]=0.
\label{eq:prop7-discrete-complementarity}
\end{equation}

The first factor in \eqref{eq:prop7-discrete-complementarity} is the number of class-$j$ requests that remain without an admission plan after planning at step $t$. After fluid scaling, it becomes
\begin{equation*}
\frac1r\left[p_{j,L_j}^r(0)+\sum_{u<t}\eta_j^r(u)-\sum_{u\leq t}y_j^r(u)\right]=\bar A_j^r(t/r)-\widehat P_j^r((t+1)/r).
\end{equation*}
For every $T<\infty$,
\begin{equation}
\begin{aligned}
&\max_{0\leq t\leq\lfloor rT\rfloor}\Biggl|\frac{1}{r}\Biggl[p_{j,L_j}^r(0)+\sum_{u<t}\eta_j^r(u)-\sum_{u\leq t}y_j^r(u)\Biggr] -\Biggl[A_j\left(f\left(\frac{t}{r}\right)\right)-\mathcal X_j\left(\frac{t}{r}\right)\Biggr]\Biggr|\\
&=\max_{0\leq t\leq\lfloor rT\rfloor}\Biggl|\bar A_j^r\left(\frac{t}{r}\right)-A_j\left(f\left(\frac{t}{r}\right)\right) -\Biggl[\widehat P_j^r\left(\frac{t+1}{r}\right)-\mathcal X_j\left(\frac{t+1}{r}\right)\Biggr] -\Biggl[\mathcal X_j\left(\frac{t+1}{r}\right)-\mathcal X_j\left(\frac{t}{r}\right)\Biggr]\Biggr|\\
&\leq\sup_{0\leq\xi\leq T}|\bar A_j^r(\xi)-A_j(f(\xi))| +\sup_{0\leq\xi\leq T+1}|\widehat P_j^r(\xi)-\mathcal X_j(\xi)| +\sup_{\substack{0\leq x,y\leq T+1\\|x-y|\leq\frac{1}{r}}}|\mathcal X_j(x)-\mathcal X_j(y)|\longrightarrow0.
\end{aligned}
\label{eq:prop7-postplan-gap-limit}
\end{equation}
Here the first two terms vanish by the uniform convergence in \eqref{eq:prop7-subsequence-uoc}, and the last by uniform continuity of $\mathcal X_j$ on $[0,T+1]$. Let $[u,v]$, with $u<v$, be a compact interval contained in
\begin{equation*}
\left\{\xi\geq0:A_j\bigl(f(\xi)\bigr)-\mathcal X_j(\xi)>0\right\}.
\end{equation*}
Since $A_j\circ f-\mathcal X_j$ is continuous and strictly positive on the compact interval $[u,v]$, there exists $\delta>0$ such that
\begin{equation*}
A_j\bigl(f(\xi)\bigr)-\mathcal X_j(\xi)\geq2\delta,\qquad \xi\in[u,v].
\end{equation*}
By \eqref{eq:prop7-postplan-gap-limit} with $T=v$, choose $r_0$ so that the maximum error in \eqref{eq:prop7-postplan-gap-limit} is less than $\delta$ and $\lceil ru\rceil<\lfloor rv\rfloor$ for every $r\geq r_0$. Then, for $t=\lceil ru\rceil,\ldots,\lfloor rv\rfloor-1$,
\begin{equation*}
\bar A_j^r(t/r)-\widehat P_j^r((t+1)/r)>A_j(f(t/r))-\mathcal X_j(t/r)-\delta \geq\delta >0.
\end{equation*}
Thus the first factor in \eqref{eq:prop7-discrete-complementarity} is strictly positive. Its second factor must be zero, giving
\begin{equation}
\sum_{\ell=1}^j\tilde m_\ell y_\ell^r(t) =\frac{w_rM}{w_r+H_{\mathrm{kv}}}. \label{eq:prop7-full-prefix-on-gap}
\end{equation}
The planning-budget bound \eqref{eq:prop7-planning-budget} now implies
\begin{equation*}
0\leq\sum_{\ell=j+1}^J\tilde m_\ell y_\ell^r(t)=\sum_{\ell=1}^J\tilde m_\ell y_\ell^r(t)-\frac{w_rM}{w_r+H_{\mathrm{kv}}}\leq0.
\end{equation*}
Since $\tilde m_\ell>0$ and $y_\ell^r(t)\geq0$, every term in the first sum is zero, so
$y_\ell^r(t)=0$ for $\ell>j$.

For $\ell>j$, summing over $t=\lceil ru\rceil,\ldots,\lfloor rv\rfloor-1$ and using the definition of $P_\ell^r$ gives
\begin{equation*}
0=\frac1r\sum_{t=\lceil ru\rceil}^{\lfloor rv\rfloor-1}y_\ell^r(t)=\frac{P_\ell^r(\lfloor rv\rfloor)-P_\ell^r(\lceil ru\rceil)}r =\widehat P_\ell^r(\lfloor rv\rfloor/r)-\widehat P_\ell^r(\lceil ru\rceil/r).
\end{equation*}
Here we used $\widehat P_\ell^r(n/r)=P_\ell^r(n)/r$ at grid points. Since $y_\ell^r(t)=0$ for $\ell>j$, summing both sides of the full-budget equality \eqref{eq:prop7-full-prefix-on-gap} over the same GPU steps and dividing by $r$ gives
\begin{equation*}
\begin{aligned}
\frac{w_rM}{w_r+H_{\mathrm{kv}}}\frac{\lfloor rv\rfloor-\lceil ru\rceil}{r} &=\frac1r\sum_{t=\lceil ru\rceil}^{\lfloor rv\rfloor-1}\sum_{\ell=1}^J\tilde m_\ell y_\ell^r(t)\\
&=\frac1r\sum_{\ell=1}^J\tilde m_\ell\sum_{t=\lceil ru\rceil}^{\lfloor rv\rfloor-1}y_\ell^r(t)\\
&=\sum_{\ell=1}^J\tilde m_\ell\frac{P_\ell^r(\lfloor rv\rfloor)-P_\ell^r(\lceil ru\rceil)}r\\
&=\sum_{\ell=1}^J\tilde m_\ell\left[\widehat P_\ell^r(\lfloor rv\rfloor/r)-\widehat P_\ell^r(\lceil ru\rceil/r)\right].
\end{aligned}
\end{equation*}
The convergence of the grid endpoints to $u$ and $v$, the limit $w_r/(w_r+H_{\mathrm{kv}})\to1$, the uniform convergence of $\widehat P_\ell^r$ to $\mathcal X_\ell$ in \eqref{eq:prop7-subsequence-uoc}, and the continuity of $\mathcal X_\ell$ together give
\begin{equation*}
\mathcal X_\ell(v) -\mathcal X_\ell(u) =0\quad(\ell>j),\qquad\text{and}\qquad\sum_{\ell=1}^J\tilde m_\ell\bigl(\mathcal X_\ell(v)-\mathcal X_\ell(u)\bigr) =M(v-u).
\end{equation*}
Substituting these equalities into the representation of $R_j\circ f$ in \eqref{eq:prop7-residual-monotone-form} gives
\begin{equation*}
R_j\bigl(f(v)\bigr) -R_j\bigl(f(u)\bigr)=\sum_{\ell=j+1}^J\mu_\ell^{-1} \bigl(\mathcal X_\ell(v)-\mathcal X_\ell(u)\bigr) +\frac bM\left[ M(v-u)-\sum_{\ell=1}^J\tilde m_\ell \bigl(\mathcal X_\ell(v)-\mathcal X_\ell(u)\bigr) \right]=0.
\end{equation*}

Since $[u,v]$ is arbitrary and $f$ is a homeomorphism, the preceding equality implies that the Lebesgue--Stieltjes measure $dR_j$ satisfies
\begin{equation*}
dR_j\!\left(\left\{s\geq0:A_j(s)-\mathcal X_j\bigl(f^{-1}(s)\bigr)>0\right\}\right)=0.
\end{equation*}
The definition of $R_j$ in \eqref{eq:prop7-candidate-residual} gives
\begin{equation*}
R_j(s) -R_{j-1}(s) +\mu_j^{-1}A_j(s)=\mu_j^{-1} \left[A_j(s)-\mathcal X_j\bigl(f^{-1}(s)\bigr)\right] \geq0.
\end{equation*}
The integrand in \eqref{eq:prop7-residual-complementarity} is therefore zero $dR_j$-almost everywhere, so for every $s\geq0$,
\begin{equation*}
\int_0^s\left[R_j(u)-R_{j-1}(u)+\mu_j^{-1}A_j(u)\right]\,dR_j(u)=0.
\end{equation*}

\par\medskip
\noindent\textbf{(iii)}
For fixed $j$, the function $R_{j-1}-\mu_j^{-1}A_j$ is continuous.
By \textup{(i)}, $R_j$ is continuous and nondecreasing, with $R_j(0)=0$, and
\begin{equation*}
R_j(s)\geq R_{j-1}(s)-\mu_j^{-1}A_j(s),\qquad s\geq0.
\end{equation*}
By \textup{(ii)}, for every $s\geq0$,
\begin{equation*}
\int_0^s\left[R_j(u)-R_{j-1}(u)+\mu_j^{-1}A_j(u)\right]\,dR_j(u)=0.
\end{equation*}
The one-dimensional reflection characterization \cite[Theorem~14.2.3]{Whitt2002}, applied on each compact interval, therefore yields
\begin{equation*}
R_j(s) =\sup_{0\leq u\leq s} \left[ R_{j-1}(u)-\mu_j^{-1}A_j(u) \right]^+.
\end{equation*}
The base case is $R_0=R_0^*$ because both are the identity function.  If $R_{j-1}=R_{j-1}^*$, the preceding display is exactly the recursion defining $R_j^*$ in \eqref{eq:target-residual}.  Induction gives
\begin{equation}
R_j(s)=R_j^*(s), \qquad j=0,\dots,J,\quad s\geq0. \label{eq:prop7-residual-identified}
\end{equation}

Substituting the definitions of $R_{j-1}$ and $R_j$ from \eqref{eq:prop7-candidate-residual} gives
\begin{equation*}
{\medmuskip=2mu\thickmuskip=2.5mu
R_{j-1}\bigl(f(\xi)\bigr) -R_j\bigl(f(\xi)\bigr)=\left[f(\xi)-\sum_{\ell=1}^{j-1}\mu_\ell^{-1}\mathcal X_\ell(\xi)\right]-\left[f(\xi)-\sum_{\ell=1}^{j}\mu_\ell^{-1}\mathcal X_\ell(\xi)\right] =\mu_j^{-1}\mathcal X_j(\xi).
}
\end{equation*}
Multiplying by $\mu_j$ yields
\begin{equation}
\mathcal X_j(\xi) =\mu_j\left[ R_{j-1}\bigl(f(\xi)\bigr)-R_j\bigl(f(\xi)\bigr) \right]. \label{eq:prop7-candidate-X-residual}
\end{equation}
Taking $j=J$ in the representation of $R_j\circ f$ in \eqref{eq:prop7-residual-monotone-form} gives
\begin{equation*}
\xi =\frac1M\sum_{j=1}^J\tilde m_j\mathcal X_j(\xi) +\frac1bR_J\bigl(f(\xi)\bigr).
\end{equation*}
The definition of $F$ in \eqref{eq:target-F}, the identification $R_j^*=R_j$ in \eqref{eq:prop7-residual-identified}, and the representation of $\mathcal X_j$ in \eqref{eq:prop7-candidate-X-residual} give
\begin{equation}
\begin{aligned}
F\bigl(f(\xi)\bigr)&=\frac1M\sum_{j=1}^J\tilde m_j\mu_j\left[R_{j-1}^*\bigl(f(\xi)\bigr)-R_j^*\bigl(f(\xi)\bigr)\right]+\frac1bR_J^*\bigl(f(\xi)\bigr)\\
&=\frac1M\sum_{j=1}^J\tilde m_j\mu_j\left[R_{j-1}\bigl(f(\xi)\bigr)-R_j\bigl(f(\xi)\bigr)\right]+\frac1bR_J\bigl(f(\xi)\bigr)\\
&=\frac1M\sum_{j=1}^J\tilde m_j\mathcal X_j(\xi) +\frac1bR_J\bigl(f(\xi)\bigr)\\
&=\xi,\qquad\xi\geq0.
\end{aligned}
\label{eq:prop7-F-f}
\end{equation}

By \textup{(i)}, $f$ is a homeomorphism of $[0,\infty)$ onto itself.  Hence \eqref{eq:prop7-F-f} implies $F=f^{-1}$.  In particular, the generalized inverse in \eqref{eq:target-clock} is the ordinary inverse, and therefore
$f^*=F^{-1}=f$.
Finally, comparing the representation of $\mathcal X_j$ in \eqref{eq:prop7-candidate-X-residual} with the definition of $\mathcal X_j^*$ in \eqref{eq:target-admission} gives $\mathcal X_j=\mathcal X_j^*$ for every $j$.
\end{proof}

\begin{proof}[Proof of Theorem~\ref{thm:moving-average-fluid-limit}]
Fix $T<\infty$.  By Lemma~\ref{lem:ctight-subsequential-limits}, every scale subsequence has a further weakly convergent subsequence.  Lemmas~ \ref{lem:ctight-subsequential-limits} and \ref{lem:residual-representation-uniqueness} show that the $\bigl(\widehat P_1^r,\dots,\widehat P_J^r, \widehat f^{\,r}\bigr)$ coordinates of every such limit equal the deterministic path $\bigl(\mathcal X_1^*,\dots,\mathcal X_J^*,f^*\bigr)$.  Hence
\begin{equation*}
\left(\widehat P_1^r,\dots,\widehat P_J^r,\widehat f^{\,r}\right)\xrightarrow[\mathrm{local}\ J_1]{d}\left(\mathcal X_1^*,\dots,\mathcal X_J^*,f^*\right).
\end{equation*}
Weak convergence to a deterministic path is convergence in probability, and at a continuous limit local $J_1$ convergence strengthens to uniform convergence on compact intervals \cite[p.~124]{Billingsley1999}.  Therefore
\begin{equation}
\sup_{0\leq\xi\leq T} \left|\widehat f^{\,r}(\xi)-f^*(\xi)\right|+\sum_{j=1}^J\sup_{0\leq\xi\leq T} \left|\widehat P_j^r(\xi)-\mathcal X_j^*(\xi)\right| \xrightarrow{\mathbb P}0.
\label{eq:prop7-proxy-target-convergence}
\end{equation}

By the clock approximation in \eqref{eq:prop7-clock-proxy-close} and the uniform convergence of $\widehat f^{\,r}$ to $f^*$ in \eqref{eq:prop7-proxy-target-convergence},
\begin{equation*}
\sup_{0\leq\xi\leq T} \left|\bar f^r(\xi)-f^*(\xi)\right| \leq \sup_{0\leq\xi\leq T} \left|\bar f^r(\xi)-\widehat f^{\,r}(\xi)\right| +\sup_{0\leq\xi\leq T} \left|\widehat f^{\,r}(\xi)-f^*(\xi)\right| \xrightarrow{\mathbb P}0.
\end{equation*}
This proves \eqref{eq:prelimit-clock-convergence}.

For each class $j$, the bound on the difference between fluid-scaled cumulative planned admissions and fluid-scaled cumulative admissions in \eqref{eq:prop7-planned-actual}, the interpolation error bound \eqref{eq:prop7-P-proxy-bounds}, and the uniform convergence of $\widehat P_j^r$ to $\mathcal X_j^*$ in \eqref{eq:prop7-proxy-target-convergence} give
\begin{equation*}
\begin{aligned}
&\sup_{0\leq\xi\leq T} \left| \frac{1}{r}\mathcal U_{j,L_j}^{r}(\lfloor r\xi\rfloor)-\mathcal X_j^*(\xi) \right|\\
&\leq \sup_{0\leq\xi\leq T} \left| \frac{1}{r}\mathcal U_{j,L_j}^{r}(\lfloor r\xi\rfloor) -\bar P_j^r(\xi) \right| + \sup_{0\leq\xi\leq T} \left| \bar P_j^r(\xi)-\widehat P_j^r(\xi) \right| +\sup_{0\leq\xi\leq T} \left|\widehat P_j^r(\xi)-\mathcal X_j^*(\xi)\right| \xrightarrow{\mathbb P}0.
\end{aligned}
\end{equation*}
This is the $k=L_j$ case of \eqref{eq:prelimit-U-convergence}. For $k=1,\dots,L_j$, the prefill-stage delay bound \eqref{eq:prop7-fixed-shift-U} gives
\begin{equation*}
\sup_{0\leq\xi\leq T} \left| \frac{1}{r}\mathcal U_{j,k}^{r}(\lfloor r\xi\rfloor)-\mathcal X_j^*(\xi) \right|\leq \sup_{0\leq\xi\leq T} \left| \frac{ \mathcal U_{j,k}^{r}(\lfloor r\xi\rfloor) -\mathcal U_{j,L_j}^{r}(\lfloor r\xi\rfloor) }{r} \right| + \sup_{0\leq\xi\leq T} \left| \frac{1}{r}\mathcal U_{j,L_j}^{r}(\lfloor r\xi\rfloor)-\mathcal X_j^*(\xi) \right| \xrightarrow{\mathbb P}0.
\end{equation*}
Similarly, for $k=1,\dots,o_j-1$, the decode-stage delay bound \eqref{eq:prop7-fixed-shift-V} gives
\begin{equation*}
\sup_{0\leq\xi\leq T} \left| \frac{1}{r}\mathcal V_{j,k}^{r}(\lfloor r\xi\rfloor)-\mathcal X_j^*(\xi) \right|\leq \sup_{0\leq\xi\leq T} \left| \frac{ \mathcal V_{j,k}^{r}(\lfloor r\xi\rfloor) -\mathcal U_{j,L_j}^{r}(\lfloor r\xi\rfloor) }{r} \right| + \sup_{0\leq\xi\leq T} \left| \frac{1}{r}\mathcal U_{j,L_j}^{r}(\lfloor r\xi\rfloor)-\mathcal X_j^*(\xi) \right| \xrightarrow{\mathbb P}0.
\end{equation*}
This proves \eqref{eq:prelimit-U-convergence} and \eqref{eq:prelimit-V-convergence}.

Finally, we prove convergence of the system states. The fluid states satisfy $p_{j,k}^{\pi^*}(\xi)=0$ for $k=1,\dots,L_j-1$ and $d_{j,k}^{\pi^*}(\xi)=0$ for $k=1,\dots,o_j-1$. By \eqref{eq:vanishing-intermediate-states}, we therefore obtain
\[
\begin{aligned}
\sup_{0\leq\xi\leq T}\left|\bar p_{j,k}^r(\xi)-p_{j,k}^{\pi^*}(\xi)\right|&\xrightarrow{\mathbb P}0,\quad k=1,\dots,L_j-1,\\
\sup_{0\leq\xi\leq T}\left|\bar d_{j,k}^r(\xi)-d_{j,k}^{\pi^*}(\xi)\right|&\xrightarrow{\mathbb P}0,\quad k=1,\dots,o_j-1.
\end{aligned}
\]
It remains to prove \eqref{eq:prelimit-p-convergence} for $k=L_j$. By definition,
\begin{equation*}
\frac1r p_{j,L_j}^r(\lfloor r\xi\rfloor)=\bar A_j^r(\xi)-\frac1r\mathcal U_{j,L_j}^r(\lfloor r\xi\rfloor)\quad \text{and} \quad p_{j,L_j}^{\pi^*}(\xi)=A_j(f^*(\xi))-\mathcal X_j^*(\xi).
\end{equation*}
By \eqref{eq:prop7-subsequence-uoc} and Lemma~\ref{lem:residual-representation-uniqueness}, every subsequential limit of $\bar A_j^r$ equals the deterministic path $A_j\circ f^*$. Hence,
\begin{equation*}
\sup_{0\leq\xi\leq T}\left|\bar A_j^r(\xi)-A_j(f^*(\xi))\right|\xrightarrow{\mathbb P}0.
\end{equation*}
Combining this with the convergence of $\mathcal U_{j,L_j}^r/r$ to $\mathcal X_j^*$ in \eqref{eq:prelimit-U-convergence}, we obtain
\begin{equation*}
\sup_{0\leq\xi\leq T}\left|\frac1r p_{j,L_j}^r(\lfloor r\xi\rfloor)-p_{j,L_j}^{\pi^*}(\xi)\right|\leq \sup_{0\leq\xi\leq T}\left|\bar A_j^r(\xi)-A_j(f^*(\xi))\right|+\sup_{0\leq\xi\leq T}\left|\frac1r\mathcal U_{j,L_j}^r(\lfloor r\xi\rfloor)-\mathcal X_j^*(\xi)\right|\xrightarrow{\mathbb P}0.
\end{equation*}
This proves the convergence of all prefill and decode states in \eqref{eq:prelimit-p-convergence} and \eqref{eq:prelimit-d-convergence}.
\end{proof}
\section{Proof of Corollary~\ref{cor:moving-average-cost-convergence}}
\label{sec:proof-moving-average-cost-convergence}

\begin{proof}
First we prove that the prelimit policy completes all requests in finitely many steps. Define
\begin{equation*}
t_H:=\min\left\{t\in\mathbb Z_+:f^r(t)\ge H\right\},\qquad t_{\mathrm c}^r:=\inf\left\{t\in\mathbb Z_+:t\ge t_H,\ z_t^r=0\right\}.
\end{equation*}
We bound $t_{\mathrm c}^r$ by adding the four step bounds in Figure~\ref{fig:cost-completion-timeline}.

\begin{figure}[H]
\centering
\begin{tikzpicture}[x=1cm,y=1cm,font=\footnotesize,>=Latex]
\draw[->,line width=0.8pt,black!70] (0,0) -- (13.5,0);
\foreach \x in {0,3.3,6.6,9.9,13.2} {\fill[black!70] (\x,0) circle (1.6pt);}
\node[above=7pt] at (1.65,0) {$t_H\le\lceil H/b^r\rceil$};
\node[above=7pt] at (4.95,0) {$n$};
\node[above=7pt] at (8.25,0) {$w_r-1$};
\node[above=7pt] at (11.55,0) {$H_{\mathrm{kv}}+1$};
\node[below=7pt] at (0,0) {Start};
\node[below=7pt] at (3.3,0) {All arrivals available};
\node[below=7pt] at (6.6,0) {All plans assigned};
\node[below=7pt] at (9.9,0) {All requests admitted};
\node[below=7pt] at (13.2,0) {All requests completed};
\end{tikzpicture}
\caption{Four step bounds for request completion.}
\label{fig:cost-completion-timeline}
\end{figure}
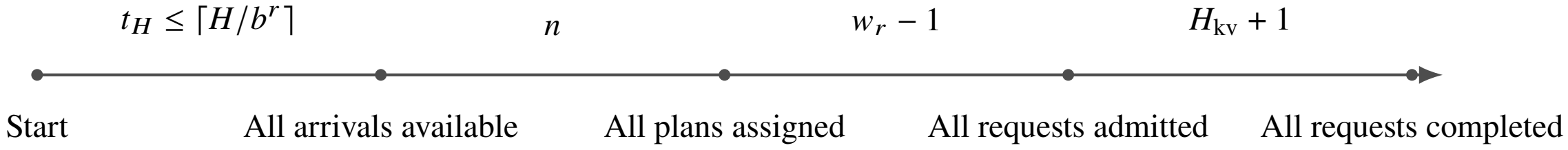
\vspace{-0.4in}
Since $\Delta^r(t)\ge b^r$,   we have
$
f^r(t)=\sum_{u<t}\Delta^r(u)\ge tb^r, t_H\le\left\lceil\frac{H}{b^r}\right\rceil.
$
For $t\ge t_H$, we have $A_j^r(t)=A_j^r(t_H)$. Applying the planning rule \eqref{eq:planned-admission-policy} recursively over the classes gives
\begin{equation}
\sum_{j=1}^J\tilde m_jy_j^r(t)=\min\left\{\sum_{j=1}^J\tilde m_j\left(A_j^r(t_H)-\sum_{u<t}y_j^r(u)\right),\frac{w_rM}{w_r+H_{\mathrm{kv}}}\right\},\qquad t\ge t_H.
\label{eq:cost-total-planning}
\end{equation}
Thus, after $n$ further planning steps,
\begin{equation}
\sum_{j=1}^J\tilde m_j\left(A_j^r(t_H)-\sum_{u<t_H+n}y_j^r(u)\right)=\left[\sum_{j=1}^J\tilde m_j\left(A_j^r(t_H)-\sum_{u<t_H}y_j^r(u)\right)-n\frac{w_rM}{w_r+H_{\mathrm{kv}}}\right]^+.
\label{eq:cost-unplanned-work}
\end{equation}
Consequently, all remaining requests receive plans within
\begin{equation}
n=\left\lceil\frac{\sum_{j=1}^J\tilde m_j\left(A_j^r(t_H)-\sum_{u<t_H}y_j^r(u)\right)}{w_rM/(w_r+H_{\mathrm{kv}})}\right\rceil\le\left\lceil\frac2M\sum_{j=1}^J\tilde m_jA_j^r(t_H)\right\rceil \text{ steps.}
\label{eq:cost-planning-steps}
\end{equation}

By the moving-average rule \eqref{eq:smoothed-admission}, a plan $y_j^r(s)$ contributes $y_j^r(s)/w_r$ to admissions in each of steps $s,\ldots,s+w_r-1$. The stage shifts \eqref{eq:prelimit-U-shift}--\eqref{eq:prelimit-V-shift} give
\begin{equation*}
\Delta\mathcal U_{j,L_j}^r(u)=
\begin{cases}
\Delta\mathcal U_{j,1}^r(u+L_j-1),&o_j=1,\\
\Delta\mathcal V_{j,1}^r(u+L_j+o_j-2),&o_j\ge2.
\end{cases}
\end{equation*}
Thus, requests admitted at step $u$ complete at the end of step $u+L_j+o_j-2$. Since $L_j+o_j-2\le H_{\mathrm{kv}}$, the four bounds yield
\begin{equation}
t_{\mathrm c}^r\le t_H+n+(w_r-1)+(H_{\mathrm{kv}}+1)\le\left\lceil\frac{H}{b^r}\right\rceil+\left\lceil\frac2M\sum_{j=1}^J\tilde m_jA_j^r(t_H)\right\rceil+w_r+H_{\mathrm{kv}}.
\label{eq:cost-completion-step-bound}
\end{equation}
Since no requests arrive after $H$ and $f^r(t_H)\ge H$, \eqref{eq:A-def} gives
\begin{equation}
A_j^r(t_H)=p_{j,L_j}^r(0)+N_j^r\!\left(r\int_0^{f^r(t_H)}\lambda_j(s)\,ds\right)=p_{j,L_j}^r(0)+N_j^r\!\left(r\int_0^\infty\lambda_j(s)\,ds\right).
\label{eq:cost-total-arrivals}
\end{equation}
Dividing \eqref{eq:cost-completion-step-bound} by $r$ and using $\lceil x\rceil\le x+1$, we obtain
\begingroup
\thinmuskip=2mu \medmuskip=2mu \thickmuskip=3mu
\begin{equation}
\frac{t_{\mathrm c}^r}r\le\frac{H}{rb^r}+\frac{w_r+H_{\mathrm{kv}}+2}{r}+\frac{2}{Mr}\sum_{j=1}^J\tilde m_jA_j^r(t_H)\lesssim1+\frac1r\sum_{j=1}^J\left[p_{j,L_j}^r(0)+N_j^r\!\left(r\int_0^\infty\lambda_j(s)\,ds\right)\right]<\infty,
\label{eq:cost-scaled-completion-bound}
\end{equation}
\endgroup
where the constant implicit in $\lesssim$ is independent of $r$ and the arrival realization. The initial-condition convergence and the Poisson law of large numbers then imply that, for a sufficiently large fixed $T<\infty$, $\mathbb P\!\left(t_{\mathrm c}^r\le\lfloor rT\rfloor\right)\longrightarrow1$.

We next prove cost convergence. For sufficiently large $T$, $p_{j,L_j}^{\pi^*}(\xi)=0$ for all $\xi\ge T$. Theorem~\ref{thm:moving-average-fluid-limit}, $\mathbb P\!\left(t_{\mathrm c}^r\le\lfloor rT\rfloor\right)\longrightarrow1$, and continuity of Stieltjes integration give
\begin{equation}
\frac1r\sum_{t=0}^{\infty}\Delta^r(t)\sum_{j=1}^Jc_j\left(\sum_{k=1}^{L_j}p_{j,k}^r(t)+\sum_{k=1}^{o_j-1}d_{j,k}^r(t)\right)\xrightarrow{\mathbb P}\int_0^T\sum_{j=1}^Jc_jp_{j,L_j}^{\pi^*}(\xi)\,df^*(\xi)=\mathcal J(\pi^*;z_0).
\label{eq:cost-sample-path-convergence}
\end{equation}

We now prove convergence of the expected costs. The memory constraint and \eqref{eq:processing_time_limit} give
\begin{equation}
\Delta^r(t)\le(a_{\mathrm{pf}}^r+a_{\mathrm{kv}}^r)M+b^r,\qquad\sup_{r\ge1,\,t\ge0}r\Delta^r(t)<\infty.
\label{eq:cost-step-time-bound}
\end{equation}
At every step, the total number of unfinished requests is at most $\sum_jA_j^r(t_H)$. All requests have completed by real time $f^r(t_{\mathrm c}^r)$, and \eqref{eq:cost-step-time-bound} gives $f^r(t_{\mathrm c}^r)\le t_{\mathrm c}^r\sup_t\Delta^r(t)\lesssim t_{\mathrm c}^r/r$. Combining this with \eqref{eq:cost-scaled-completion-bound} gives
\begin{equation}
\begin{aligned}
&\frac1r\sum_{t=0}^{\infty}\Delta^r(t)\sum_jc_j\left(\sum_{k=1}^{L_j}p_{j,k}^r(t)+\sum_{k=1}^{o_j-1}d_{j,k}^r(t)\right)\\
&\qquad\le(\max_jc_j)\frac{\sum_jA_j^r(t_H)}r\,f^r(t_{\mathrm c}^r)\lesssim\frac{\sum_jA_j^r(t_H)}r\,\frac{t_{\mathrm c}^r}r\lesssim\left(1+\frac1r\sum_jA_j^r(t_H)\right)^2.
\end{aligned}
\label{eq:cost-pathwise-bound}
\end{equation}
Using $p_{j,L_j}^r(0)/r\to p_{j,L_j}(0)$, \eqref{eq:cost-total-arrivals}, \eqref{eq:cost-pathwise-bound}, and the Poisson fourth-moment formula, we obtain
\begin{equation*}
\begin{aligned}
&\sup_{r\ge1}\mathbb E\!\left[\left\{\frac1r\sum_{t=0}^{\infty}\Delta^r(t)\sum_jc_j\left(\sum_{k=1}^{L_j}p_{j,k}^r(t)+\sum_{k=1}^{o_j-1}d_{j,k}^r(t)\right)\right\}^{\!2}\right]\lesssim\sup_{r\ge1}\mathbb E\!\left[\left(1+\frac1r\sum_jA_j^r(t_H)\right)^4\right]<\infty.
\end{aligned}
\end{equation*}
Thus these costs are uniformly integrable. Together with \eqref{eq:cost-sample-path-convergence}, this yields
\begin{equation*}
\mathcal J^r(\widehat\pi^r;z_0^r)=\mathbb E\!\left[\frac1r\sum_{t=0}^{\infty}\Delta^r(t)\sum_jc_j\left(\sum_{k=1}^{L_j}p_{j,k}^r(t)+\sum_{k=1}^{o_j-1}d_{j,k}^r(t)\right)\right]\longrightarrow\mathcal J(\pi^*;z_0).
\end{equation*}
\end{proof}

\section{Numerical Implementation Details}
\label{appendix:implentation}
We set vLLM's KV-cache blocks to 64 tokens and cap prefill chunks at 2,048 tokens. To reduce preemptions, vLLM admits a request only if its full prompt fits in the KV cache.

To construct the output-length predictor, we randomly split the requests into training, validation, and test sets, with proportions of 40\%, 20\%, and 40\%, respectively.
All requests from the same conversation are assigned to the same set.
We use a multilayer perceptron (MLP) to estimate the first and second moments of output length conditional on the request's input features.
We map each of the first 32 input-block identifiers to an 8-dimensional embedding and concatenate them into a 256-dimensional vector. We log-transform and standardize the prompt length and conversation turn number, and one-hot encode the dataset identifier. The resulting 262-dimensional feature vector serves as input to an MLP with two hidden layers of 64 units each and ReLU activations, followed by a linear output layer with two units.

The two output units are trained to predict standardized \(O_i\) and \(O_i^2\), respectively, with each target standardized using its mean and standard deviation computed separately for each source dataset. When making predictions, we reverse the standardization and clip the estimates using the admissible output-length range and the requirement of nonnegative variance.

We select the number of training epochs based on the validation error, then retrain the model from scratch on the combined training and validation sets for the selected number of epochs.

\end{document}